\documentclass[10pt,twoside]{amsart}
\usepackage{latexsym,amssymb,graphicx,fdsymbol,enumitem}
\usepackage{mathrsfs,}
\usepackage[all]{xy}

\newdir{ >}{{}*!/-10pt/\dir{>}}

\def\smono[#1]{\ar@{->}[#1]|@{|}}
\newcommand{\xymono}{\ar@{ >->}}
\newcommand{\xyepi}{\ar@{->>}}

\let\oldtocsection=\tocsection

\let\oldtocsubsection=\tocsubsection

\let\oldtocsubsubsection=\tocsubsubsection

\renewcommand{\tocsection}[2]{\hspace{0em}\oldtocsection{#1}{#2}}
\renewcommand{\tocsubsection}[2]{\hspace{1em}\oldtocsubsection{#1}{#2}}
\renewcommand{\tocsubsubsection}[2]{\hspace{2em}\oldtocsubsubsection{#1}{#2}}

\newtheorem{theorem}{Theorem}[section]
\newtheorem{lemma}[theorem]{Lemma}
\newtheorem{corollary}[theorem]{Corollary}
\newtheorem{proposition}[theorem]{Proposition}

\theoremstyle{definition}

\newtheorem{definition}[theorem]{Definition}
\newtheorem{example}[theorem]{Example}
\newtheorem{remark}[theorem]{Remark}

\numberwithin{equation}{section}

\newcommand{\Z}{\mathbb{Z}}

\newcommand{\A}{\mathcal{A}}
\newcommand{\B}{\mathcal{B}}
\newcommand{\C}{\mathcal{C}}

\newcommand{\E}{\mathscr{E}}

\renewcommand{\H}{\mathscr{H}}
\newcommand{\M}{\mathcal{M}}
\newcommand{\N}{\mathbb{N}}

\newcommand{\ffi}{\varphi}
\newcommand{\eps}{\varepsilon}

\newcommand{\colim}{\operatorname{colim}}

\newcommand{\im}{\operatorname{Im}}

\newcommand{\Hom}{\operatorname{Hom}}

\newcommand{\Fun}{\operatorname{Fun}}

\newcommand{\Ch}{\operatorname{Ch}}
\newcommand{\can}{\operatorname{can}}
\newcommand{\Cone}{\operatorname{Cone}}

\newcommand{\Ar}{\operatorname{Ar}}
\newcommand{\RR}{\mathcal{R}}

\renewcommand{\P}{\mathcal{P}}

\newcommand{\coker}{\operatorname{coker}}

\newcommand{\quis}{\operatorname{quis}}

\newcommand{\U}{\mathcal{U}}

\newcommand{\ev}{\operatorname{ev}}

\newcommand{\Ab}{\operatorname{Ab}}

\newcommand{\Quad}{\operatorname{Quad}}

\newcommand{\cone}{\operatorname{cone}}
\newcommand{\Mor}{\operatorname{Mor}}

\newcommand{\proj}{\operatorname{proj}}

\renewcommand{\SS}{\mathcal{S}}

\newcommand{\n}{\mathcal{N}}
\newcommand{\FormCatW}{\operatorname{FormCatW}}
\newcommand{\lax}{\operatorname{lax}}
\newcommand{\str}{\operatorname{str}}

\newcommand{\Ac}{\operatorname{Ac}}
\newcommand{\Q}{\mathcal{Q}}
\newcommand{\QF}{\mathcal{Q}\operatorname{Form}}
\newcommand{\scN}{\mathcal{N}}

\DeclareRobustCommand\longtwoheadrightarrow
 {\xymatrix{\xyepi[r]&}}
\DeclareRobustCommand\longrightarrowtail
     {\xymatrix{\xymono[r]&}}
\DeclareRobustCommand\longtwoheadleftarrow
 {\xymatrix{&\xyepi[l]}}

\newcommand{\downsim}{\hspace{.4ex} \downarrow \hspace{-.8ex} \wr \hspace{.8ex}}

\title[Higher $K$-theory of forms II]{Higher $K$-theory of forms II.\\ From exact categories to chain complexes}
 \author{Marco Schlichting}
 \address{Marco Schlichting, Mathematics Institute,
Zeeman Building,
University of Warwick,
Coventry CV4 7AL, UK} 

\thanks{}

\email{m.schlichting@warwick.ac.uk}

\subjclass{}

\keywords{}

\begin{document}
\bibliographystyle{alpha}

\begin{abstract}
We prove basic statements about the Hermitian $K$-theory of exact form categories with weak equivalences providing the background for the comparison in \cite{ourForm3} of the classical $1$-categorical version of the Hermitian $K$-theory of exact categories as defined in \cite{myForm1}  and the $\infty$-categorical version of \cite{9authors}.
Notably, in the current paper, we extend a quadratic functor with values in abelian groups from an exact category to its category of bounded chain complexes in a way that does not change Grothendieck-Witt spaces.
 \end{abstract}

\maketitle

\tableofcontents

\section{Introduction}

This is the second paper in the series \cite{myForm1}, \cite{myForm2}, \cite{ourForm3} announced in \cite{myForm1} which lead to \cite{myKSPZ}, \cite{myCRAScorrect}. 
In \cite{myForm1} we defined the Hermitian $K$-theory $GW(\E,Q)$ of an exact form category with strong duality $(\E,\sharp,\can,Q)$ and showed that if all admissible exact sequences in $\E$ split, then that definition agrees with the group completion definition of quadratic spaces in $(\E,\sharp,\can,Q)$ \cite[Theorem 1.1]{myForm1}. 
This is the analog in Hermitian $K$-theory of Quillen's $Q=+$ Theorem \cite[Theorem p.7]{grayson:HigherKII}.

In this paper, we will study the Hermitian $K$-theory of exact form categories with weak equivalences and prove basic results such as Additivity (Theorem \ref{thm:AddtyForWeakExFunctors}), a Fibration Theorem (Theorem \ref{thm:ChgOfWkEq}) and a Cofinality Theorem (Theorem \ref{thm:CofinilityForWeakEq}) the analogs of which are well-known in algebraic $K$-theory.
Our main result (Theorem \ref{prop:gilletWald}) associates to any exact form category with strong duality $(\E,\sharp,\can,Q)$ in the sense of \cite{myForm1}, an exact form category with weak equivalences  $(\Ch_b\E,\quis,\sharp, \can,Q)$ of bounded chain complexes in $\E$ and states that the inclusion $\E \subset \Ch_b\E$ as complexes concentrated in degree $0$ induces a homotopy equivalence of associated Grothendieck-Witt spaces
\begin{equation}
\label{eqn:introMainResult}
GW(\E,Q) \stackrel{\sim}{\longrightarrow} GW(\Ch_b\E,\quis,Q).
\end{equation}
This generalises \cite[Proposition 6]{myMV} from symmetric forms to arbitrary quadratic forms with values in abelian groups.
 Our main contribution here is our definition of the extension of the quadratic functor $Q:\E^{op} \to \Ab$ defined on an exact category $\E$ to a quadratic functor $Q:(\Ch_b\E)^{op} \to \Ab$ defined on all bounded chain complexes in $\E$ such that the homotopy equivalence (\ref{eqn:introMainResult}) holds.

We note that results similar to our Additivity, Fibration and Cofinality Theorems have been obtained in  \cite{9authors} and its sequals in the setting of infinity categories. 
To date, however, agreement of the approach of \cite{9authors} with the classical $1$-categorical approach of \cite{myForm1} has only been proved in special cases such as the split exact case \cite{HebestreitSteimle}, \cite{myForm1}, the case when $2$ is invertible \cite[Corollary B.2.5]{9authors1/2Agreement}, or the case when Zariski-descent holds \cite{CalmesA1rep}, \cite{myMV}.
In \cite{ourForm3}, 
we will associate to any complicial exact form category with weak equivalences a Poincare infinity category in the sense of \cite{9authors} and show that the 1-categorical Hermitian $K$-theory of the former as defined in \cite{myForm1} and this paper agrees with the infinity categorical Hermitian $K$-theory of the latter as defined in \cite{9authors}.
Together with the homotopy equivalence (\ref{eqn:introMainResult}) of Theorem \ref{prop:gilletWald}, this shows that the approach of \cite{myForm1} agrees with that of \cite{9authors} generalising all previous agreement statements cited above.

In addition to the results already mentioned, we will construct a sequence of (usually non-connective) spectra $GW^{[n]}(\E,w,Q)$, $n\geq 0$, for any exact form category with weak equivalences and strong symmetric cone $(\E,w,\sharp,\can,Q)$ whose negative homotopy groups $GW_i^{[n]}(\E,w,Q) = \pi_iGW^{[n]}(\E,w,Q)$, $i<0$, are the Witt groups $W(\E^{[n-i]},w,Q)$ of certain exact form categories with weak equivalences $(\E^{[n]},w,\sharp,\can,Q)$, where $(\E^{[0]},w,\sharp,\can,Q)=(\E,w,\sharp,\can,Q)$.
The infinite loop space $\Omega^{\infty}GW^{[n]}(\E,w,Q)$ is the Grothendieck-Witt space $GW(\E^{[n]},w,Q)$, and the sequence of spectra fit into a homotopy fibration of spectra with connective algebraic $K$-theory for $n\geq 0$
$$GW^{[n]}(\E,w,Q) \stackrel{F}{\longrightarrow} K(\E,w)  \stackrel{H}{\longrightarrow} GW^{[n+1]}(\E,w,Q)$$
generalising the Bott sequence of \cite[Theorem 6.1]{myJPAA}.
This is proved in Theorem \ref{thm:BottSeq}.

Most of the results in this paper were obtained in conjuction with \cite{myForm1}.
Often, the proofs are direct adaptations of the case of symmetric forms and are {\em mutatis mutandis} the same as in \cite{myHermKex}, \cite{myMV} and \cite{myJPAA}. 
Our main contribution is to equip everything in sight with the correct quadratic functor. 
At places, we have taken the opportunity to simplify proofs ({\em e.g.}, Theorem \ref{thm:exAddty}), or to shorten proofs relying on the symmetric forms case ({\em e.g.}, Theorem \ref{thm:filt_locn} and Lemma \ref{lem:filt18v2}).

\newpage

\part{Exact form categories}
\label{part:ExFormCat}

\section{The form $Q$-construction}
\label{sec:FormQ}

Recall \cite[Definition 2.1]{myHermKex} that
an {\em exact category with duality} is a tuple $(\E,\sharp,\can)$ where $\E$ is an exact category \cite{quillen:higherI}, $\sharp:\E^{op} \to \E$ is an exact functor, $\can:1 \to \sharp\sharp$ is a natural transformation such that $1_{X^{\sharp}} = \can_X^{\sharp}\can_{X^{\sharp}}$ for all objects $X \in \E$.
We denote by $\sigma$ the map 
$$\sigma:\E(X,Y^{\sharp}) \to \E(Y,X^{\sharp}): f \mapsto f^{\sharp}\can_Y$$
natural in $X$ and $Y$.
Occasionally, we may write $\bar{f}$ for $\sigma(f)$.
Note that $\sigma^2=1$.
We call the duality {\em strong} if $\can:1 \to \sharp\sharp$ is an isomorphism, and we call the duality {\em strict} if $X^{\sharp\sharp}=X$ and $\can_X=1_X$ for all objects $X$ of $\E$.
Clearly, a strict duality is also strong.

\begin{definition}[{\cite[Definitions 2.1 and 2.22]{myForm1}}\footnote{In \cite{myForm1} the functor $\sharp$ was always (implicitly) assumed to be exact though this wasn't stated explicitly in \cite[Definition 2.22]{myForm1}.}]
\label{dfn:ExFormCat}
An {\em exact form category} is a tuple $(\E,\sharp,\can,Q)$ where $(\E,\sharp,\can)$ is an exact category with duality and $Q:\E^{op} \to \Ab$ is a functor with values in abelian groups together with a functorial $C_2$-equivariant diagram of abelian groups
$$(\E(X,X^{\sharp}),\sigma) \stackrel{\tau}{\longrightarrow} Q(X) \stackrel{\rho}{\longrightarrow} (\E(X,X^{\sharp}),\sigma)$$
where $Q(X)$ carries the trivial $C_2$-action.
We require the following.
\begin{enumerate}
\item 
\label{item1:dfn:ExFormCat}
$\rho\tau=1+\sigma$.
\item
\label{item2:dfn:ExFormCat}
For all $f,g\in \E(X,Y)$ and $\xi\in Q(Y)$ we have 
$$(f+ g)^{\bullet}(\xi) =  f^{\bullet}(\xi) +g^{\bullet}(\xi)  + \tau(g^{\sharp} \circ \rho(\xi) \circ f)$$
where $f^{\bullet}(\xi)$ denotes $Q(f)(\xi)$.
\item
\label{item3:dfn:ExFormCat}
For every admissible exact sequence $X \stackrel{i}{\rightarrowtail} Y \stackrel{p}{\twoheadrightarrow} Z$ in $\E$, the following sequence of abelian groups is exact
$$\xymatrix{0 \ar[r] & Q(Z) \ar[r]^{p^{\bullet}} &  Q(Y) \ar[rr]^{\hspace{-6ex}(i^{\sharp}\circ \rho(\underline{\phantom{x}}),\ i^{\bullet})} && \A(Y,X^{\sharp}) \times Q(X)}.$$
\end{enumerate}
The elements of $Q(X)$ are called quadratic forms on $X$.
\end{definition}

\begin{definition}
Let $(\E,\sharp,\can)$ be an exact form category with strong duality.
A {\em quadratic space}\footnote{The notion of quadratic space will be generalised in Definition \ref{dfn:QspaceWeEq} below.}
in $(\E,Q)$ is a pair $(X,\xi)$ where $X$ is an object of $\E$ and $\xi\in Q(X)$ is a quadratic form on $X$ such that $\rho(\xi)$ is an isomorphism.
A map of quadratic spaces $f:(X,\xi) \to (Y,\zeta)$ is a map $f:X \to Y$ in $\E$ such that $\xi=f^{\bullet}(\zeta)$.
Composition of maps of quadratic spaces is composition in $\E$.
An isomorphism of quadratic spaces is also called {\em isometry}.
We denote by 
$$i\Quad(\E,Q)$$
the category of quadratic spaces and isometries in $(\E,Q)$.
\end{definition}

Recall \cite{quillen:higherI} that Quillen associates to any exact category $\E$ a new category $\Q\E$ and defines the $K$-theory space of $\E$ as $K(\E) = \Omega |\Q\E|$, the loop space of the classifying space $|\Q\E|$ of $\Q\E$.
The objects of $\Q\E$ are the objects in $\E$ and the arrows $X \to Y$ are equivalences classes  $[X \stackrel{p}{\twoheadleftarrow} U \stackrel{s}{\rightarrowtail} Y] $ of data where $p$ is an admissible epimorphism and $s$ is an admissible monomorphism in $\E$.
We have equality of equivalence classes $[X \stackrel{p}{\twoheadleftarrow} U \stackrel{s}{\rightarrowtail} Y] = [X \stackrel{p'}{\twoheadleftarrow} U' \stackrel{s'}{\rightarrowtail} Y] $ if and only if there is an isomorphism $u:U \cong U'$ such that $p'=up$ and $s'=us$.
Composition is defined by
$$[Y \stackrel{q}{\twoheadleftarrow} V \stackrel{r}{\rightarrowtail} Z] \circ [X \stackrel{p}{\twoheadleftarrow} U \stackrel{s}{\rightarrowtail} Y] = [X \stackrel{pq'}{\twoheadleftarrow} U\times_YV \stackrel{rs'}{\rightarrowtail} Z]$$
where $q':U\times_YV \twoheadrightarrow U$ and $s':U\times_YV \rightarrowtail V$ are the canonical projections.
Note that $[X \stackrel{p}{\twoheadleftarrow} U \stackrel{s}{\rightarrowtail} Y] = s_!\circ p^!$ where $p^!=[X \stackrel{p}{\twoheadleftarrow} U \stackrel{1}{\rightarrowtail} U]$ and $s_!=[U \stackrel{1}{\twoheadleftarrow} U \stackrel{s}{\rightarrowtail} Y]$. 
We have $s_!\circ p^!=q^!\circ t_!$ where $V$ is the pushout of $p$ along $s$ with $q:Y \twoheadrightarrow V$ and $t:X \rightarrowtail V$ the induced maps.

\begin{definition}[The form $Q$-construction]
\label{dfn:FormQ}
Let $(\E,\sharp,\can,Q)$ be an exact form category with strong duality.
We define a new category 
$$\QF(\E,Q) = \QF(\E,\sharp,\can,Q),$$
called {\em form $Q$-construction}, as follows.
Objects are the quadratic spaces in $(\E,Q)$.
An arrow $(X,\xi) \to (Y,\zeta)$ is an arrow $[X \stackrel{p}{\twoheadleftarrow} U \stackrel{s}{\rightarrowtail} Y]$ in Quillen's $Q\E$ such that 
\begin{enumerate}
\item
$p^{\bullet}(\xi) = s^{\bullet}(\zeta)$ and
\item
\label{dfn:FormQ:item2}
 $\ker(p) = U^{\perp}$ where $U^{\perp} = \ker(s^{\sharp}\rho(\zeta))$ is the orthogonal of $U$ in $Y$.
 \end{enumerate}
Composition is composition in Quillen's $\Q\E$.
Note that (\ref{dfn:FormQ:item2}) is equivalent to the square
$$\xymatrix{
U\  \ar@{>->}[r]^s \ar@{->>}[d]_p & Y\ar@{->>}[d]^{s^{\sharp}\rho(\zeta)}\\
X \  \ar@{>->}[r]_{p^{\sharp}\rho(\xi)} & U^{\sharp}}
$$
being bicartesian.
\end{definition}

In \cite[Definition 6.1]{myForm1}, for every poset $\P$ with strict duality $\P^{op} \to \P: i\mapsto i'$ and every exact form category $(\E,\sharp,\can,Q)$, we have equipped the category of functors $\Fun(\P,\E)$ with the structure $(\sharp,\can,Q)$ of an exact form category.
The dual of a functor $A:\P \to \E$ is $(A^{\sharp})_i=(A_{i'})^{\sharp}$.
A {\em quadratic form} on a functor $A:\P \to \A:i\mapsto A_i$ is a pair $(\xi,\ffi)$ where  $\xi = (\xi_i)_{i \leq i'}$ is a family  of quadratic forms $\xi_i \in Q(A_i)$ indexed over those $i\in \P$ satisfying $i\leq i'$, and $\ffi:A \to A^{\sharp}$ is a symmetric form in $\Fun(\P,\A)$ such that $A_{i \leq j}^{\bullet}(\xi_j)= \xi_i$ whenever $i\leq j \leq j'\leq i'$, and  $\rho(\xi_i) = \ffi_{i'}\circ A_{i\leq i'}$.
In particular, the categories $S_n\E$ of Waldhausen \cite{wald:spaces} are equipped with the structure $(\sharp,\can,Q)$ of an exact form category as a fully exact form subcategory of 
$\Fun(\Ar([n]),\E)$ where $[n]=\{0<1<\dots n\}$ has unique strict duality $i \mapsto n-i$; see \cite[\S 6]{myForm1}.

\begin{example}
\label{ex:S2explicit}
Let us make the structure of the exact form category on $S_2\E$ more explicit as this will be needed in the Additivity Theorem \ref{thm:exAddty} below.
Similar observations can be made for $S_n\E$.
Let $(\E,\ast,\can,Q)$ be an exact form category.
Objects of $S_n\E$ are functors $A:\Ar[n] \to \E: (i\leq j) \mapsto A_{ij} $ such that $A_{ii}=0$ and $A_{ij} \rightarrowtail A_{ik} \twoheadrightarrow A_{jk}$ is admissible exact whenever $i\leq j\leq k$.

The exact category $S_2\E$ is therefore the category of admissible exact sequences $A_{01} \rightarrowtail A_{02} \twoheadrightarrow A_{12}$ in $\E$ with point-wise exact structure.
In $\Ar[2]$, only the pairs $(i,j)$ of the form $(0,0)$, $(0,1)$, $(0,2)$, $(1,1)$ satisfy $(i,j) \leq (i,j)' =(2-j,2-i)$.
Thus, to give a quadratic form on $A\in S_2\E$ is the same as to give a symmetric bilinear form $\ffi:A \to A^*$ and quadratic forms $\xi_{0,0} \in Q(A_{0,0})$, $\xi_{0,1} \in Q(A_{0,1})$, $\xi_{1,1} \in Q(A_{1,1})$, $\xi_{1,2} \in Q(A_{1,2})$ satisfying the compatibilities of \cite[Definition 6.1]{myForm1} recalled above.
Since $A_{0,0}=A_{1,1}=0$, we have $\xi_{0,0}=\xi_{1,1}=0$. 
Since $(0,1) \leq (1,1),(0,2)$, the form $\xi_{0,1}$ is the restriction of $\xi_{0,2}$ and of $\xi_{1,1}=0$.
In particular, $0=\xi_{0,1}=(\xi_{0,2})_{|A_{0,1}}$.
Moreover, $\rho(\xi_{0,2})=\ffi_{0,2}$.
Note that $\ffi_{1,2}=\ffi_{0,1}\circ \can_{1,2}$, that $\ffi_{0,1}$ is determined by $\ffi_{0,2}$ and that $\ffi_{0,1}$ exists if and only if 
$A_{(0,1)\to (0,2)}^*\circ  \ffi_{0,2}\circ A_{(0,1) \to (0,2)} =0$.

Thus, $S_2\E$ is the category of admissible exact sequences $X \stackrel{a}{\rightarrowtail} Y \twoheadrightarrow Z$ in $\E$ with point-wise exact structure.
A quadratic form on such an object is a quadratic form $\xi\in Q(Y)$ such that $a^{\bullet}\xi=0$.
Its associated symmetric form is determined by $\rho(\xi):Y \to Y^*$. 
Indeed, it satisfies $a^*\circ \rho(\xi)\circ a = \rho(a^{\bullet}\xi)=0$ and thus defines a symmetric form on the exact sequence.
If $\E$ has a strong duality, then the form $\xi$ on $X \stackrel{a}{\rightarrowtail} Y \twoheadrightarrow Z$ is non-degenerate if and only if the map of exact sequences induced by $\rho(\xi):Y \to Y^*$ is an isomorphism.
The latter is equivalent to $\rho(\xi):Y \to Y^*$ being an isomorphism such that 
$(a, a^*\rho(\xi))$ are the maps in an admissible exact sequence $X \rightarrowtail Y \twoheadrightarrow X^*$.
Note that, up to equivalence of form categories, we do not need to specify the quotient $Z$ is it is determined up to isomorphism by the admissible monomorphism $a:X \rightarrowtail Y$.
\end{example}

For a simplicial object $X$, denote by $X^e$ the edge-wise subdivision \cite[\S 1.9]{wald:spaces}, \cite[\S 6]{myForm1}.
Its object of $n$-simplicies is $X^e[n] = X([n]^{op}[n])$ where $[n]^{op}[n]$ is the totally ordered set $\{n^{op} < \dots < 1^{op} <0^{op} < 0 < \dots n\}$.
The natural inclusion $[n] \to [n]^{op}[n]: i \mapsto i$ induces the natural transformation of simplicial objects $X^e \to X$ which induces a weak equivalence $|X^e| \to |X|$ of topological realisations for all simplicial sets $X$ \cite[Lemma 1]{myMV}.

Recall \cite[Definition 6.3]{myForm1} that the Grothendieck-Witt space of an exact form category with strong duality $(\E,\sharp,\can,Q)$ is defined as the space $GW(\E,Q)$ which sits in a homotopy fibration
\begin{equation}
\label{GWFibDfnSe}
GW(\E,Q) \longrightarrow | i\Quad(S_{\bullet}^e\E,Q)| \longrightarrow |iS_{\bullet}\E|
\end{equation}
where the second map is the composition of the map $i\Quad(S_{\bullet}^e\E,Q) \to iS^e_{\bullet}\E: (A,\xi) \mapsto A$ forgetting quadratic forms followed by the natural transformation $iS^e_{\bullet}\E \to iS_{\bullet}\E$.
Since $S_0\E=0$ and $S_1\E=\E$, the natural zero simplex inclusion $X^e[0]=X_1 \to X^e$ 
induces, for $X= i\Quad(S_{\bullet}^e\E,Q)$, a map $$| i\Quad(\E,Q)| \to GW(\E,Q).$$
If all admissible exact sequences in $\E$ split then this map is a group completion of the homotopy commutative $H$-space $(|i\Quad(\E,Q)|,\perp)$; see  \cite[Theorem 6.6]{myForm1}.

Below, we will need to iterate the $S^e_{\bullet}$-construction. 
To avoid awkward notation, we may write $\RR_n$ for $S_n^e$.
So
$$\RR_n\E = S^e_n\E.$$
In this notation, the homotopy fibration (\ref{GWFibDfnSe})
becomes the homotopy fibration
\begin{equation}
\label{GWFibDfnR}
GW(\E,Q) \longrightarrow | i\Quad(\RR_{\bullet}\E,Q)| \longrightarrow |iS_{\bullet}\E|.
\end{equation}

\begin{proposition}
\label{prop:SeAndQ}
For any exact form category with strong duality $(\E,\sharp,\can,Q)$, there is natural zigzag of
homotopy 
equivalences between the following two topological spaces 
$$|\QF(\E,Q)| \sim |i\Quad(\RR_{\bullet}\E,Q)|.$$
\end{proposition}

\begin{proof}
The proof is {\em mutatis mutandis} the same as in \cite[Proposition 2]{myMV}, and we omit the details.
\end{proof}

\begin{remark}
As in  \cite[Proposition 2]{myMV}, Proposition \ref{prop:SeAndQ} implies that there is a homotopy fibration 
$$GW(\E,Q) \to \QF(\E,Q) \to \Q\E,$$
but we won't need this.
\end{remark}

\begin{definition}[The form category of exact functors]
\label{dfn:FunExFomCat}
Let $(\A,*,\can_{\A}, Q_{\A})$ and $(\B,*,\can_{\B}, Q_{\B})$ be two small exact form categories.
We make the category of exact functors $\Fun_{ex}(\A,\B)$
into an exact form category 
\begin{equation}
\label{eqn:FunExFomCat}
(\Fun_{ex}(\A,\B),\sharp,\eta,Q)
\end{equation}
as follows.
For exact functors $F,G:\A \to \B$, denote by $[F,G]$ the set of maps from $F$ to $G$ in the category $\Fun_{ex}(\A,\B)$, that is, the set of natural transformations $F \to G$.
\begin{enumerate}[labelindent=0pt,leftmargin=*]
\item
\label{dfn:FunEx:item1}
A sequence $F_0 \to F_1 \to F_2$ of functors $\A \to \B$ is called admissible exact if for all objects $A\in \A$ the sequence $F_0A \rightarrowtail F_1A \twoheadrightarrow F_2A$ is admissible exact in $\B$.
\item
\label{dfn:FunEx:item2}
The dual of a functor $F:\A \to \B$ is $F^{\sharp}=\ast F \ast$.
The dual of a natural transformation $f:F \to G$ is the natural transformation $f^{\sharp}:G^{\sharp} \to F^{\sharp}$ which at $A\in \A$ is 
$(f^{\sharp})_A = (f_{A^*})^*$.
The double dual identification $\eta_F:F \to F^{\sharp\sharp}$ at $A\in \A$ is
$(\eta_F)_A=\can_{F(A^{**})}\circ F(\can_A): FA \to F^{\sharp\sharp}A=F(A^{**})^{**}$.
\item
\label{dfn:FunEx:item3}
There is an isomorphism of symmetric bilinear functors
$$\hat{\phantom{f}}: [F,G^{\sharp}] \stackrel{\cong}{\longrightarrow} [F*,*G]: g \mapsto \hat{g}=(G\can)^*\circ g_*$$
with inverse
the map $f \mapsto f_*\circ F(\can)$,
where the symmetry on the right is
$$[F*,*G] \to [G*,*F]: \ffi \mapsto \tilde{\ffi} = (F\can)^*\circ \ffi_*^*\circ \can_{G*}.$$
We have $\tilde{\tilde{\ffi}}=\ffi$ and $\ffi^*_A\circ \can_{GA} = \tilde{\ffi}_{A^*}\circ G(\can_A)$ for all $A\in \A$.
\item
\label{dfn:FunEx:item4}
For any exact functor $F:\A \to \B$, any natural transformation $\ffi:F* \to *F$  and any arrow $a:A \to A^*$ in $\A$ we have
$$\overline{\ffi_A\circ F(a)} = \tilde{\ffi}_A\circ F(\bar{a}).$$
\item
\label{dfn:FunEx:item5}
On objects $G \in \Fun_{ex}(\A,\B)$, the set $Q(G)$ of quadratic forms on the exact functor $G:\A \to \B$ is the set 
$$Q(G) = \{(q,\ffi)|\  (G, q,\ffi): (\A,Q_{\A}) \to (\B,Q_{\B})\  \text{exact form functor}\}$$
of pairs $(q,\ffi)$ 
of natural transformations $q:Q_{\A} \to Q_{\B}G$ and $\ffi:G \to G^{\sharp}$ 
such that $(G,q,\ffi)$ is an exact form functor \cite[Definition 2.6]{myForm1} from $(\A,*,\can_{\A}, Q_{\A})$ to $(\B,*,\can_{\B}, Q_{\B})$.
Recall that this means that $\ffi^{\sharp}\eta_G=\ffi$, and the following diagram commutes for all $A\in \A$
$$\xymatrix{\A(A,A^*) \ar[r]^{\tau} \ar[d]_{\hat{\ffi}_A\circ G(\phantom{a})} & Q(A) \ar[d]^{q_A} \ar[r]^{\rho} & \A(A,A^*) \ar[d]^{\hat{\ffi}_A\circ G(\phantom{a})}\\
\B(GA,\ast GA) \ar[r]_{\hspace{2ex}\tau} & Q(GA)  \ar[r]_{\hspace{-3ex}\rho} & \B(GA,\ast GA).
}$$
On morphisms $f:F \to G$ in $\Fun_{ex}(\A,\B)$, the map $Q(f): Q(G) \to Q(F)$ sends  $(q,\ffi) \in Q(G)$ to $(Q_{\B}(f)\circ q, f^{\sharp}\ffi  f) \in Q(F)$.
\item
\label{dfn:FunEx:item6}
For $G\in  \Fun_{ex}(\A,\B)$, the functorial Mackey-functor of quadratic forms 
$$[G, G^{\sharp}]\stackrel{\tau}{\longrightarrow} Q(G) \stackrel{\rho}{\longrightarrow} [G, G^{\sharp}]$$
has restriction $\rho(q,\ffi) = \ffi$ and transfer $\tau(\ffi) = (q,\ffi+\ffi^{\sharp}\eta_G)$, where the natural transformation $q:Q_{\A} \to Q_{\B}G$  is defined by
$$q_A(\xi)= \tau(\hat{\ffi}_A\circ G\rho(\xi))\ \in Q_{\B}(GA),$$
for $A\in \A$ and $\xi \in Q_{\A}(A)$.
\end{enumerate}

\begin{lemma}
\label{lem:FunExIsFormCat}
The tuple (\ref{eqn:FunExFomCat}) is an exact form category.
\end{lemma}

\begin{proof}
It is clear that the triple $(\Fun_{ex}(\A,\B),\sharp,\can)$ is an exact category with duality and that Definition \ref{dfn:ExFormCat} (\ref{item1:dfn:ExFormCat}) holds.
In order to verify Definition \ref{dfn:ExFormCat} (\ref{item2:dfn:ExFormCat}), we let $f,g:F \to G$ be two natural transformations of exact functors $F,G:\A\to \B$.
Then, for $(q,\ffi)\in Q(G)$, we have to check the equality
\begin{equation}
\label{eqn:lem:FunExIsFormCat1}
[Q(f+g)-Q(f)-Q(g)](q,\ffi)=\tau(g^{\sharp}\circ \rho(q,\ffi)\circ f)\ \  \in Q(F).
\end{equation}
The left hand side is
$$
\renewcommand\arraystretch{1.5}
\begin{array}{rl}
&(Q_{\B}(f+g)q-Q_{\B}(f)q-Q_{\B}(g)q,\ (f+g)^{\sharp}\ \ffi\ (f+g)-f^{\sharp}\ffi f-g^{\sharp}\ffi g) \\
=&(\tau(g^*\ \rho(q)\ f), \ f^{\sharp} \ffi  g + g^{\sharp} \ffi  f),
\end{array}
$$
whereas the right hand side of (\ref{eqn:lem:FunExIsFormCat1}) is 
$$
\renewcommand\arraystretch{1.5}
\begin{array}{rcl}
\tau(g^{\sharp}\ffi f)&=&(\tau(\widehat{g^{\sharp}\ffi f}\circ F(\rho)), \ g^{\sharp}\ffi f + (g^{\sharp}\ffi f)^{\sharp}\eta_F)\\
&=&(\tau[ (F\can)^*\circ g^{\sharp}_*\ffi_*f_*\circ F(\rho)], g^{\sharp}\ffi f + f^{\sharp}\ffi g)
\end{array}
$$
in view of $(g^{\sharp}\ffi f)^{\sharp}\eta_F = f^{\sharp}\ffi^{\sharp}g^{\sharp\sharp}\eta_F = f^{\sharp}(\ffi^{\sharp}\eta_G) g$ and $\ffi^{\sharp}\eta_G=\ffi$.
In particular, the two expressions in (\ref{eqn:lem:FunExIsFormCat1}) agree on the second component.
At $A\in \A$ and $\xi\in Q_{\A}(\xi)$, the first component of the left hand side of (\ref{eqn:lem:FunExIsFormCat1}) is $\tau_{FA}$ applied to
$$g^*_A\circ  \rho(q_A(\xi))\circ f_A = g^*_A \circ \hat{\ffi}_A \circ G(\rho(\xi))\circ f_A=g^*_A \circ G(\can_A)^*\ffi_{A^*} \circ G(\rho(\xi))\circ f_A,$$
whereas the first component of the right hand side of (\ref{eqn:lem:FunExIsFormCat1}) is $\tau_{FA}$ applied to
$$(F\can_A)^*\circ g^{\sharp}_{A^*}\ffi_{A^*}f_{A^*}\circ F(\rho(\xi)).$$
The last two expressions agree since we have 
$G(\rho(\xi))\circ f_A = f_{A^*}\circ F(\rho(\xi))$ and $g^*_A \circ G(\can_A)^*=(F\can_A)^*\circ g^{\sharp}_{A^*}$ in view of the naturality of $f$ and $g$.

Finally, we need to check condition (\ref{item3:dfn:ExFormCat}) of Definition \ref{dfn:ExFormCat}, that is, for an exact sequence of functors $0 \to F_1 \stackrel{f}{\to} F_2 \stackrel{g}{\to} F_3 \to 0$ from $A$ to $\B$, we have to verify the exactness of the sequence
\begin{equation}
\label{eqn:lem:FunExIsFormCat2}
0 \to Q(F_3) \stackrel{Q(g)}{\longrightarrow} Q(F_2) \stackrel{(Q(f),\rho \circ f)}{\longrightarrow} Q(F_1)\times [F_1,F_2^{\sharp}].
\end{equation}
The map $Q(g)$ is injective because its component functions $Q_{\B}(g)$ and $\ffi \mapsto g^{\sharp}\ffi g$ are injective as $g$ is an admissible epimorphism.
For exactness at $Q(F_2)$, it is clear that the composition of the two maps in (\ref{eqn:lem:FunExIsFormCat2}) is zero.
Now consider $(q,\ffi) \in Q(F_2)$ such that $0=(Q_{\B}(f)\circ q, f^{\sharp}\ffi f) \in Q(F_1)$ and $0=\ffi f\in [F_1,F_2^{\sharp}]$.
Then there is a unique natural map $\psi:F_3 \to F_3^{\sharp}$ such that $g^{\sharp} \psi g =\ffi$ since $F_3$ is a cokernel of $f$ and $F_3^{\sharp}$ is a kernel of $f^{\sharp}$. 
By uniqueness and the equation $\ffi^{\sharp}\eta_{F_2} = \ffi$, we have $\psi^{\sharp}\eta_{F_3} = \psi$.
In view of the exact sequence 
$$0 \to Q_{\B}(F_3A) \stackrel{Q_{\B}(g)}{\longrightarrow} Q_{\B}(F_2A) \stackrel{(Q_{\B}(f), \rho \circ f)}{\longrightarrow} Q_{\B}(F_1A) \times \B(F_1A,(F_2A)^*)$$
for $A \in \A$,
the natural transformation $q:Q_{\A} \to Q_{\B}\circ F_2$ uniquely factors through $Q_{\B}(g): Q_{\B}\circ F_3 \hookrightarrow Q_{\B}\circ F_2$ since $0=Q_{\B}(f)\circ q$, by assumption, and for $\xi\in Q_{\A}(A)$ we have
$$\renewcommand\arraystretch{1.5}
\begin{array}{rcl}
\rho_{F_2A}(q_A(\xi))\circ f_A &=& \hat{\ffi}_A\circ  F_2(\rho_A\xi)\circ f_A \\
&=& F_2(\can_A)^*\circ \ffi_{A^*}\circ F_2(\rho_A\xi)\circ f_A\\
&=& F_2(\can_A)^*\circ F_2^{\sharp}(\rho_A\xi)\circ  \ffi_{A} \circ f_A\\
&=&0
\end{array}$$
as $\ffi_A f_A=0$.
Hence, there is a unique natural transformation $r: Q_{\A} \to Q_{\B}\circ F_3$ such that $Q_{\B}(g)\circ r = q$.
Then $(r,\psi)\in Q(F_3)$ with $Q(g)(r,\psi)=(q,\ffi)$.
Indeed, the diagram 
$$\xymatrix{\A(A,A^*) \ar[r]^{\tau_A} \ar[d]_{\widehat{\psi}_A\circ F_3(\phantom{a})} & Q_{\A}(A) \ar[d]^{r_A} \ar[r]^{\rho_A} & \A(A,A^*) \ar[d]^{\widehat{\psi}_A\circ F_3(\phantom{a})}\\
\B(F_3A,\ast F_3A) \ar[r]_{\hspace{2ex}\tau_{F_3A}} & Q_{\B}(F_3A)  \ar[r]_{\hspace{-3ex}\rho_{F_3A}} & \B(F_3A,\ast F_3A),}$$
 when composed with the injective map to $(\tau_{F_2A},\rho_{F_2A})$ with components $\B(g,g^*)$, $Q_{\B}(g)$ and $\B(g,g^*)$, yields the map from $(\tau_{A},\rho_{A})$ to $(\tau_{F_2A},\rho_{F_2A})$ with components $\hat{\ffi}_A\circ F_2$, $q_A$ and $\hat{\ffi}_A\circ F_2$.
\end{proof}

\begin{remark}[Terminology]
According to Lemma \ref{lem:FunExIsFormCat}, if $\A$ and $\B$ are exact form categories then so is $\Fun_{ex}(\A,\B)$.
Note that if the dualities on $\A$ and $\B$ are strong then so is the duality on $\Fun_{ex}(\A,\B)$.
In this case, a non-singular exact form functor from $\A$ to $\B$, that is, a non-singular form functor \cite[Definition 2.5]{myForm1} which is exact, is the same is quadratic space in $\Fun_{ex}(\A,\B)$, see Definition \ref{dfn:FunExFomCat}.
Therefore, any standard fact and terminology that apply to any exact form category with strong duality also apply to $\Fun_{ex}(\A,\B)$.
For instance, we have the notion of a {\em hyperbolic exact form functor}, that is, a hyperbolic object in $(\Fun_{ex}(\A,\B),\sharp,\eta,Q)$.
We have the notion of a {\em metabolic exact form functor}, that is, a metabolic object in $(\Fun_{ex}(\A,\B),\sharp,\eta,Q)$.
Similarly, we have the notions of a {\em Lagrangian and a sublagrangian functor of an exact form functor} \cite[Definition 2.24]{myForm1}, {\em othogonal decompositions} of non-singular exact form functors \cite[Lemma 2.28]{myForm1}, {\em sublagrangian construction} \cite[Lemma 2.29]{myForm1}...
\end{remark}
\end{definition}

\begin{lemma}
\label{lem:QofCartSq2}
Let
$$\xymatrix{
A \ar[r]^s \ar[d]_f & B \ar[d]^g \\ C \ar[r]_t & D
}$$
be a bicartesian square in an exact form category $(\E,\ast,\can,Q)$, 
that is,
the following sequence is admissible exact
$$\xymatrix{
0 \ar[r] & A \ar[r]^{\hspace{-3ex}\left(\begin{smallmatrix} s \\ -f\end{smallmatrix}\right)} & B \oplus C \ar[r]^{\hspace{2ex}(g \ t)} & D \ar[r] & 0.}$$
Then the following square of abelian groups is cartesian (that is, a pull-back)
$$\xymatrix{
Q(D) \ar[d]_{\left(\begin{smallmatrix} g^{\bullet} \\ t^{\bullet}\end{smallmatrix}\right)} \ar[r]^{t^*\circ\rho(\phantom{x})\circ  g}  &  \E(B,C^*)\ar[d] \\
Q(B) \oplus Q(C) \ar[r]  & Q(A) \oplus \E(B,A^*) \oplus \E(A,C^*)
}$$
where the right vertical and lower horizontal maps are the obvious one's (no sign):
$$
\left(\begin{matrix}  \tau(f^*\circ(\phantom{x})\circ s) \\ 
f^*(\phantom{x}) \\ 
(\phantom{x})\circ s  \end{matrix}\right)
\hspace{3ex}\text{and}\hspace{3ex}
\left(\begin{matrix} s^{\bullet} & f^{\bullet} \\ 
s^*\circ \rho(\phantom{x}) & 0 \\ 
0 & \rho(\phantom{x}) \circ f 
\end{matrix}\right)$$
\end{lemma}

\begin{proof}
This is a consequence of Definition \ref{dfn:ExFormCat} (\ref{item3:dfn:ExFormCat}).
\end{proof}

\section{Additivity for exact form categories}
\label{sec:AddtyExForm}

\noindent
{\bf Hyperbolic and forgetful form functors}.
For any exact category $\E$ we have defined in \cite[\S 3.5]{myHermKex} an exact category with strict duality $\H\E$, called the {\em hyperbolic category} of $\E$, whose $GW$-space is equivalent to the $K$-theory space of $\E$.
As exact category, $\H\E$ is $\E \times \E^{op}$, and the duality functor is $(X,Y) \mapsto (Y,X)$.
Any exact category with duality $(\A,\sharp,\can)$ can be considered as an exact form category $(\A,\sharp,\can,Q)$ where $Q(X) = \A(X,X^{\sharp})^{C_2}$ is the abelian group of symmetric forms on $X$; see \cite[Example 2.24]{myForm1}.
In this sense, $\H\E$ is considered an exact form category.
Recall from \cite[Remark 4.3 and Proposition 4.7]{myHermKex} the homotopy equivalences 
$$\QF(\H\E) \simeq \Q(\E)\hspace{3ex}\text{and}\hspace{3ex} GW(\H\E) \simeq K(\E)$$
where the first is the equivalence of categories $\QF(\H\E) \simeq \Q\E: (X,Y,\ffi) \mapsto X$.

For any exact form category $(\E,\sharp,\can,Q)$ we have exact form functors
$$\H\E \stackrel{H}{\longrightarrow} (\E,\sharp,\can,Q) \stackrel{F}{\longrightarrow} \H\E,$$
called {\em hyperbolic form functor} and {\em forgetful form functor}, respectively.
The hyperbolic form functor $H$ sends the object $(X,Y)$ to $X \oplus Y^{\sharp}$, an arrow $(f,g)$ to $f\oplus g^{\sharp}$, has duality compatibility $\left(\begin{smallmatrix} 0 & 1\\ \can & 0\end{smallmatrix}\right): Y \oplus X^{\sharp} \to X^{\sharp} \oplus Y^{\sharp\sharp}$, and on quadratic forms it is
the composition
$$\xymatrix{
Q_{\H\E}(X,Y) = \E(X,Y) \ar[rr]^{f\mapsto \left(\begin{smallmatrix} 0 & 0\\ \can_Yf & 0 \end{smallmatrix}\right)} && \E(X\oplus Y^{\sharp}, X^{\sharp}\oplus Y^{\sharp\sharp}) \ar[r]^{\hspace{5ex}\tau} & Q_{\E}(X\oplus Y^{\sharp}).}$$
The forgetful form functor $F$ sends the object $X$ of $\E$ to $(X,X^{\sharp})$, the arrow $f$ of $\E$ to $(f,f^{\sharp})$, has duality compatibility $(1,\can):(X^{\sharp},X^{\sharp\sharp}) \to (X^{\sharp},X)$, and on quadratic forms it is the map
$$Q_{\E}(X) \stackrel{\rho}{\longrightarrow} Q_{\H\E}(X,X^{\sharp}) = \E(X,X^{\sharp}).$$
 
For any exact form category $(\E,\sharp,\can,Q)$ we define the exact form functors
\begin{equation}
\label{eqn:S2toH}
(S_2\E,Q) \stackrel{F}{\longrightarrow} \H(S_2\E) \stackrel{\ev_{01}}{\longrightarrow} \H\E
\end{equation}
where the second functor is obtained from the evaluation $S_2\E \to \E: X \mapsto X_{01}$ at the object $(0, 1) \in \Ar[2]$ by functoriality of the hyperbolic category construction.

\begin{theorem}[Additivity for exact categories]
\label{thm:exAddty}
Let $(\E,*,\can,Q)$ be an exact form category with strong duality.
Then the following functor is a weak equivalence
$$\QF (S_2\E,Q) \stackrel{\sim}{\longrightarrow} \Q\E: (X,\xi) \mapsto X_{01}.$$
In particular, the form functor (\ref{eqn:S2toH}) 
induces a weak equivalence of spaces
$$GW(S_2\E,Q) \stackrel{\sim}{\longrightarrow} GW(\H\E) \simeq K(\E).$$
\end{theorem}

\begin{proof}
We will write $\QF(S_2\E)$ for $\QF(S_2\E,Q)$ omitting in the notation the quadratic functor $Q$. 
In view of the description of the exact form category $(S_2\E,Q)$ given in Example \ref{ex:S2explicit},
the category $\QF S_2\E$ is equivalent to the category of tuples $(a:A_0 \rightarrowtail A, \alpha)$ where $(A,\alpha)$ is a quadratic space in $(\E,Q)$ with Lagriangian $a:A_0 \rightarrowtail A$.
Recall \cite[Definition 2.5]{myForm1} that the latter means that $a^{\bullet}(\alpha) = 0$, and the admissible monomorphism $a:A_0 \rightarrowtail A$ induces an isomorphism $A_0 \cong \ker(a^*\circ\ffi)$ where $\ffi=\rho(\alpha)$.
An arrow $(a:A_0 \rightarrowtail A, \alpha) \to (b:B_0 \rightarrowtail B, \beta)$ in $\QF S_2\E$ is an arrow in Quillen's $\Q S_2\E$ from  the exact sequence $(a,a^*\ffi)$ to $(b,b^*\psi)$ represented by a diagram
\begin{equation}
\label{eqn:QhS2Maps}
\xymatrix{
A_0 \ar[d]_a 
  & C_0 \ar@{->>}[l]_{p_0} \ar[d]^c \hspace{1ex} \ar@{>->}[r]^{j_0} 
     & B_0 \ar[d]^b\\
   A \ar[d]_{a^*\ffi}
    & C \ar[d]^{c_1} \hspace{1ex}  \ar@{->>}[l]^{p}  \ar@{>->}[r]_j 
       & B \ar[d]^{b^*\psi}\\
A_0^*
   & C_1 \ar@{->>}[l]^{p_1} \hspace{1ex} \ar@{>->}[r]_{j_1} 
     & B_0^*
}
\end{equation}
where the columns are admissible exact sequences and  $\psi = \rho(\beta)$.
The diagram defines a map in $\QF S_2\E$ if and only if 
$j^{\bullet}(\beta) = p^{\bullet}(\alpha)$, and the following two squares are bicartesian
$$\xymatrix{
(x) 
 & C \hspace{1ex} \ar@{>->}[r]^j  \ar@{->>}[d]_p 
  & B \ar@{->>}[d]^{j^*\psi} 
    & & (y)& C_0 \hspace{1ex} \ar@{>->}[r]^{j_0}  \ar@{->>}[d]_{p_0}
      &  B_0 \ar@{->>}[d]^{j_1^*\can} \\
& A \hspace{1ex} \ar@{>->}[r]_{p^*\ffi}  \ar@{}[ur]|-{\Box}
  & C^*, 
    & & & A_0 \hspace{1ex} \ar@{>->}[r]_{p^*_1\can} \ar@{}[ur]|-{\Box}
      & C_1^*.
}$$ 
Composition is composition in ${\Q} S_2\E$.

We want to show that the functor 
$$G:{\QF}S_2\E \to \Q\E: (a:A_0 \rightarrowtail A, \alpha) \mapsto A_0$$
is a weak equivalence.
We will use Quillen's Theorem A and show that for every object $X\in \Q\E$, the comma category
$\C = (G\downarrow X)$ is non-empty and contractible.
Let $\C'\subset \C$ be the full subcategory of objects $\M = (M,f:GM \to X)$ such that $f= p^{!}$ for some admissible epimorphism $p:X \twoheadrightarrow GM$ in $\E$.
The category $\C'$ has an initial object, namely $(0,0^!:0 \to X)$.
Hence $\C'$ is contractible (and $\C$ is non-empty).
We will show that the inclusion $\C' \subset \C$ has a left adjoint and hence is a weak equivalence finishing the proof.
To show the existence of a left adjoint, we have to construct for every object $\M \in \C$ a universal arrow  $u:\M \to \scN$ in $\C$ with $\scN \in \C'$, that is, an arrow $u:\M \to \scN$ with target in $\C'$ such that for every $\M_0 \in \C'$  every arrow $g:\M \to \M_0$ in $\C$ factors uniquely through $u$.

Given an object $\M = (M,f:GM\to X)$ in $\C$ where  $M=(k:K \rightarrowtail M,\alpha)$ is an object of ${\QF}S_2\E$, and $f:GM=K \to X$ is a map in $\Q\E$.
The map $f$ has a representation as $f= e^!\circ t_!$ where  $t:K \rightarrowtail L$ is an admissible monomorphism and $e:X \twoheadrightarrow L$ is an admissible epimorphism in $\E$.
We construct the object 
$$\scN = (N,e^!:GN=L \to X)$$
of $\C'$ as follows.
It consists of the object $N=(l:L\rightarrowtail N,\beta)$ of ${\QF}S_2\E$ with quadratic form $\beta \in Q(N)$, associated symmetric form $\rho(\beta) = \psi$  and Lagrangian $l:L \subset N$.
Their definitions can be read off the diagrams below where a $\Box$ denotes a bicartesian square:
\begin{equation}
\label{eqn2:thm:exAddty}
\xymatrix{
K \hspace{1ex} \ar@{>->}[r]^s \ar@{>->}[dr]_k \ar[ddd]_t
  & P \ar@{->>}[d]^m \ar@{->>}[r]^r 
     & L^* \ar@{->>}[ddd]^{t^*}
      &
       &  K \hspace{1ex} \ar@{>->}[r]^s \ar[ddd]_t \ar@{}[dr]_{\Box}
         & P \ar[d]^n \ar@{->>}[r]^r 
           & L^* \ar@{->>}[ddd]^{t^*}\\
& M \ar[d]^{\ffi}_{\wr}
  && &
    &N \ar[d]_{\psi}^{\wr}
       &\\
& M^* \ar[d]_{m^*} \ar@{->>}[dr]^{k^*} \ar@{}[uur]_{\Box}
   &&
     && N^* \ar[d]^{n^*} \ar@{->>}[uur]_{l^*}   \ar@{}[dr]^{\Box}
      &\\
 L \hspace{1ex} \ar@{>->}[r]_{r^*\can}  \ar@{}[uur]^{\Box}
 &  P^* \ar@{->>}[r]_{s^*} 
   & K^* ,
     &
      & L \hspace{1ex} \ar@{>->}[r]_{r^*\can}  \ar[uur]^{l}
       & P^* \ar@{->>}[r]_{s^*} 
         & K^*.
}\end{equation}
In detail, let $P$ be the pull-back of $k^*\ffi$ along the admissible epimorphism $t^*$, and denote by $r: P \to L^*$ and $m:P \to M$ the induced maps; see the upper right part of the left diagram.
The maps $r$ and $m$ are admissible epimorphisms, since $k^*\ffi$ and $t^*$ are.
Since $(k,k^*\ffi)$ is an admissible exact sequence, there is a unique admissible monomorphism $s:K \to P$ such that $k=ms$ and $(s,r)$ is an admissible exact sequence.
Since $\ffi^*\can=\ffi$, the left diagram commutes.
Let $N$ be the push-out of $s$ along $t$ with induced maps $n:P \to N$ and $l:L \to N$ being admissible monomorphisms since $s$ and $t$ are.
The dual of the push-out square is a pull-back square, and we obtain a unique map $\psi:N \to N^*$ such that $n^*\psi \ n=m^* \ffi\ m$ and $l^*\psi\ l=0$; making the right diagram commute.
The map $\psi$ is an isomorphism, by the Five Lemma, since it is part of the map of admissible exact sequences $(1,\psi,1)$ from $(l,l^*\psi)$ to $(\psi l, l^*)$.
It satisfies $\psi^*\can=\psi$, by uniqueness of $\psi$ and the symmetry of the right diagram, comparing that diagram with its dual.
Lemma \ref{lem:QofCartSq2} applied to the cocartesian square defining $N$ (upper left square in the right diagram) implies that there is a unique form $\beta \in Q(N)$ such that
$$l^*\circ \rho(\beta)\circ n = r, \hspace{2ex}n^{\bullet}(\beta) = m^{\bullet}(\alpha),\hspace{2ex} l^{\bullet}(\beta) =0.$$
By the uniqueness of $\psi$ satisfying the same equations (Lemma \ref{lem:QofCartSq2} appied to symmetric forms), it follows that $\rho(\beta)=\psi$, and $\beta$ is non-degenerate.
Hence, $\scN$ is an object of $\C'$.
The middle diagram in 
\begin{equation}
\label{eqn:mapMuN}
\xymatrix{
&&K \ar[d]_k & K \ar@{->>}[l]_{1} \ar[d]^s \hspace{1ex} \ar@{>->}[r]^t  & L \ar[d]^l & &P\ \ \ar@{>->}[r]^{n} \ar@{->>}[d]_m & N \ar@{->>}[d]^{n^*\psi} \\
\M \stackrel{u}{\longrightarrow} \scN: && M \ar[d]_{k^*\ffi}& P \ar@{}[ur]|-{\Box} \ar[d]^r \hspace{1ex}  \ar@{->>}[l]_m  \ar@{>->}[r]_n & N \ar[d]^{l^*\psi}&& M \  \ar@{>->}[r]^{m^*\ffi} & P^*\\
&&K^* \ar@{}[ur]|-{\Box} & L^* \ar@{->>}[l]^{t^*} \hspace{1ex} \ar@{>->}[r]_1 & L^*,&&
}
\end{equation}
defines an arrow
$u: \M \to \scN$ in $\C$ since the right commutative diagram in (\ref{eqn:mapMuN}) is bicartesian, verifying condition $(x)$; condition $(y)$ being clear.
Indeed, the right diagram in (\ref{eqn:mapMuN}) is bicartesian because the horizontal maps are admissible monomorphisms, the vertical maps are admissible epimorphisms, and the induced map between vertical kernels is an isomorphism as this is the case when composed with the isomorphism induced by the map $r$, by the diagrams in (\ref{eqn2:thm:exAddty}). 

So, $u : \M \to \scN$ defines a map in $\C$ with target in $\C'$, and we will show that this map is universal.
To that end, let 
$$g:\M \to \M_0 = (M_0,f_0:GM_0 \to X)$$
 be an arrow in $\C$ with target in $\C'$.
So, $M_0 = (k_0:K_0 \rightarrowtail M_0,\alpha_0)$ is an object of ${\QF}S_2\E$ and $f_0:GM_0 = K_0 \to X$ is an arrow in $Q\E$ of the form $f_0 = e_0^!$ where $e_0:X \twoheadrightarrow K_0$ is an admissible epimorphism in $\E$.
Moreover, $\alpha_0 \in Q(M_0)$ has associated non-degenerate symmetric form $\ffi_0=\rho(\alpha_0)$, and $k_0:K_0 \subset M_0$ is a Lagrangian for $\alpha_0$.
The arrow $g$ in $\C$ is is given by a map $\M \to \M_0$ in ${\QF}S_2\E$ represented by a diagram
\begin{equation}
\label{eqn:MM0map}
\xymatrix{
&&K \ar[d]_k 
  & K_1 \ar@{->>}[l]_{p_K} \ar[d]^{k_1} \hspace{1ex} \ar@{>->}[r]^{j_K} 
     & K_0 \ar[d]^{k_0}\\
\M \stackrel{g}{\longrightarrow} \M_0: 
  && M \ar[d]_{k^*\ffi}
    & M_1 \ar[d]^{q_1} \hspace{1ex}  \ar@{->>}[l]^{p}  \ar@{>->}[r]_j 
       & M_0 \ar[d]^{k_0^*\ffi_0}\\
&&K^* 
   & Q_1 \ar@{->>}[l]^{p_Q} \hspace{1ex} \ar@{>->}[r]_{j_Q} 
     & K_0^*
}
\end{equation}
such that $f_0\circ G(g)=f$ in $\Q\E$.
Diagram (\ref{eqn:MM0map}) defines a map in ${\QF}S_2\E$ if and only if the columns are admissible exact sequences in $\E$, 
we have $j^{\bullet}(\alpha_0) = p^{\bullet}(\alpha)$, and the following two squares are bicartesian
$$\xymatrix{
(a) 
 & M_1 \hspace{1ex} \ar@{>->}[r]^j  \ar@{->>}[d]_p 
  & M_0 \ar@{->>}[d]^{j^*\ffi_0} 
    & & (b)& K_1 \hspace{1ex} \ar@{>->}[r]^{j_K}  \ar@{->>}[d]_{p_K}
      &  K_0 \ar@{->>}[d]^{j_Q^*\can} \\
& M \hspace{1ex} \ar@{>->}[r]_{p^*\ffi}  \ar@{}[ur]|-{\Box}
  & M_1^*, 
    & & & K \hspace{1ex} \ar@{>->}[r]_{p^*_Q\can} \ar@{}[ur]|-{\Box}
      & Q_1^*.
}$$ 
The equality $f_0\circ G(g)=f$ in $\Q\E$ means that there is an admissible epimorphism $e_1: K_0 \twoheadrightarrow L$ such diagram $(c)$ commutes and the  square $(d)$ is bicartesian:

$$\xymatrix{
(c) & 
X  \ar@{->>}[r]^{e_0}   \ar@{->>}[dr]_{e}
  & K_0  \ar@{->>}[d]^{e_1}  
   && (d)
    & K_1 \hspace{1ex} \ar@{>->}[r]^{j_K}  \ar@{->>}[d]_{p_K} 
      & K_0 \ar@{->>}[d]^{e_1} \\
& & L,
      &&
       & K \hspace{1ex} \ar@{>->}[r]_{t}  \ar@{}[ur]|-{\Box}
         & L.
}$$ 
From the bicartesian squares $(b)$ and $(d)$, we see that there is a unique isomorphism $h: L \cong Q_1^*$ such that $h \circ e_1 = j_Q^*\can$ and $h\circ t = p_Q^*\can$.
Using this isomorphism, we can replace $Q_1$ in diagram (\ref{eqn:MM0map}) by $L^*$ and obtain the following diagram representing the arrow $g:\M \to \M_0$ 
\begin{equation}
\label{eqn:MM0map2}
\xymatrix{
&&K \ar[d]_k 
  & K_1 \ar@{->>}[l]_{p_K} \ar[d]^{k_1} \hspace{1ex} \ar@{>->}[r]^{j_K} 
     & K_0 \ar[d]^{k_0}\\
\M \stackrel{g}{\longrightarrow} \M_0: 
  && M \ar[d]_{k^*\ffi}
    & M_1 \ar[d]^{q} \hspace{1ex}  \ar@{->>}[l]^{p}  \ar@{>->}[r]_j 
       & M_0 \ar[d]^{k_0^*\ffi_0}\\
&&K^* 
   & L^* \ar@{->>}[l]^{t^*} \hspace{1ex} \ar@{>->}[r]_{e_1^*} 
     & K_0^*
}
\end{equation}
where $q=h^*\can_{Q_1} q_1$.

Let us first prove that there is at most one arrow $v:\scN \to \M_0$ in $\C$ such that $g=v\circ u$.
So, assume we have such an arrow.
Then $v$ is represented by a diagram
\begin{equation}
\label{eqn:NM0map}
\xymatrix{
&&L \ar[d]_l 
  & K_2 \ar@{->>}[l]_{p_L} \ar[d]^{k_2} \hspace{1ex} \ar@{>->}[r]^{j_L} 
     & K_0 \ar[d]^{k_0}\\
\scN \stackrel{v}{\longrightarrow} \M_0: 
  && N \ar[d]_{l^*\psi}
    & M_2 \ar[d]^{q_2} \hspace{1ex}  \ar@{->>}[l]^{p_N}  \ar@{>->}[r]_{j_N} 
       & M_0 \ar[d]^{k_0^*\ffi_0}\\
&&L^* 
   & Q_2 \ar@{->>}[l]^{p_{Q_2}} \hspace{1ex} \ar@{>->}[r]_{j_{Q_2}} 
     & K_0^*
}
\end{equation}
with columns admissible exact sequences such that $j_N^{\bullet}(\alpha_0) = p_N^{\bullet}(\beta)$ and the following two diagrams are bicartesian
$$\xymatrix{
(e) 
 & M_2 \hspace{1ex} \ar@{>->}[r]^{j_N}  \ar@{->>}[d]_{p_N} 
  & M_0 \ar@{->>}[d]^{j^*_N\cdot \ffi_0} 
    & & (f)& K_2 \hspace{1ex} \ar@{>->}[r]^{j_L}  \ar@{->>}[d]_{p_L}
      &  K_0 \ar@{->>}[d]^{j_{Q_2}^*\can} \\
& N \hspace{1ex} \ar@{>->}[r]_{p_N^*\cdot\psi} \ar@{}[ur]|-{\Box}
  & M_2^* 
    & & & L \hspace{1ex} \ar@{>->}[r]_{p^*_{Q_2}\can} \ar@{}[ur]|-{\Box}
      & Q_2^*
}$$ 
Since $\scN, \M_0\in \C'$, we have the equality $e_0^!\circ G(v) = e^!$ in $\Q \E$ which implies that $G(v)=e_2^!$ for some admissible epimorphism in $\E$.
Thus, we can choose $j_L=1_{K_0}$ and $p_L=e_2:K_0 \twoheadrightarrow L$ such that $e_2e_0=e$.
Since $e_0$ is an admissible epimorphism, this uniquely determines the map $e_2$, and we have $e_2=e_1$.
Since the diagram $(f)$ is bicartesian, we can replace the lower horizontal line in (\ref{eqn:NM0map}) with $L^* \stackrel{1}{\leftarrow} L^* \stackrel{e_1^*}{\rightarrowtail} K_0^*$ and obtain a representation of $v$ as follows
\begin{equation}
\label{eqn:NM0mapV2}
\xymatrix{
&&L \ar[d]_l 
  & K_0 \ar@{->>}[l]_{e_1} \ar[d]^{k_2} \hspace{1ex} \ar@{>->}[r]^{1} 
     & K_0 \ar[d]^{k_0}\\
\scN \stackrel{v}{\longrightarrow} \M_0: 
  && N \ar[d]_{l^*\psi} \ar@{}[ur]|-{\Box}
    & M_2 \ar[d]_{q_2} \hspace{1ex}  \ar@{->>}[l]^{p_N}  \ar@{>->}[r]^{j_N} 
       & M_0 \ar[d]^{k_0^*\ffi_0}\\
&&L^* 
   & L^* \ar@{->>}[l]^{1} \hspace{1ex} \ar@{>->}[r]_{e_1^*}  \ar@{}[ur]|-{\Box}
     & K_0^*.
}
\end{equation}
Since $g=v\circ u$, it follows that the exact sequence $(k_1,q)$ is an admissible subobject of the exact sequence $(k_2,q_2)$ as in the following diagram
\begin{equation}
\label{eqn:DfnM2}
\xymatrix{
K_1 \hspace{1ex} \ar@{>->}[r]^{k_1}  \ar[d]_{j_K} 
  &  M_1 \ar@{->>}[r]^{q} \ar[d]^{j_M} 
     & L^* \ar[d]^1 \\
K_0 \hspace{1ex} \ar@{>->}[r]_{k_2} \ar@{}[ur]|-{\Box}
  & M_2 \ar@{->>}[r]_{q_2}
    & L_*
}\end{equation}
in which the vertical arrows are admissible monomorphisms, and $j = j_Nj_M$.
Since the right vertical arrow in (\ref{eqn:DfnM2}) is an isomorphism, $M_2$ is the push-out of $k_1$ along $j_K$, and $k_2$, $j_N$, $j_M$ are the induced maps.
Up to unique isomorphism, this determines $M_2$, and once $M_2$ is chosen, this determines $k_2$ and $j_N$.
Again, from the equation $g = v\circ u$, there is an admissible epimorphism $\bar{q}:M_1 \twoheadrightarrow P$ and a bicartesian diagram
\begin{equation}
\label{eqn:PNasPushout}
\xymatrix{
M_1 \hspace{1ex} \ar@{>->}[r]^{j_M}  \ar@{->>}[d]_{\bar{q}}  
  & M_2 \ar@{->>}[d]^{p_N}  \\
P  \hspace{1ex} \ar@{>->}[r]_{n} \ar@{}[ur]|-{\Box}
  & N 
}\end{equation}
such that 
\begin{equation}
\label{eqn:dfnBarq}
p=m\bar{q}, \hspace{4ex}q = r\bar{q}.
\end{equation}
By (\ref{eqn:mapMuN}), $P$ is a pullback with induced maps $r$ and $m$.
Hence, the map $\bar{q}$ is determined by the last two equations.
By (\ref{eqn:DfnM2}), $M_2$ is a pushout, and the map $p_N:M_2 \to N$ is determined by the equations
\begin{equation}
\label{eqn:dfnPN}
p_N\circ k_2 = l \circ e_1, \hspace{4ex}p_N\circ j_M  = n\circ \bar{q}
\end{equation}
from diagrams  (\ref{eqn:NM0mapV2}) and (\ref{eqn:PNasPushout}).
Hence, the map $p_N$ is uniquely determined.
Finally, the map $q_2$ is determined by the right lower square in diagram (\ref{eqn:NM0mapV2}).
This shows that the arrow $v$ is unique if it exists.

Now we prove that $v$ exists.
So, we are given the maps $u$ and $g$ represented by the diagrams (\ref{eqn:mapMuN}) and (\ref{eqn:MM0map2}) such that the squares $(a)$ and $(d)$ are bicartesian and the triangle $(c)$ commutes.
Taking pull-back of the lower left corner and push-out of the upper right corner of diagram (\ref{eqn:MM0map2}) yields the following factorization of that diagram
\begin{equation}
\label{eqn:MM0map2Factorization}
\xymatrix{
 K \ar[d]_k
  & K  \ar@{->>}[l]_1 \ar[d]_s
   & K_1 \ar@{->>}[l]_{p_K} \ar[d]^{k_1} \hspace{1ex} \ar@{>->}[r]^{j_K} 
    & K_0 \ar[d]^{k_2} \hspace{1ex} \ar@{>->}[r]^{1} 
     & K_0 \ar[d]^{k_0}\\
M \ar[d]_{k^*\ffi}
  & P \ar[d]^{r}   \ar@{->>}[l]_{m} \ar@{}[ur]|-{\Box}
    & M_1 \ar[d]^{q} \hspace{1ex}  \ar@{->>}[l]^{\bar{q}}  \ar@{>->}[r]_{j_M} \ar@{}[ur]|-{\Box}
     & M_2 \ar[d]_{q_2} \hspace{1ex}   \ar@{>->}[r]^{j_N}
       & M_0 \ar[d]^{k_0^*\ffi_0}\\
K^* \ar@{}[ur]|-{\Box}
   & L^* \ar@{->>}[l]^{t^*} 
     & L^* \hspace{1ex} \ar@{>->}[r]_{1}  \ar@{->>}[l]^{1}
      & L^* \hspace{1ex} \ar@{>->}[r]_{e_1^*}  \ar@{}[ur]|-{\Box}
     & K_0^*
}
\end{equation}
in which the columns are admissible exact sequences in $\E$, $p=m\circ \bar{q}$ and $j  = j_N\circ j_M$.
Since $M_2$ and $N$ are defined as push-outs, we can define the map $p_N:M_2 \to N$ as the map of horizontal push-outs of the left hand diagram in (\ref{eqn:DfnMapPN}) and obtain the right hand diagram of bicartesian squares in (\ref{eqn:DfnMapPN}):
\begin{equation}
\label{eqn:DfnMapPN}
\xymatrix{
M_1 \ar@{->>}[d]_{\bar{q}} 
  &  \text{\hspace{1ex}$K_1$}  \ar@{>->}[l]_{k_1}  \hspace{1ex} \ar@{>->}[r]^{j_K}  \ar@{->>}[d]_{p_K}  
    & K_0  \ar@{->>}[d]^{e_1}
     && M_1  \hspace{1ex} \ar@{>->}[r]^{j_M}  \ar@{->>}[d]_{\bar{q}} 
      & M_2 \ar@{->>}[d]_{p_N}
       &  \text{\hspace{1ex}$K_0$}   \ar@{>->}[l]_{k_2} \ar@{->>}[d]^{e_1} \\
P \ar@{}[ur]|-{\Box}
  &  \text{\hspace{1ex}$K$}  \ar@{>->}[l]^{s}  \hspace{1ex} \ar@{>->}[r]_{t}  \ar@{}[ur]|-{\Box}
    & L,
     & & P \ar@{}[ur]|-{\Box}  \hspace{1ex} \ar@{>->}[r]_{n} 
      & N \ar@{}[ur]|-{\Box}
       & \text{\hspace{1ex}$L$.}  \ar@{>->}[l]^{l} 
}
\end{equation}    
In particular, $p_N$ is an admissible epimorphism.
This defines diagram (\ref{eqn:NM0mapV2}) representing a map $v$ in $\Q S_2\E$ such that $v\circ u = g$ in $\Q S_2\E$.
It remains the show that $v:\scN \to \M_0$ defines a map in ${\QF}$. 
Then we will have $g=v\circ u$ in $\C$.
So, we have to show that $p_N^{\bullet}(\beta) = j_N^{\bullet}(\alpha_0)$ and that the diagrams 
\begin{equation}
\label{eqn:BicartCheck}
\xymatrix{
  M_2 \hspace{1ex} \ar@{>->}[r]^{j_N}  \ar@{->>}[d]_{p_N} 
  & M_0 \ar@{->>}[d]^{j^*_N\cdot \ffi_0} 
     & & K_0 \hspace{1ex} \ar@{>->}[r]^{1}  \ar@{->>}[d]_{e_1}
      &  K_0 \ar@{->>}[d]^{e_1^{**}\can} \\
 N \hspace{1ex} \ar@{>->}[r]_{p_N^*\cdot\psi} 
  & M_2^* 
   & & L \hspace{1ex} \ar@{>->}[r]_{\can} 
     & L^{**}
}
\end{equation}
corresponding to $(x)$ and $(y)$ are bicartesian.
To verify the equation $p_N^{\bullet}(\beta) = j_N^{\bullet}(\alpha_0)$, we apply Lemma \ref{lem:QofCartSq2} to the cocartesian square defining $M_2$ (third top square in (\ref{eqn:MM0map2Factorization})).
Thus, we have to show
$$
\renewcommand\arraystretch{1.5}
\begin{array}{cccc}
(1) & j_M^{\bullet}\ p_N^{\bullet}(\beta) & = & j_M^{\bullet}\ j_N^{\bullet}(\alpha_0)\\
(2) & k_2^{\bullet}\ p_N^{\bullet}(\beta) & = & k_2^{\bullet}\ j_N^{\bullet}(\alpha_0)\\
(3) & k_2^* \circ \rho(p_N^{\bullet}(\beta) ) \circ j_M& = & k_2^*  \circ \rho( j_N^{\bullet}(\alpha_0)) \circ j_M.
\end{array}
$$
Recall that we have $p^{\bullet}(\alpha) = j^{\bullet}(\alpha_0)$ and $m^{\bullet}(\alpha) = n^{\bullet}(\beta)$ from the fact that $u$ and $g$ are arrows in ${\QF}S_2\E$.
Now $(1)$ follows from
$$ j_M^{\bullet}\ p_N^{\bullet}(\beta) = \bar{q}^{\bullet}\ n^{\bullet}(\beta) =  \bar{q}^{\bullet}\ m^{\bullet}(\alpha) =  p^{\bullet}(\alpha) =  j^{\bullet}(\alpha_0) = j_M^{\bullet}\ j_N^{\bullet}(\alpha_0).$$
Equation $(2)$ follows from
$$k_2^{\bullet}\ p_N^{\bullet}(\beta) = e_1^{\bullet}\ l^{\bullet}(\beta) = 0,\hspace{1ex}
k_2^{\bullet}\ j_N^{\bullet}(\alpha_0) = k_0^{\bullet}\ (\alpha_0) = 0.$$
Equation $(3)$ follows from
$$\renewcommand\arraystretch{1.5}
\begin{array}{rl}
& k_2^* \circ \rho(p_N^{\bullet}(\beta) ) \circ j_M = k_2^*\ p_N^*\ \psi\  p_N\ j_M = e_1^*\ l^*\ \psi \ n \ \bar{q} = e_1^*\ r\ \bar{q} = e_1^*\ q \\
= &k_0^*\ \ffi_0\  j  = k_2^*\ j_N^*\ \ffi_0\  j_N \ j_M = k_2^*  \circ \rho( j_N^{\bullet}(\alpha_0)) \circ j_M.
\end{array}
$$
The right diagram in (\ref{eqn:BicartCheck}) is clearly bicartesian, and the left one is isomorphic to the left square in the following commutative diagram
$$\xymatrix{
M_2 \ar@{->>}[r]^{\psi\cdot p_N} \ar[d]_{j_N} 
  & N^* \ar@{->>}[r]^{l^*} \ar[d]^{p_N^*} 
    & L^* \ar[d]^{e_1^*} \\
M_0 \ar@{->>}[r]_{j_N^*\cdot \ffi_0} 
  & M_2^*  \ar@{->>}[r]_{k_2^*} 
    & K_0^*
    }$$
 in which the right square is the bicartesian square dual to the right most square in (\ref{eqn:DfnMapPN})
 and the composition is the bicartesian square which is the lower right square of (\ref{eqn:MM0map2Factorization}).
 To see that the horizontal compositions are as claimed we note that,
 in view of the push-out definition of $M_2$ (third top square in (\ref{eqn:MM0map2Factorization})), we have
 $$l^*\ \psi\  p_N =q_2,\hspace{3ex}k_2^*\ j_N^*\ \ffi_0 = k_0\ \ffi_0$$
 since
 $l^*\cdot \psi \cdot  l \cdot e_1 =0$, $q_2\cdot k_2 = 0$, $k^* \psi \cdot n \cdot \bar{q} = r\cdot \bar{q} = q = q_2\cdot j_M$, and $k_0=j_N\cdot k_2$.
  It follows that the left square in  (\ref{eqn:BicartCheck}) is also bicartesian.
 \end{proof}

For any $i\leq j$, the map $[0] \to \Ar[n]: 0 \mapsto (i\leq j)$ induces by restriction
the exact evaluation functor $S_n\A\to \A: A \mapsto A_{ij}$ evaluating at $i\leq j$.
The map $[0] \to \Ar[n]: 0 \mapsto (i\leq j)$ preserves dualities if $j=n-i$ in which case restriction along $[0] \to \Ar[n]$ defines an exact form functor $(S_n\A,Q)\to (\A,Q): A \mapsto A_{ij}$.
For a form functor $G:(\A,Q) \to S_n(\B,Q)$, we will denote by $G_{ij}: (\A,Q) \to (\B,Q)$ the composition of $G$ and evaluation at $i\leq j$ if $j=n-i$.
For an exact functor $G:\A \to \B$ we denote by $H(G):(\A,\sharp,\can,Q) \to (\B,\sharp,\can,Q)$ the hyperbolic form functor, that is, the form functor corresponding to the hyperbolic form on $G$ in the exact form category $\Fun_{ex}(\A,\B)$ of exact functors $\A \to \B$ of Lemma \ref{lem:FunExIsFormCat}.
This is also the composition
$$(\A,\sharp,\can,Q) \stackrel{F}{\longrightarrow} \H\A \stackrel{\H(G)}{\longrightarrow} \H\B \stackrel{H}{\longrightarrow} (\B,\sharp,\can,Q).$$
A direct definition can be found in \cite[Example 2.13]{myForm1}.

\begin{theorem}[Additivity for functors]
\label{thm:AddtyForExFunctors}
Let $(A,\sharp,\can,Q)$ and $(\B,\sharp,\can,Q)$ be exact form categories with strong duality.
\begin{enumerate}
\item
\label{prop:AddtyForFunctors:1}
Let $G=(G,\ffi_q,\ffi): (\A,Q) \to (S_2\B,Q)$ be a non-singular exact form functor.
Then the two non-singular exact form functors 
$G_{02}$ and $H(G_{01})$
induce homotopic maps of  Grothendieck-Witt spaces
$$G_{02} \sim H(G_{01}): GW(\A,Q) \to GW(\B,Q).$$
\item
\label{prop:AddtyForFunctors:2}
Let $G = (G,\ffi_q,\ffi):  (\A,Q) \to (S_3\B,Q)$ be a non-singular exact form functor. 
Then the two non-singular exact form functors 
$G_{03}$ and  $G_{12} \perp H(G_{01})$
induce homotopic maps of Grothendieck-Witt spaces
$$G_{03}\sim (G_{12} \perp H(G_{01})): GW(\A,Q) \to GW(\B,Q).$$
\end{enumerate}
\end{theorem}

\begin{proof}
This is {\em mutatis mutandis} the same as \cite[Proposition 6.7]{myForm1} replacing \cite[Theorem 5.1]{myForm1} with Theorem \ref{thm:exAddty}. Note that the Grothendieck-Witt spaces of exact form categories with strong duality are group complete commutative $H$-spaces, by \cite[Theorem 6.5]{myForm1}.
\end{proof}

\begin{corollary}
\label{cor:SnAddAddty}
Let $(\A,\sharp,\can,Q)$ be an exact form category with strong duality.
Then the form functors
{\footnotesize
$$
\renewcommand\arraystretch{2} 
\begin{array}{ll}
(S_{2n}\A,Q) \to  (\H\A)^n: &  (A,\xi) \mapsto
H(A_{01},A_{12},...,A_{n-1,n})\\
(S_{2n+1}\A,Q) \to (\H\A)^n\times (\A,Q): &
(A,\xi)\mapsto
H(A_{01},A_{12},...,A_{n-1,n}),(A_{n,n+1},\xi_{n,n+1}) 
\end{array}
$$
}
induce homotopy equivalences of Grothendieck-Witt spaces
$$GW(S_{2n}\A,Q) \stackrel{\sim}{\longrightarrow} GW(\H\A)^n \simeq K(\A)^n$$
$$GW(S_{2n+1}\A,Q)  \stackrel{\sim}{\longrightarrow}  GW(\H\A)^n\times GW(\A,Q)\simeq K(\A)^n \times GW(\A,Q).$$
\end{corollary}

\begin{proof}
More precisely, the first functor is the forgetful functor $F: (S_{2n}\A,Q) \to \H (S_{2n}\A)$ followed by $\H$ applied to the exact functor $A \mapsto (A_{01},A_{12},...,A_{n-1,n})$.
The second functor on the first $n$ factors is the forgetful functor $F: (S_{2n+1}\A,Q) \to \H (S_{2n+1}\A)$ followed by $\H$ applied to the exact functor $A \mapsto (A_{01},A_{12},...,A_{n-1,n})$, and on the last factor it is the evaluation $(A,\xi)\mapsto (A_{n,n+1},\xi_{n,n+1})$. 
The corollary immediately follows from Theorem \ref{thm:AddtyForExFunctors} and is {\em mutatis mutandis}
 the same as \cite[Proposition A.8]{myJPAA}.
\end{proof}

\section{Cofinality for exact form categories}

A full inclusion $\A \subset \U$ of exact categories is called {\em fully exact} if it is extension closed, preserves and detects admissible exact sequences.
It is called {\em cofinal} if it is fully exact and for every object $X$ of $\U$ there is an object $Y$ of $\U$ such that $X\oplus Y$ is isomorphic to an object of $\A$.

Let $(\U,\sharp,\can,Q)$ be an exact form category with strong duality.
Let $\A \subset \U$ be a cofinal exact subcategory closed under the duality.
We define the abelian group 
$$GW_0(\U,\A,Q)$$
as the quotient monoid of the monoid of isometry classes of quadratic spaces in $\U$ modulo the submonoid of isometry classes of quadratic spaces in $\A$.
The abelian monoid $GW_0(\U,\A,Q)$ is indeed a group since for any quadratic space $(X,q)$ in $
\U$, the quadratic space $(X,q)\perp (X,q) \perp HY$ lies in $\A$ whenever $X\oplus Y \in \A$.
In particular, 
$$-[X,q]= [X,q] + [HY] \in GW_0(\U,\A,Q).$$
Note that the group $GW_0(\U,\A,Q)$ is independent of the exact structure of $\U$.
Similarly, one defies the abelian group $K_0(\U,\A)$ as the quotient abelian monoid of the monoid of isomorphism classes of objects in $\U$ modulo the submonoid of isomorphism classes in $\A$.

Let $(\E,Q)$ be an exact form category with strong duality.
By \cite[Definition 2.27 and Theorem 6.5]{myForm1},  the group $GW_0(\E,Q) = \pi_0GW(\E,Q)$ is the abelian group generated by isometry classes $[X]$ of quadratic spaces in $(\E,Q)$ modulo the relations $[X\perp Y]=[X]+[Y]$ and $[M]=[H(L)]$ for a metabolic space $M$ in $(\E,Q)$ with Lagrangian $L$.

Let $(M,\xi)$ be a metabolic space in $(\U,Q)$ with Lagrangian $L$.
Choose $K\in \U$ such that $K\oplus L \in \A$.
Then $M \oplus H(K) \in \A$ since it is an extension of $L\oplus K$ and $(L\oplus K)^{\sharp}$ both of which are in $\A$.
In $GW_0(\U,\A,Q)$ we therefore have 
$[M,\xi] = -[H(K)] = [H(L)]$ which shows that the map $GW_0(\U,Q) \to GW_0(\U,\A,Q): [M,\xi] \mapsto [M,\xi]$ in the following lemma is well-defined.

\begin{lemma}
\label{lem:beforeExCof}
Let $(\U,\sharp,\can,Q)$ be an exact form category with strong duality.
Let $\A \subset \U$ be a cofinal exact subcategory closed under the duality.
Then the following sequence of abelian groups is exact
$$GW_0(\A,Q) \to GW_0(\U,Q) \to GW_0(\U,\A,Q) \to 0$$
\end{lemma}

\begin{proof}
Let $B$ be the cokernel of the first map.
Since the composition is zero, we have a natural map $B \to GW_0(\U,\A,Q): [M] \mapsto [M]$.
But the map $GW_0(\U,\A,Q) \to B: [M] \to [M]$ is clearly well-defined and an inverse to the map above.
\end{proof}

The left map in the lemma is also injective.
This is a special case of the following.

\begin{theorem}[Cofinality for exact form categories]
\label{thm:CofinalityExCats}
Let $(\U,\sharp,\can,Q)$ be an exact form category with strong duality.
Let $\A \subset \U$ be a cofinal exact subcategory closed under the duality.
Then the inclusion $\A \subset \U$ induces an isomorphism 
$$GW_i(\A,Q) \cong GW_i(\U,Q)\hspace{3ex} \text{for }i\geq 1$$ and a monomorphism 
$GW_0(\A,Q) \rightarrowtail GW_0(\U,Q)$
in degree $0$.
\end{theorem}

\begin{proof}
In \cite[Definition 2.15]{myForm1}, the group completion of $i\Quad(\E,Q)$ of the category of quadratic spaces in an exact form category with strong duality $(\E,\sharp,\can,Q)$ 
was denoted $GW^{\oplus}(\E,Q)$.
We will use that if the realisation of a simplicial commutative $H$-space is group complete then group completing degree-wise that simplicial $H$-space doesn't change the homotopy type of the realisation.
In particular, group completing degree-wise the commutative square of simplicial categories
$$\xymatrix{
i\Quad(\RR_{\bullet}\A,Q) \ar[r] \ar[d] & i\Quad(\RR_{\bullet}\U,Q) \ar[d] \\
i\RR_{\bullet}\A \ar[r]  & i\RR_{\bullet}\U}$$
yields (after topological realisation) the homotopy equivalent square which is the (topological realisation of the) left square of the commutative diagram of simplicial spaces
$$\xymatrix{
GW^{\oplus}(\RR_{\bullet}\A,Q) \ar[r] \ar[d] & GW^{\oplus}(\RR_{\bullet}\U,Q) \ar[d] \ar[r] & GW_0(\RR_{\bullet}\U,\RR_{\bullet}\A,Q) \ar[d]\\
K^{\oplus}(S_{\bullet}\A) \ar[r]  & K^{\oplus}(S_{\bullet}\U) \ar[r] & K_0(S_{\bullet}\U,S_{\bullet}\A) }$$
in which the rows are homotopy fibrations (degree-wise and after realisation), by Cofinality for $GW^{\oplus}$  \cite[Lemma 2.18]{myForm1} and its $K$-theoretical analogue.
By Additivity (Corollary \ref{cor:SnAddAddty}) for $GW_0$ and $K_0$, the right vertical map is a surjective map of simplicial abelian groups with kernel the constant simplicial abelian group $GW_0(\U,\A,Q)$.
Taking vertical homotopy fibres of the map of homotopy fibrations yields the result.
\end{proof}

\begin{example}[Idempotent and semi-idempotent competion]
\label{ex:IdempCompl}
An exact category is called idempotent complete if every idempotent has a kernel.
The idempotent completion $\tilde{\E}$ of an exact category with duality $(\E,\sharp,\can)$ is canonically an exact category with duality such that $\E \subset \tilde{\E}$ is cofinal \cite[\S 5.1]{myHermKex}. 
Recall that an object of $\tilde{\E}$ is a pair $(A,p)$ with $A\in \E$ and $p^2=p:A \to A$ an idempotent.
Arrows $f:(A,p) \to (B,q)$ in $\tilde{\E}$ are arrows $f:A \to B$ such that $f=fp=qf$.
Composition is composition in $\E$.
We have a fully faithful embedding $\E \subset \tilde{\E}: A \mapsto (A,1)$.
A sequence in $\tilde{\E}$ is called exact if it is a retract of an exact sequence in $\E$.

Let $(\E,\sharp,\can,Q)$ be an exact form category.
We extent $Q$ to $\tilde{\E}$ as follows.
For an object $(A,p)$ of $\tilde{E}$,
the Mackey functor of quadratic forms on $(A,p)$ is the image of the Mackey functor of $A$ under the idempotent induced by $p$.
For instance $Q(A,p) = \im( Q(p): Q(A) \to Q(A))$.
By the Cofinality Theorem \ref{thm:CofinalityExCats}, we have
$$GW_i(\E,Q) \cong GW_i(\tilde{\E},Q),\hspace{2ex}i\geq 1,\hspace{2ex} GW_0(\E,Q) \subset GW_0(\tilde{\E},Q).$$

An exact category is called {\em semi-idempotent complete} if every arrow which has a section also has a kernel.
The {\em semi-idempotent completion} of $\E$ is the full exact subcategory $\tilde{\E}^{0}$ of $\tilde{\E}$ of objects $X\in \tilde{E}$ for which there is $E\in \E$ such that $X \oplus E$ is isomorphic to an object of $\E$.
A sequence in $\tilde{\E}^0$ is exact if it is exact in $\tilde{\E}$.
The inclusion $\E \subset \tilde{\E}^0$ is cofinal and the Cofinality Theorem yields an equivalence of Grothendieck-Witt spaces
$$GW(\E,Q) \stackrel{\sim}{\longrightarrow} GW(\tilde{\E}^0,Q)$$
since $(X,\xi)\perp HA$ is in $\E$ whenever $A, X\oplus A \in E$.

In a semi-idempotent complete category, if $fg$ as an admissible monomorphism, then so is $g$, and if $fg$ is an admissible epimorphism, then so is $f$ \cite[Appendix A]{Keller:chain}.
\end{example}

\section{The filtering localization theorem}

Recall from \cite{mydeloop} that a fully exact inclusion $\A \subset
  \U$ of exact categories is called {\it s-filtering} 
if
\begin{enumerate}
\item
\label{subsec:sFilt:enum:1}
every map $A  \to U$ with $A \in \A$ factors through an admissible subobject of
  $U$ belonging to $\A$, 
\item
\label{subsec:sFilt:enum:2}
every map $U \to A$ with $A$ in $\A$ factors through an admissible quotient
object of $U$ belonging to $\A$,
\item
\label{subsec:sFilt:enum:3}
for every deflation $U \twoheadrightarrow A$ in $\U$
  with $A$ in $\A$, there is an inflation $B \rightarrowtail U$ with $B\in \A$
  such that the composition $B \to A$ is a deflation in $\A$, and
\item
\label{subsec:sFilt:enum:4}
for every inflation $A \rightarrowtail U$ in $\U$ with $A$ in $\A$, there is a
deflation $U \twoheadrightarrow B$ with $B \in \A$ such that the composition
$A \to B$ is an inflation in $\A$.
\end{enumerate}
Items (\ref{subsec:sFilt:enum:1}) and (\ref{subsec:sFilt:enum:2}) imply that
$\A$ is closed under admissible sub-objects and quotient-objects
in $\U$.

Let $\A \subset \U$ be an s-filtering inclusion of exact
  categories.
 To simplify we will assume that $\A$ is idempotent complete.
 A {\em weak isomorphism} \cite[Definition 1.12]{mydeloop} is a finite composition of deflations with kernel in
$\A$ and inflations with cokernel in $\A$.
By definition, weak isomorphisms are closed under composition.
They satisfy a calculus of fractions \cite[Lemma 1.13]{mydeloop}.
In diagrams, weak isomorphisms are displayed as $\stackrel{_{\sim}}{\to}$,
inflations with cokernel in $\A$ as $\stackrel{_{\sim}}{\rightarrowtail}$, and
deflations with kernel in $\A$ as $\stackrel{_{\sim}}{\twoheadrightarrow}$.
By \cite[Proposition 1.16]{mydeloop}, the exact quotient category
  $\U/\A$ exists and is obtained from $\U$ by formally inverting
the set of weak isomorphisms.  
A sequence of $\U/\A$ is a conflation in $\U/\A$ if and only if it is
  isomorphic to the image of a conflation of $\U$.

\begin{lemma}
\label{lem:FiltINcl}
Let $\A \subset \U$ be an $s$-filtering inclusion of exact categories with $\A$ idempotent complete.
Then the following holds.
\begin{enumerate}
\item
\label{item0:lem:FiltINcl}
A map in $\U$ is a weak isomorphism if and only if it is an isomorphism in $\U/\A$.
\item
\label{item1:lem:FiltINcl}
For every weak isomorphism $f:X \stackrel{\sim}{\to} Y$ in $\U$ there is an inflation $s:U \stackrel{\sim}{\rightarrowtail} X$ with cokernel in $\A$ such that $fs$ is an inflation with cokernel in $\A$.
\item
\label{item3:lem:FiltINcl}
For any admissible monomorphisms $u:X \rightarrowtail Y$ and $s:Z \stackrel{\sim}{\rightarrowtail} Y$ in $\U$ with $\coker(s)$ in $\A$ there are admissible monomorphisms $t:W \stackrel{\sim}{\rightarrowtail} X$ and $v:W \rightarrowtail Z$ in $\U$ with $\coker(t)$ in $\A$ such that $ut=sv$.
\item
\label{item4:lem:FiltINcl}
If in a map $f=(f_0,f_1,f_2)$ of admissible exact sequences in $\U$ two of the three maps are weak isomorphisms then so is the third.
\item
\label{item2:lem:FiltINcl}
For every admissible monomorphism $v:V_0 \stackrel{\sim}{\rightarrowtail} V$ in $\U$ with cokernel in $\A$ and every admissible epimorphism $g:V \twoheadrightarrow U$ in $\U$ there is an admissible epimorphism $g_0:V_0 \twoheadrightarrow U_0$ and weak isomorphism $u:U_0 \stackrel{\sim}{\to} U$ in $\U$ such that $gv = ug_0 \in \U(V_0,U)$.
\label{item5:lem:FiltINcl}
Let $f:X \to Y$ be a map in $U$ such that there is a map $g:Y\to Z$with $gf$ an inflation in $\U$. If either $X$ is in $\A$ or $f$ is a weak isomorphism, then $f$ is an inflation.
\end{enumerate}
\end{lemma}

\begin{proof}
Part (\ref{item0:lem:FiltINcl}) is \cite[Lemma 1.17 (6)]{mydeloop}.
Part (\ref{item1:lem:FiltINcl}) is \cite[Lemma 1.17 (3)]{mydeloop}.
Part (\ref{item3:lem:FiltINcl}) is  \cite[Lemma 8.5 (e)]{myHermKex}.
Part (\ref{item4:lem:FiltINcl}) is \cite[Lemma 1.17 (7)]{mydeloop}.
Part (\ref{item2:lem:FiltINcl}) follows from  (\ref{item3:lem:FiltINcl}) and  (\ref{item4:lem:FiltINcl}).
Part (\ref{item5:lem:FiltINcl}) is \cite[Lemma 8.5 (g)]{myHermKex}.
\end{proof}

Let $(\U,\sharp,\can,Q_{\U})$ be an exact form category, and let
$\A\subset \U$  be a duality preserving s-filtering inclusion
of exact categories.
Assume $\A$ is idempotent complete and closed under the duality functor $\sharp$.
Then $(\A,\sharp,\can,Q_{\A})$ is an exact form category and the inclusion
$I:\A \subset \U$ defines a fully faithful form functor
$$(I,1,1): (\A,\sharp,\can,Q_{\A}) \subset (\U,\sharp,\can,Q_{\U})$$
where $Q_{\A}$ is the restriction of $Q_{\U}$ to $\A$. 
The duality $(\sharp,\can)$ on $\U$ induces an exact duality and a natural double dual
identification on the quotient $\U/\A$, and we
consider $\U/\A$ equipped with this structure of an exact category with
duality. 
For an object $X$ of $\U/\A$, we define the abelian group $Q_{\U/\A}(X)$ of forms on $X$ as
$$Q_{\U/\A}(X) = \colim_{Y \stackrel{\sim}{\to}X}Q_{\U}(Y)$$
where the indexing category of the colimit has objects the weak isomorphisms with target $X$ and morphisms the commutative triangles in $\U$.
The indexing category is filtering.
For an arrow $[gs^{-1}]:X \to Y$ in $\U/\A$ represented by the $\U$-diagram $gs^{-1}:X \stackrel{\sim}{\leftarrow} X_0 \to Y$ with $s:X_0 \stackrel{\sim}{\to} X$ a weak isomorphism and $[\xi,r]\in Q_{\U/\A}(Y)$ represented by $\xi \in Q(U)$ for some weak isomorphism $r:U \stackrel{\sim}{\to} Y$, choose a weak isomorphism $r_0:U_0 \stackrel{\sim}{\to} X_0$ and a $\U$-arrow $g_0:U_0 \to U$ such that $rg_0 = gr_0$ and set
$$[gs^{-1}]^{\bullet}\left([\xi,r]\right) = [g_0^{\bullet}(\xi),sr_0]\in Q_{\U/\A}(X).$$
One checks that this is independent of the choices and makes $Q_{\U/\A}$ into a functor from $(\U/\A)^{op}$ to abelian groups.

Since $\U/\A(X,X^{\sharp}) = \colim_{Y \stackrel{\sim}{\to} X}\U(Y, Y^{\sharp})$, the Mackey functor of quadratic forms in $\U/\A$
$$\U/\A(X,X^{\sharp}) \stackrel{\tau}{\longrightarrow} Q_{\U/\A}(X) \stackrel{\rho}{\longrightarrow} \U/\A(X,X^{\sharp})$$
is the colimit 
$$\colim_{Y \stackrel{\sim}{\to}X} \left(\U(Y,Y^{\sharp}) \stackrel{\tau}{\longrightarrow} Q_{\U}(Y) \stackrel{\rho}{\longrightarrow} \U(Y,Y^{\sharp})\right).$$

\begin{lemma}
\label{lem:UAformcat}
Assume that $\A$ is idempotent complete. 
Then the following tuple defines an exact form category:
$$(\U/\A,\sharp,\can,Q_{\U/\A}).$$
\end{lemma}

\begin{proof}
We already know that $(\U/\A,\sharp,\can)$ is an exact category with duality \cite[\S 8.2]{myHermKex} and it is clear from the definition that the conditions in Definition \ref{dfn:ExFormCat} (\ref{item1:dfn:ExFormCat}) and (\ref{item2:dfn:ExFormCat}) hold. 
So we are left with checking condition (\ref{item3:dfn:ExFormCat}) of Definition \ref{dfn:ExFormCat}.
Since every admissible exact sequence in $\U/\A$ is isomorphic to the image of an admissible exact sequence of $\U$ we only need to check the condition in Definition \ref{dfn:ExFormCat} (\ref{item3:dfn:ExFormCat}) for such sequences.
So, let $X \stackrel{i}{\rightarrowtail} Y \stackrel{p}{\twoheadrightarrow} Z$ be an admissible exact sequence in $\U$ and consider the sequence
$$0 \to Q_{\U/\A}(Z) \stackrel{p^{\bullet}}{\longrightarrow} Q_{\U/\A}(Y) \stackrel{(i^{\bullet},\rho(\phantom{\xi}) \circ i)}{\longrightarrow}  Q_{\U/\A}(X)\times \U/\A(Y,X^{\sharp})$$
which we need to show is exact.

Injectivity of $p^{\bullet}$ follows from Lemma \ref{lem:FiltINcl} (\ref{item2:lem:FiltINcl}).
Let $[\xi,t]\in Q_{\U/\A}(Z)$ be represented by $\xi\in Q_{\U}(U)$ for some admissible monomorphism $t:U \stackrel{\sim}{\rightarrowtail} Z$ with cokernel in $\A$ (Lemma \ref{lem:FiltINcl} (\ref{item1:lem:FiltINcl})) such that $p^{\bullet}[\xi,t] =0 \in Q_{\U/\A}(Y)$.
Since $p^{\bullet}[\xi,t] = [g^{\bullet}\xi,t_0]$ is represented by $g^{\bullet}(\xi)\in Q_{\U}(V)$ for $t_0:V \stackrel{\sim}{\rightarrowtail} Y$ the pull-back of $t$ along $p$ and $g: V \twoheadrightarrow U$ the induced admissible epimorphism, there is an admissible monomorphism $v:V_0 \stackrel{\sim}{\rightarrowtail} V$ 
such that $v^{\bullet}g^{\bullet}(\xi)=0\in Q_{\U}(V_0)$.
By Lemma \ref{lem:FiltINcl} (\ref{item2:lem:FiltINcl}), there are an admissible epimorphism $g_0:V_0 \twoheadrightarrow U_0$ and a weak isomorphism $u:U_0 \stackrel{\sim}{\to} U$ such that 
$gv = ug_0 \in \U(V_0,U)$.
Then $0=v^{\bullet}g^{\bullet}(\xi) = g_0^{\bullet}u^{\bullet}(\xi) \in Q_{\U}(V_0)$ implies 
$0=u^{\bullet}(\xi) \in Q_{\U}(U_0)$ since $g_0^{\bullet}$ is injective.
Then $[\xi,t] = [u^{\bullet}(\xi) , tu] =0 \in Q_{\U/\A}(Z)$.

For exactness at $Q_{\U/\A}(Y)$, let $[\xi,t]\in Q_{\U/\A}(Y)$ be represented by $\xi\in Q_{\U}(U)$ and the admissible monomorphism $t:U \stackrel{\sim}{\rightarrowtail} Y$ in $\U$ with cokernel in $\A$.
Assume that $i^{\bullet}[\xi,t]=0$ and $\rho[\xi,t]\circ i=0$.
Using Lemma \ref{lem:FiltINcl} (\ref{item3:lem:FiltINcl}) we find
admissible monomorphisms $t_0:V \stackrel{\sim}{\rightarrowtail} X$ and $j:V\rightarrowtail U$ with $\coker t_0$ in $\A$ and $it_0=tj \in \U(V,Y)$ such that $j^{\bullet}(\xi) =0 \in Q_{\U}(V)$ and $\rho(\xi)j=0$.
Denote by $r:U \to U/V=\coker(j)$ the quotient map and by $t_1:U/V \to Z$ the induced map on quotients.
By Lemma \ref{lem:FiltINcl}  (\ref{item4:lem:FiltINcl}), $t_1$ is a weak isomorphism.
By quadratic left-exactness of $Q_{\U}$, there is a unique $\zeta \in Q_{\U}(U/V)$ such that 
$r^{\bullet}(\zeta) = \xi$. 
Then $[\zeta,t_1] \in Q_{\U/\A}(Z)$ and $p^{\bullet}[\zeta,t_1] = [\xi,t] \in Q_{\U/\A}(Y)$.
\end{proof}

We have the sequence exact form functors
\begin{equation}
\label{eqn:FiltLocSeq}
(\A,\sharp,\can,Q_{\A}) \longrightarrow (\U,\sharp,\can,Q_{\U}) \longrightarrow (\U/\A,\sharp,\can,Q_{\U/\A})
\end{equation}
where the first is the canonical full inclusion, $Q_{\A}$ is the restriction of $Q_{\U}$ to $\A$,  and the second form functor is the localisation functor 
which on quadratic forms is the map
$$Q_{\U}(X) \longrightarrow Q_{\U/\A}(X): \xi \mapsto [\xi,1_X].$$

The rest of the section is devoted to the proof of the following theorem and its Grothendieck-Witt space version in Theorem \ref{thm:cor:filt_locn} below.

\begin{theorem}
\label{thm:filt_locn}
Let $(\U,\sharp,\can,Q_{\U})$ be an exact form category with strong duality, and 
let $\A \subset \U$ be an s-filtering inclusion of exact
categories closed under the duality.
Assume that $\A$ is idempotent complete.
Then the sequence (\ref{eqn:FiltLocSeq})  induces a homotopy fibration on
classifying spaces of
$$\QF(\A,Q_{\A}) \longrightarrow \QF(\U,Q_{\U}) {\longrightarrow} \QF(\U/\A,Q_{\U/\A}).$$
\end{theorem}

\begin{lemma}
\label{lem:UAsymmSpace}
Assume the hypothesis of Theorem \ref{thm:filt_locn}.
Let $(X,\xi)$ be a form in $(\U,Q_{\U})$ with $\rho(\xi)$ a weak isomorphism which is
metabolic in $(\U/\A,Q_{\U/\A})$ with Lagrangian the admissible monomorphism $L_0 \subset X$ in $\U/\A$. 
Then there is a weak isomorphism of admissible exact sequences
in $\U$ of the form
\begin{equation}
\label{eq:UAsymmSpace}
\xymatrix{L \xymono[r]^i \ar[d]_l^{\wr} & X \xyepi[r]^p
  \ar[d]^{\rho(\xi)}_{\wr} & X/L \ar[d]_{\wr}^{l^{\sharp}\can}\\
(X/L)^{\sharp} \xymono[r]^{\hspace{3ex} p^{\sharp}} & X^{\sharp} \xyepi[r]^{i^{\sharp}} & L^{\sharp}}
\end{equation}
such that $i^{\bullet}(\xi)=0 \in Q_{\U}(L)$ and $L$ is isomorphic to $L_0$ as subobjects of $X$ in $\U/\A$.
Moreover, there is a metabolic space $(Y,\zeta)$ in $(\U,Q_{\U})$ with Lagrangian $M$ and an isometry $(X,\xi)\cong (Y,\zeta)$ in $\U/\A$ which takes $L$ to $M$.
\end{lemma}

\begin{proof}
For the underlying symmetric form $\rho(\xi)$, the diagram in the lemma was constructed in \cite[Lemma 8.8 (b)]{myHermKex}.
For this diagram we have $i^{\bullet}(\xi)=0$ in $\U/\A$.
In view of the calculus of fractions and Lemma \ref{lem:FiltINcl} (\ref{item1:lem:FiltINcl}), we can replace $L$ with an admissible subobject in $\U$ which is isomorphic to $L$ in $\U/\A$ such that $i^{\bullet}(\xi)=0$.

To construct the metabolic $Y$ in $\U$, consider the quadratic form 
$$\nu = (1\ 0)^{\bullet}(\xi) + \tau\left(\begin{smallmatrix}0 & p^{\sharp}\\0 & 0\end{smallmatrix}\right) \in Q_{\U}\left(X \oplus (X/L)^{\sharp}\right)$$
on $X \oplus (X/L)^{\sharp}$. 
Denote by $Y$ the cokernel of the admissible monomorphism 
$$\left(\begin{smallmatrix}-i\\ l\end{smallmatrix}\right): L \longrightarrowtail X \oplus (X/L)^{\sharp}$$
in $\U$ and let $q:X \oplus (X/L)^{\sharp} \twoheadrightarrow Y$ be the corresponding quotient map.
Note that 
$$\renewcommand\arraystretch{2}
\begin{array}{rcl}
\nu_{| L} &= & 
\left(\begin{smallmatrix}-i\\ l\end{smallmatrix}\right)^{\bullet}
\left( (1\ 0)^{\bullet}(\xi) + \tau\left(\begin{smallmatrix}0 & p^{\sharp}\\0 & 0\end{smallmatrix}\right) \right)\\
&=& (-i)^{\bullet}(\xi) + \tau\left( (-i^{\sharp}\ l^{\sharp}) \circ
\left(\begin{smallmatrix}0 & p^{\sharp}\\0 & 0\end{smallmatrix}\right) \circ
\left(\begin{smallmatrix}-i\\ l\end{smallmatrix}\right)
\right)\\
&=& 0
\end{array}$$
and 
$$(-i^{\sharp}\ l^{\sharp})\circ \rho(\nu) 
= (-i^{\sharp}\ l^{\sharp}) \circ \left(\begin{matrix}\rho(\xi) & p^{\sharp}\\ \can\circ p  & 0\end{matrix}\right) =0$$
and that therefore there is a unique quadratic form $\zeta \in Q_{\U}(Y)$ on $Y$ such that $q^{\bullet}(\zeta) = \nu$.
The form $\zeta$ is non-degenerate since $\rho(\zeta)$ is the canonical isomorphism from the push-out of the upper left corner of (\ref{eq:UAsymmSpace}) to the pull-back of the lower right corner of (\ref{eq:UAsymmSpace}).
In particular, we have an admissible exact sequence
$$
(X/L)^{\sharp} \stackrel{j}{\longrightarrowtail} Y \stackrel{\rho(\zeta)\circ j^{\sharp}}{\longtwoheadrightarrow} (X/L)^{\sharp\sharp}
$$
where 
$$j  = q\circ \left(\begin{smallmatrix}0\\ 1\end{smallmatrix}\right)\hspace{3ex}\text{and}\hspace{3ex}
j^{\bullet}(\zeta) = \left(\begin{smallmatrix}0\\ 1\end{smallmatrix}\right)^{\bullet}(\nu) = 0$$
which shows that $(Y,\zeta)$ is metabolic in $(\U,Q)$ with Lagrangian $j$.
Finally, in $\U/\A$, the map $l$ is an isomorphism and thus, the canonical map from $X$ to the push-out of the left upper corner of (\ref{eq:UAsymmSpace}) is an isomorphism and we can choose $M=(X/L)^{\sharp}$.
\end{proof}

Recall \cite[Definition 2.27]{myForm1} that the Witt group $W(\E,Q)$ of an exact form category with strong duality is the abelian monoid of isometry classes of quadratic spaces in $(\E,Q)$ under orthogonal sum modulo the submonoid of metabolic spaces.
$W(\E,Q)$ is an abelian group since $(X,\xi) \perp (X,-\xi)$ is metabolic for any quadratic space $(X,\xi)$.
By \cite[Lemma 6.4]{myForm1}, we have an isomorphism $$\pi_0\QF(\E,Q) \cong W(\E,Q): [X,\xi]\mapsto [X,\xi].$$

\begin{lemma}
\label{lem:AUW0ex}
Assume the hypothesis of Theorem \ref{thm:filt_locn}.
Then the following sequence of abelian groups is exact
$$W(\A,Q_{\A}) \to W(\U,Q_{\U}) \to W(\U/\A,Q_{\U/\A}).$$
\end{lemma}

\begin{proof}
The composition is clearly zero.
Let $(X,\xi)$ be a quadratic space in $(\U,Q_{\U})$ which is stably metabolic in
$\U/\A$.
Then there is a metabolic space $(M,\mu)$ in $(\U/\A,Q_{\U/\A})$ such that $(X,\xi)\perp
(M,\mu)$ is metabolic in $(\U/\A,Q_{\U/\A})$.
By definition of $Q_{\U/\A}$, we can assume
$(M,\mu)$ to be a form in $(\U,Q_{\U})$ with $\rho(\mu)$ a weak isomorphism.
The quadratic form $(M,-\mu)$ in $(\U,Q_{\U})$ is also metabolic in $(\U/\A,Q_{\U/\A})$ with same Lagrangian as $(M,\mu)$.
The orthogonal sum of forms  
$(X,\xi) \perp (M,\mu) \perp (M,-\mu)$ in $(\U,Q_{\U})$ is therefore metabolic in $\U/\A$.
Consider the map of forms in $(\U,Q_{\U})$
\begin{equation}
\label{eqn:lem:AUW0ex}
\left(\begin{smallmatrix}-1&1\\ \rho(\mu) & 0\end{smallmatrix}\right): (M,-\mu)
\perp (M,\mu) \longrightarrow H^{\mu}(M)= (M\oplus
M^{\sharp},p_M^{\bullet}(\mu) + h_M)
\end{equation}
where  $h_M=\tau( p_M^{\sharp}\circ p_{M^{\sharp}})$ is the hyperbolic form and 
$p_M:M \oplus M^{\sharp} \to M$,  $p_{M^{\sharp}}:M \oplus M^{\sharp} \to M^{\sharp}$ are the canonical projections.
This map does indeed preserve forms because
$$\renewcommand\arraystretch{2}
\begin{array}{rl}
 &\left(\begin{smallmatrix} -1 & 1 \\ \rho(\mu) & 0 \end{smallmatrix}\right)^{\bullet} 
\left( p_M^{\bullet}(\mu) + h_M\right)\\
= &
\left(\begin{smallmatrix} -1 & 1 \\ \rho(\mu) & 0 \end{smallmatrix}\right)^{\bullet} 
\left( (1\ 0)^{\bullet}(\mu) + \tau\left( \begin{smallmatrix}0 & 1\\ 0 & 0 \end{smallmatrix} \right) \right)\\
=&
 (-1\ 1)^{\bullet}(\mu) 
+ \tau\left(
\left( \begin{smallmatrix}-1 & \rho(\mu)^{\sharp}\\ 1 & 0 \end{smallmatrix} \right)
\left( \begin{smallmatrix}0 & 1\\ 0 & 0 \end{smallmatrix} \right)
\left( \begin{smallmatrix}-1 & 1\\ \rho(\mu) & 0 \end{smallmatrix} \right)
\right)\\
=&
 (1\ 0)^{\bullet}(-1)^{\bullet}(\mu) +  (0\ 1)^{\bullet}(\mu) 
 + \tau\left( 
\left( \begin{smallmatrix}0 \\  1 \end{smallmatrix} \right) \rho(\mu) (-1\ 0)\right)
+\tau
\left( \begin{smallmatrix}- \rho(\mu)& 0\\ \rho(\mu) & 0 \end{smallmatrix} \right)\\
=&
(1\ 0)^{\bullet}(-1)^{\bullet}(\mu) +  (0\ 1)^{\bullet}(\mu) 
+\tau
\left( \begin{smallmatrix}- \rho(\mu)& 0\\ 0 & 0 \end{smallmatrix} \right)\\
=&
(1\ 0)^{\bullet}(-\mu) +  (0\ 1)^{\bullet}(\mu) 
\end{array}
$$
since $(-1)^{\bullet}=-1+\tau\rho$; see \cite[Remark 2.3]{myForm1}.
Note that $H^{\mu}(M)$ is a quadratic space in $(\U,Q_{\U})$ with associated symmetric bilinear form
$\left( \begin{smallmatrix} \rho(\mu)& 1\\ \can & 0 \end{smallmatrix} \right)$.
As a weak isomorphism, the map (\ref{eqn:lem:AUW0ex}) of forms defines an isometry in $(\U/\A,Q_{\U/\A})$.
It follows that the quadratic space
$(Y,\psi) = (X,\xi) \perp H^{\mu}(M)$ in $(\U,Q)$ is a metabolic in $\U/\A$.
By Lemma \ref{lem:UAsymmSpace} there is an
admissible monomorphism $i:L \rightarrowtail Y$ in $\U$ such that $i^{\bullet}(\psi)= 0 \in Q_{\U}(L)$ 
and such that the canonical map $L \to (Y/L)^{\sharp} \cong L^{\perp}$ is a weak isomorphism.
By \cite[Lemma 8.5 (g)]{myHermKex}, the map $L \to
L^{\perp}$ is an admissible monomorphism in $\U$ with cokernel in $\A$. 
In particular, $L$ is a sublagrangian of $(Y,\psi)$, and the induced quadratic form $\zeta$ on $L^{\perp}/L\in \A$ is non-degenerate \cite[Lemma 2.30]{myForm1}.
Since $H^{\mu}(M)$ is metabolic in $(\U,Q_{\U})$, \cite[Lemma 2.30]{myForm1} implies 
the  equalities  
$[X,\xi] = [Y,\psi] = [L^{\perp}/L,\zeta]$ 
in $W(\U,Q_{\U})$.
This shows that $[X,\xi]$  is in the image of 
the map $W(\A,Q_{\A}) \to W(\U,Q_{\U})$.
\end{proof}

For an exact form category with strong duality $(\E,\sharp,\can,Q)$, let $\Q_{m}(\E) \subset \QF(\E,Q)$  be the
full subcategory of metabolic spaces in $(\E,Q)$.
Under the hypothesis of Theorem \ref{thm:filt_locn}, let
$\Q_{\A}(\U) \subset \QF(\U,Q_{\U})$ be the full subcategory of those
quadratic spaces $(X,\xi)$ in $(\U,Q_{\U})$ for which there is a map 
$(A,\alpha) \to (X,\xi)$ in $\QF(\U,Q_{\U})$ with $A$ in $\A$.
Note that we have inclusions 
\begin{equation}
\label{eqn:QFIncl}
\Q_{m}(\U) \subset \Q_{\A}(\U) \subset \QF(\U,Q_{\U}).
\end{equation}

\begin{lemma}
\label{lem:QaU=QUcofinal}
On homotopy groups, the inclusion $\Q_m(\U/\A) \to \QF(\U/\A,Q_{\U/\A})$ and the two inclusions in (\ref{eqn:QFIncl}) induce monomorphisms on $\pi_0$ and isomorphisms on $\pi_i$ for $i\geq 1$.
\end{lemma}

\begin{proof}
The proof is {\em mutatis mutandis} the same as \cite[Lemma 8.10]{myHermKex}.
\end{proof}

The functor $\QF(\U,Q_{\U}) \to \QF(\U/\A,Q_{\U/\A})$ maps $\Q_{\A}(\U)$ into the
subcategory $\Q_m(\U/\A)$ of metabolic objects of $(\U/\A,Q_{\U/\A})$ and thus defines a functor 
$$\lambda: \Q_{\A}(\U) \longrightarrow \Q_m(\U/\A).$$
In fact, $\Q_{\A}(\U)$ is precisely the full subcategory of $\QF(\U,Q)$ of those
objects which have image in $\Q_m(\U/\A)$ in view of Lemma \ref{lem:UAsymmSpace} and \cite[Lemma 8.5 (g)]{myHermKex} and \cite[Lemma 2.30]{myForm1}.
Theorem \ref{thm:filt_locn} follows from Lemmas \ref{lem:AUW0ex},
\ref{lem:QaU=QUcofinal} and the following proposition the proof of which will
occupy the rest of the section.

\begin{proposition}
\label{prop:filt_locn}
Assume the hypothesis of Theorem \ref{thm:filt_locn}.
Then the sequence (\ref{eqn:FiltLocSeq}) induces a homotopy fibration on
classifying spaces of
$$\QF(\A,Q_{\A}) \longrightarrow \Q_{\A}(\U) \stackrel{\lambda}{\longrightarrow} \Q_m(\U/\A).$$
\end{proposition}

\begin{proof}
We will write $\QF(\A)$ for $\QF(\A,Q_{\A})$ and similarly for $\U$ and $\U/\A$.
In order to identify the homotopy fibre of $\lambda$, we will use Theorem B of
  Quillen.
We will show that for a map $g:M \to N$ in $\Q_m(\U/\A)$, the induced
  functor $g^*:(N\downarrow \lambda) \to (M \downarrow \lambda)$ is a
  homotopy equivalence on classifying spaces.
Since every object in $\Q_m(\U/\A)$ is metabolic, every object
  is the target of a map from $0$.
Thus it suffices to show that $g^*$ is a homotopy equivalence for 
  $g: 0 \to N$.
Recall that the category $(N \downarrow \lambda)$ has objects pairs
  $(X,\alpha)$ with $X$ an object of $\Q_{\A}(\U)$ and $\alpha: N \to
  \lambda X$ a map in $\Q_m(\U/\A)$.
Maps from $(X,\alpha)$ to $(Y,\beta)$ correspond to maps $s:X \to Y$ in
$\Q_{\A}(\U)$ such that $\beta = \lambda(s) \circ \alpha$.
Composition is composition of the maps $s$ in $\Q_{\A}(\U)$.

Consider the following diagram of categories and functors
\begin{equation}
\label{diagram1}
\xymatrix{E_g \ar[r]^{\hspace{-3ex} p} \ar[d]_{\gamma} & (N {\downsim} \lambda)
\ar[r]^{i_N} & (N {\downarrow} \lambda) \ar[d]^{\hspace{.5ex} g^*}\\
{\QF(\A)} 
\ar@{=}[r] 
& (0 
{\downsim} 
\lambda) \ar[r]^{i_0} & (0 {\downarrow} \lambda),}
\end{equation}
where the categories and functors are as follows.
\begin{itemize}
\item
The functor $i_N$ is the inclusion
  of the full subcategory $(N {\downsim} \lambda)$ of $(N {\downarrow} \lambda)$
  whose objects are those pairs 
  $(X,\alpha)$ for which $\alpha:N\to \lambda(X)$ is an isomorphism.
Note that $\QF(\A) = (0 {\downsim} \lambda)$ since $\A$ is closed under subquotients in $\U$ and its objects are precisely those of $\U$ which are isomorphic to $0$ in $\U/\A$. 
Note that the categories $(N {\downsim} \lambda)$ are non-empty, by Lemma \ref{lem:UAsymmSpace}.
\item
The category $E_g$ has objects the pairs $(a , \alpha)$ with 
  $a: A \to X$ a map in $\Q_{\A}(\U)$, $A \in \QF(\A)$, and $\alpha:
  N \to \lambda (X)$ an isomorphism in $\QF(\U/\A)$ such that
  $\alpha \circ g=\lambda(a)$ where $\lambda (A)$ is identified with $0$ via the
  unique isomorphism $0 \to \lambda (A)$ in $\QF(\U/\A)$.
A map $(a,\alpha) \to (b:B \to Y,\beta)$ is a pair of maps $c:A \to B$
  and $s:X \to Y$ in $\Q_{\A}(\U)$  such that $sa=bc$ and $\beta = \lambda(s) \alpha$.
Composition is composition of the individual maps of the pairs.
\item
The functor $p: E_g \to (N {\downsim} \lambda)$ 
  sends $(a:A \to X, \alpha)$ to $(X,\alpha)$.
\item
The functor $\gamma: E_g \to \QF(\A)$  sends
  $(a:A \to X, \alpha)$ to $A$.
\end{itemize}
There is a natural transformation of functors\phantom{z} $i_0 \circ \gamma \to
g^* \circ i_N \circ p$\phantom{z} given by 
$$i_0 \circ \gamma (a:A \to X, \alpha) =
 (A,0) \stackrel{a}{\longrightarrow} g^* \circ i_N \circ p (a,\alpha) = (X,
 \alpha \circ g).$$ 
Thus, on classifying spaces, diagram (\ref{diagram1}) commutes up to homotopy.
The proof of Proposition \ref{prop:filt_locn} is complete once we show that
$i_N$ (hence $i_0$), $p$ and $\gamma$ induce 
 homotopy equivalences on classifying spaces.
This will be done in Lemmas \ref{lem:Filt1}, \ref{lem:Filt2} and \ref{lem:Filt3} below.
\end{proof}

We will repeatedly use the Sublagrangian Construction \cite[Lemma 2.30]{myForm1}.
Recall that a sublagrangian of a quadratic space $(X,\xi)$ in an exact form category with strong duality $(\E,\sharp,Q)$ is an admissible subobject $i:L \subset X$ such that $i^{\bullet}(\xi)=0$ and the map $L \to L^{\perp}=\ker(i^{\sharp}\rho(\xi))$ of subobjects of $X$ is an admissible monomorphism.
In this case, there is a unique non-degenerate form $\bar{\xi}\in Q(L^{\perp}/L)$ on $L^{\perp}/L$ such that $p^{\bullet}(\bar{\xi} )= \xi_{|L^{\perp}}\in Q(L^{\perp})$ where $p:L^{\perp} \to L^{\perp}/L$ is the quotient map.
In particular, to give an arrow $Y \to X$ in $\QF(\E,Q)$ is the same as to give a sublagrangian $L\subset X$ of $X$ and an isometry $L^{\perp}/L\cong Y$.

\begin{lemma}
\label{lem:Filt1}
For every object $N$ of $\Q_m(\U/\A)$, the full inclusion 
  $$i_N: (N {\downsim} \lambda) \to (N \downarrow \lambda)$$ is a homotopy
  equivalence on classifying spaces.
\end{lemma}

\begin{proof}
The lemma follows from Theorem A of Quillen once we show that for all
  objects $(X,\alpha)$ of $(N \downarrow \lambda)$, the comma category 
  $\C = (i_N \downarrow (X,\alpha))$ is contractible.
The map $\alpha: N \to \lambda (X)$ is represented by a totally isotropic subspace
$L_0 \subset X$  of $X$ in $\U/\A$  together with an
 isometry $L_0^{\perp}/L_0 \cong N$ in $\U/\A$.
By the calculus of fractions and \cite[Lemmas 8.5 (a) and (f)]{myHermKex}
we can assume that $L_0 \subset X$ is actually a totally isotropic subspace of
$X$ in $(\U,Q_{\U})$.
This shows that $\C$ is non-empty, as
$L_0^{\perp}/L_0$ together with the canonical map in $\QF(\U)$ to $X$  and the
isometry $L_0^{\perp}/L_0 \cong N$ in $\U/\A$ defines an object of $\C$.

The comma category $\C$ is equivalent to the category
with objects the totally isotropic subspaces $L \subset X$ in $(\U,Q_{\U})$ such that
$L$ is isomorphic, as a subobject of $X$ in $\U/\A$, to $L_0$.
There is a (unique) map $L \to L'$ in $\C$ if $L' \subset L$.
Note that in this case, $L'\subset L$ is an admissible monomorphism, by Lemma \ref{lem:FiltINcl} (\ref{item5:lem:FiltINcl}).
Using the calculus of fractions for the set of weak isomorphisms and Lemma
\cite[Lemma 8.5 (a)]{myHermKex}, we see that
$\C$ is a filtering 
category  so that the classifying space of $\C$ is contractible
\cite[\S 1 Corollary 2]{quillen:higherI}.
\end{proof}

\begin{lemma}
\label{lem:Filt2}
The functor $p: E_g \to (N {\downsim} \lambda)$ is a homotopy
  equivalence on classifying spaces.
\end{lemma}

\begin{proof}
In view of Quillen's Theorem A it suffices to show that for any object
$(X,\alpha)$ of $(N\downsim \lambda)$, the comma 
category $\C= (p \downarrow (X,\alpha))$ is contractible.
An object of $\C$ is a sequence 
$A \stackrel{a}{\to} U \stackrel{s}{\to} X$  
of composable maps $a$, $s$ in $\Q_{\A}(\U)$ with $A \in \QF(\A)$
and $s$ an isomorphism in $\QF(\U/\A)$ such that $ \lambda (s a ) = \alpha g$.
A map from $(a ,s)$ to $B \stackrel{b}{\to} V  \stackrel{t}{\to} X$ is a pair
of maps $c:A \to B$, $r: U \to V$ in $\Q_{\A}(\U)$ such that $ra=bc$
  and $s=tr$.
Composition is composition of commutative diagrams in $\Q_{\A}(\U)$.

Let $J: \C' \subset \C$ be
  the full subcategory of those sequences which are of the form $A
  \stackrel{a}{\to} X \stackrel{1}{\to} X$. 
The functor $G: \C \to \C' $ which sends an object $(a,s)$ to $(s\circ a, id)$
and a map $(c,r)$ to $(c,id_X)$ is a homotopy inverse to $J$ because
  $GJ=id$ and because $(id,s): (a,s) \to (s \circ a, id)$ defines a natural
  transformation $id \to JG$.

The category $\C'$ is non-empty and
 filtered, and thus contractible, by exactly
  the same arguments as in the proof of Lemma \ref{lem:Filt1}.
In detail, the map $\alpha g: 0 \to X$ in $\QF(\U/\A)$ is represented by a
Lagrangian 
$L_0$ of $X$ in $\U/\A$ which, by \cite[Lemmas 8.5 (a) and (f)]{myHermKex} we
can assume to be a totally isotropic subspace $L_0$ of $X$ in $(\U,Q_{\U})$ such that
the quotient $L_0^{\perp}/L_0$ is in $\A$.
The category $\C'$ is equivalent to the category whose objects are totally
isotropic subspaces $L \subset X$ in $(\U,Q_{\U})$ such that $L$ is isomorphic to $L_0$
as subobjects of $X$ in $\U/\A$.
Maps are inclusions of totally isotropic subspaces.
This category is filtered by the calculus of fractions and \cite[Lemmas 8.5 (a)]{myHermKex}.
\end{proof}

The following lemma 
completes the proof of Proposition \ref{prop:filt_locn} and thus of Theorem
\ref{thm:filt_locn}.

\begin{lemma}
\label{lem:Filt3}
The functor $\gamma: E_g \to \QF(\A)$ is a homotopy equivalence on
  classifying spaces.
\end{lemma}

\begin{proof}
Again, the lemma follows from Theorem A of Quillen, once we show that for all
  objects $A$ of $\QF(\A)$, the comma category $\C=(A \downarrow \gamma )$ is
  contractible. 
The category $\C$ has objects $(A \stackrel{b}{\to}
  B \stackrel{u}{\to} X, N \stackrel{s}{\to} \lambda X)$ with $b$, $u$ maps
  in $\Q_{\A}(\U)$, $B \in \QF(\A)$ and $s$ an isomorphism in
  $\QF(\U/\A)$ such that $\lambda (u) = sg$.
A map from $(b,u,s)$ to $(A \stackrel{c}{\to}
   C \stackrel{v}{\to} Y, t)$ is a pair of maps $a: B \to C$ and $w:
  X \to Y$ in $\Q_{\A}(\U)$ such that $c=ab$, $va=wu$ and $t = \lambda(w)
  s$. 
Composition is composition of commutative diagrams in $\Q_{\A}(\U)$.

Let $J:\C' \subset \C$ be the
  full subcategory of objects $(b,u,s)$ with $b=id_A$.
The functor $G:\C \to \C'$
  which sends an object $(b,u,s)$ to $(id,ub,s)$ and a map $(a,w)$ to
  $(id_A,w)$ is a homotopy inverse of $J$ because $GJ=id$ and because
  the map $(b,id_X): (id_A,ub,s) \to (b,u,s)$ defines a natural
  transformation of functors $JG \to id$.
We will show that $\C'$ is a non-empty filtered (projective) system, so that
$\C'$ and hence $\C$ will be contractible. 

Omitting the information $b=id_A$, the objects of $\C'$ are pairs
$(A\stackrel{u}{\to} X, N \stackrel{s}{\to}\lambda X)$ with $u$ a map in
$\Q_{\A}(\U)$ and $s$ an isomorphism in $\QF(\U/\A)$ such that $\lambda(u)=sg$.
The category $\C'$ is non-empty.
Indeed, by definition of $Q_m(\U/\A)$ the map $g: 0 \to N$ is represented by a Lagrangian $L\subset N$ in $\U/A$. 
By Lemma \ref{lem:UAsymmSpace}, there is a metabolic space $Y$ with Lagrangian $M$ in $(\U,Q_{\U})$ and an isometry $s:N \cong Y$ in $\U/\A$ carrying $L$ to $M$.
The Lagrangian $M$ of $Y$ defines a map $h:0 \to Y$ in $\Q_{\A}(\U)$.
Since $\lambda (A) = 0$, the pair 
$$(\hspace{.5ex} A \stackrel{1_{A}\perp h}{\longrightarrow} A \perp Y, 
\hspace{1ex} N \stackrel{s}{\longrightarrow} \lambda (Y))$$
defines an object of $\C'$.
This shows that the category $\C'$ is non-empty.
The fact that $\C'$ is a filtered projective system is shown in Lemma \ref{lem:filt18v2} below.
\end{proof}
\vspace{2ex}

Let $(X,\xi)$ be a quadratic space in an exact form category with strong duality $(\E,Q)$.
A map $u=i_!\circ p^!:V \to X$ in $\Q(\E)$ defines an arrow in $\QF(\E,Q)$ with target $(X,\xi)$ if and only if $\ker(p) \subset X$ is a sublagrangian for $\xi$ with orthogonal the source of $p$. 
In this case, there is a unique form on $V$, given by the Sublagrangian Construction \cite[Lemma 2.30]{myForm1}, such that $u:V \to X$ is a map in $\QF(\E,Q)$.
We may write $u^{\bullet}\xi\in Q(V)$ for this unique non-degenerate form and it the restriction of $\xi$ to $V$. 

\begin{lemma}
\label{lem:FromSymToQ}
Let $(\E,Q)$ be an exact form category with strong duality, and let $(\A,Q_{\A})$ and $(\U,Q_{\U})$ as in Theorem \ref{thm:filt_locn}.
\begin{enumerate}
\item
\label{item1:lem:FromSymToQ}
Consider a string of maps 
\begin{equation}
\label{eqn:StringQU}
Y \to V \to X
\end{equation}
in Quillen's $\Q(\E)$.
Let $(X,\xi)$ by a quadratic space in $(\E,Q)$.
Assume that (\ref{eqn:StringQU}) defines a string of maps in the form $Q$-construction for symmetric forms with target $(X,\rho(\xi))$
and that the composition $Y \to X$ defines a map in $\QF(\E,Q)$ with target $(X,\xi)$.
Then (\ref{eqn:StringQU}) defines a string of maps in $\QF(\E,Q)$ with target $(X,\xi)$
\item
\label{item2:lem:FromSymToQ}
Let $u:A \to X$ be an arrow in Quillen's $\Q(\U)$ with $A\in \A$.
Let $\xi_i\in Q_{\U}(X)$, $i=1,2$, be non-degenerate forms such that $u$ defines a map in $\QF(\U,Q_{\U})$ with target $(X,\xi_i)$, $i=1,2$.
Assume that $\rho(\xi_1)=\rho(\xi_2)$, $\xi_1=\xi_2 \in Q_{\U/\A}(X)$ and $u^{\bullet}(\xi_1)=u^{\bullet}(\xi_2)\in Q_{\A}(A)$.
Then there is a factorisation of $u$ as $u=sa$ with $a:A \to Y$ and $s:Y \to X$ such that $a$ and $s$ define maps in $\QF(\U,Q_{\U})$, $s$ is an isomorphism in $\U/\A$, and $s^{\bullet}(\xi_1)=s^{\bullet}(\xi_2)\in Q_{\U}(Y)$.
\end{enumerate}
\end{lemma}

\begin{proof}
For (\ref{item1:lem:FromSymToQ}), we notice that the maps $Y \to V \to X$ are given by sublagrangians $L_0 \subset L_1 \subset X$ where $L_0 \subset X$ is a sublagrangian for $\rho(\xi)$ giving rise to the map $V \to X$ and $L_1$ is a sublagrangian of $\xi$ giving rise to the map $Y \to X$. 
Since $\xi_{|L_1}=0$, we have $\xi_{|L_0}=0$, and $Y \to V \to X$ defines a string of maps in $\QF(\E,Q)$ with target $(X,\xi)$.

We are left with proving (\ref{item2:lem:FromSymToQ}).
The map $u=j_!\circ p^!:A \to X$ is given by a subobject $i:L \subset X$ which is a sublagrangian for $\xi_1$ and $\xi_2$ where $j: L^{\perp} =\ker(i^{\sharp}\ffi) \rightarrowtail X$ for  $\ffi=\rho(\xi_1)=\rho(\xi_2)$ and $p:L^{\perp} \twoheadrightarrow A$ is the quotient map.
The forms $\xi_1$ and $\xi_2$ restrict to the same form on $L^{\perp}$ since 
$j^{\bullet}(\xi_1) = p^{\bullet}(u^{\bullet}\xi_1) = p^{\bullet}(u^{\bullet}\xi_2) =j^{\bullet}(\xi_2)$.
We will construct a factorisation $j=m\circ n$ with $n:L^{\perp} \rightarrowtail M$, $m:M \stackrel{\sim}{\rightarrowtail} X$ admissible monomorphisms and $\coker(m) \in \A$ such that $m^{\bullet}\xi_1=m^{\bullet}\xi_2$ as then $M^{\perp}\subset L$ will be a sublagrangian for $\xi_1$ and $\xi_2$ defining the required factorisation $A \to Y=M/M^{\perp} \to X$.

Since $\xi_1=\xi_2\in Q_{\U/\A}(X)$, there is an admissible monomorphism $K \stackrel{\sim}{\rightarrowtail} X$ with cokernel in $\A$ such that $(\xi_1)_{|K}=(\xi_1)_{|K} \in Q_{\U}(K)$.
By Lemma \ref{lem:FiltINcl} (\ref{item3:lem:FiltINcl}), $K$ and $L^{\perp}$ have a common admissible subobject $K_0$ such that $K_0 \rightarrowtail L^{\perp}$ has cokernel in $\A$.
Let $P$ be the pushout of $K$ and $L$ over their subobject $K_0$.
The induced map $q:P \to X$ is an isomorphism in $\U/\A$ and $q^{\bullet}(\xi_1)=q^{\bullet}(\xi_2)$ in view of Lemma \ref{lem:QofCartSq2}.
The map on quotients $P/L^{\perp} \to X/L^{\perp}$ is an isomorphism in $\U/\A$.
In view of Lemma \ref{lem:FiltINcl} (\ref{item1:lem:FiltINcl}), there is an admissible monomorphism $N \stackrel{\sim}{\rightarrowtail} P/L^{\perp} $ with cokernel in $\A$ such that the composition $N \to  X/L^{\perp}$ is also an admissible monomorphism with cokernel in $\A$.
Let $M$ be the pull-back of $N \to  X/L^{\perp}$ along the admissible epimorphism $P \to P/L^{\perp}$.
Then the induced map $M \to X$ is an admissible monomorphism with cokernel in $\A$ as a pullback of such, $(\xi_1)_{|M}=(\xi_2)_{|M}$ as $M \to X$ factors through $P$, and $M$ contains $L^{\perp}$ as admissible subobject.
\end{proof}

\begin{lemma}
\label{lem:filt18v2}
$\phantom{ds}$
\begin{enumerate}
\item
Let  $u_i:A \to X_i$, $i=1,2$, be two maps in $\QF(\U)$
with $A$ an object of $\QF(\A)$.
Assume that there is an isomorphism $s:\lambda (X_1) \stackrel{\cong}{\to} \lambda
(X_2)$  in $\QF(\U/\A)$ such that $s \circ \lambda(u_1)=\lambda(u_2)$ in $\QF(\U/\A)$.  
Then there are maps $u:A \to V$  and $s_i:V \to X_i$ in $\QF(\U)$ such that 
$s_iu=u_i$, $s\circ \lambda(s_1)=\lambda(s_2)$, and such that $s_i$ is an isomorphism in $\QF(\U/\A)$, $i=1,2$. 
\item
Let $x: A \to X$, $u_i: X \to Y$, $i=1,2$ be maps in $\QF(\U)$ with $A$ in
$\QF(\A)$ and $\lambda(u_i)$ isomorphisms in $\QF(\U/\A)$ satisfying
 $u_1x=u_2x$ in $\QF(\U)$ and $\lambda(u_1)=\lambda(u_2)$ in $\QF(\U/\A)$. 
Then there are maps $z:A \to Z$, $v: Z \to X$ in $\QF(\U)$ such
  that $\lambda(v)$ is an isomorphism in $\QF(\U/\A)$, $u_1v=u_2v$, and $x=vz$.
\end{enumerate}
\end{lemma}

\begin{proof}
Using Lemma \ref{lem:FromSymToQ}, 
Lemma \ref{lem:filt18v2} follows from the symmetric forms case proved in \cite[Lemmas 8.16 and 8.17]{myHermKex}.
\end{proof}

\begin{theorem}
\label{thm:cor:filt_locn}
Let $(\U,\sharp,\can,Q)$ be an exact form category with strong duality, and 
let $\A \subset \U$ be an s-filtering fully exact inclusion of exact
categories closed under the duality.
Assume that $\A$ is idempotent complete.
Then the sequence (\ref{eqn:FiltLocSeq})  induces a homotopy fibration of Grothendieck-Witt spaces
$$GW(\A,Q) \longrightarrow GW(\U,Q) {\longrightarrow} GW(\U/\A,Q_{\U/\A}).$$
\end{theorem}

\begin{proof}
In view of the definition of the Grothendieck-Witt space, this follows from
Theorem \ref{thm:filt_locn} and its $K$-theory analog \cite[Theorem
2.1]{mydeloop}.
\end{proof}

\section{The cone construction}
\label{sec:ConeConstr}

In this section we generalise from symmetric forms to arbitrary forms the cone construction of \cite[\S 9]{myHermKex}.
Let $(\U,\sharp,\can,Q)$ be an exact form category and $\A \subset \U$ be a fully exact subcategory closed under the duality. 
We consider $\A$ as an exact form category by restricting the functors $\sharp$ and $Q$ to $\A$.
In this section, we will construct an exact form category $\C(\U,\A)$ together with a duality preserving fully faithful form functor
 $\U \to \C(\U,\A)$ such that the sequence $\A \to \U \to \C(\U,\A)$ induces a homotopy fibration of Grothendieck-Witt spaces.

\subsection{The form category $\C_0(\U,\A)$}
Let $(\U,\sharp,\can,Q)$ be an exact form category.
Let  $\n$ be the totally ordered set $\n=\N\sqcup \N^{op}$ in which every element
of $\N^{op}$ is greater than every element of $\N$ and the set of non-negative intergers $\N$ is equipped with its usual order.
I will write $0,1,2,...$ for elements of $\N$ and $0^{op},1^{op},2^{op},...$ for elements in $\N^{op}$.
So $0<1< 2 < ... < 2^{op}< 1^{op} < 0^{op}$.
The poset $\n$ is equipped with the strict duality $n \leftrightarrow n^{op}$.
As explained in \S \ref{sec:FormQ} and \cite[Lemma 6.2]{myForm1}, the functor category $\Fun(\n,\U)$ is an exact form category.
For an object $U$ of $\Fun(\n,\U)$ we may write $U_n$ for $U(n)$ and $U^n$ for $U(n^{op})$, $n\in \N$.
Let $\A \subset \U$ be fully exact subcategory 
closed under the duality.
Recall from \cite[\S 9.1]{myHermKex} the category $\C_0(\U,\A)$ which is the full subcategory of the functor category $\Fun(\n,\U)$ of those objects $U:\n \to \U$ such that
\begin{enumerate}
\item
for all $n\leq m \in \N$, the maps $U_n \stackrel{_{\sim}}{\rightarrowtail} U_{m}$
and $U^m \stackrel{_{\sim}}{\twoheadrightarrow} U^n$ 
are admissible monomorphisms with cokernel in $\A$ and admissible epimorphisms with kernel in $\A$,
respectively, and
\item
\label{C0(U,A):enum:2}
there is a $k \in \N$ such that for all $n\in \N$, the maps
$U_n \stackrel{_{\sim}}{\rightarrowtail} U^{n+k}$ are admissible monomorphisms with cokernel in
$\A$, and the maps $U_{n+k} \stackrel{_{\sim}}{\twoheadrightarrow} U^n$
are admissible epimorphisms with kernel in $\A$.
\end{enumerate}
The category $\C_0(\U,\A)$ is closed under extensions and duality in $\Fun(\n,\U)$ and is thus equipped with a structure of exact form category such that the inclusion $\C_0(\U,\A) \subset \Fun(\n,\U)$ preserves and reflects admissible exact sequences, duality and forms.
In particular, for an object $U$ of $\C_0(\U,\A)$, the abelian group of forms $Q(U)$ on $U$ is the set of pairs $(\xi,\ffi)$ where $\xi = (\xi_n)_{n\in \N}$ is a family of forms $\xi_n\in Q(U_n)$ such that $\xi_{n} = U(n \leq m)^{\bullet} \xi_m$ for all $n\leq m \in \N$ and $\ffi: U \to U^{\sharp}$ is a symmetric form such that $\rho(\xi_n) = U(n < n^{op})^{\sharp}\circ \ffi_n$.
When $\U=\A$, we may write $\C_0(\A)$ in place of $\C_0(\A,\A)$.

\subsection{The cone form category $\C(\U,\A)$}
For $k\in \N$, the functors $\eps_k: \n \to \n: n\mapsto n+k,\ n^{op}\mapsto n^{op}$ and $\eps^k:\n \to \n: n\mapsto n,\ n^{op}\mapsto (n+k)^{op}$ together with the unique natural transformations $1 \to \eps_k$ and $\eps^k \to 1$ induce exact functors $ \C(\U,\A) \to \C(\U,\A): U \mapsto U_{[k]}=U\circ \eps_k$, $U \mapsto U^{[k]}=U\circ \eps^k$ and natural transformations $e_{[k]}:U \to U_{[k]}$ and $e^{[k]}: U^{[k]} \to U$.
Note that $e_{[i]}\circ e_{[j]} = e_{[i+j]}$, $e^{[i]}\circ e^{[j]} = e^{[i+j]}$ and $e_{[i]}\circ e^{[j]} = e^{[j]}\circ e_{[i]}$.
The category $\C(\U,\A)$ is obtained from $\C_0(\U,\A)$ by formally inverting the arrows $e_{[k]}$ and $e^{[k]}$ for all objects $U$ of $\C_0(\U,\A)$.
In particular, the objects of $\C(\U,\A)$ are the same as the objects of $\C_0(\U,\A)$, and for $U,V \in \C(\U,\A)$, the set of arrows $U \to V$ in $\C(\U,\A)$ is
$$\Hom_{\C} (U,V)= \colim_{i,j}\Hom_{\C_0}(U^{[i]},V_{[j]})$$
where we abbreviated $\C = \C(\U,\A)$ and $\C_0=\C_0(\U,\A)$.
The colimit is taken over the projective system 
$\cdots \to U^{[3]} \to U^{[2]} \to U^{[1]} \to U$ and over the inductive system
$V \to V_{[1]} \to V_{[2]} \to V_{[3]} \to \cdots$ induced by the maps $e^{[1]}$ and $e_{[1]}$.
Composition of maps is defined as follows.
If the maps $[f]:U \to V$ and $[g]: V \to W$ in $\C(\U,\A)$
are represented by the maps $f:U^{[i]} \to V_{[j]}$ and $g: V^{[k]} \to W_{[l]}$
in $\C_0(\U,\A)$, their composition $[g]\circ [f]$ is represented by the map 
$g_{[j]}\circ f^{[k]}: U^{[i+k]} \to V^{[k]}_{[j]} \to W_{[l+j]}$.
By \cite[Lemma 9.2]{myHermKex}, the category $\C(\U,\A)$ is an exact category with duality where a sequence is admissibly exact if it isomorphic in $\C(\U,\A)$ to the image of an admissible exact sequence of $\C_0(\U,\A)$ under the functor $\C_0(\U,\A)\to\C(\U,\A)$.
The duality functor is induced by that of $\C_0(\U,\A)$ so that $\C_0(\U,\A) \to \C(\U,\A)$ preserves duality.
Since $(U_{[k]})^{\sharp} = (U^{\sharp})^{[k]}$, $(U^{[k]})^{\sharp} = (U^{\sharp})_{[k]}$, 
$e_{[k]}(U)^{\sharp}=e^{[k]}(U^{\sharp})$, $e^{[k]}(U)^{\sharp}=e_{[k]}(U^{\sharp})$ we have
$$\renewcommand\arraystretch{1.5}\begin{array}{rcl}
\Hom_{\C}(U,U^{\sharp})& =& \colim_{i,j}\Hom_{\C_0}(U^{[i]},(U^{\sharp})_{[j]})\\
&=& \colim_{i,j}\Hom_{\C_0}(U^{[i]},(U^{[j]})^{\sharp})\\
&=& \colim_{i}\Hom_{\C_0}(U^{[i]},(U^{[i]})^{\sharp})
\end{array}$$
where the last colimit is over the restrictions along $e^{[1]}$.
Thus, taking colimit over the restrictions along $e^{[1]}$ of the functorial Mackey functors $(\tau,\rho)$ on $\C_0(U,\A)$ defines a functorial Mackey functor
$$\Hom_{\C}(U,U^{\sharp}) \stackrel{\tau}{\longrightarrow} Q_{\C}(U) \stackrel{\rho}{\longrightarrow} \Hom_{\C}(U,U^{\sharp})$$
on $\C(\U,\A)$ where
$$Q_{\C}(U) = \colim_iQ_{\C_0}(U^{[i]})$$
is the colimit over the inductive system
$$Q_{\C_0}(U) \stackrel{e^{\bullet}}{\longrightarrow} Q_{\C_0}(U^{[1]}) \stackrel{e^{\bullet}}{\longrightarrow} Q_{\C_0}(U^{[2]}) \stackrel{e^{\bullet}}{\longrightarrow} \cdots
$$
where we abbreviated $e=e^{[1]}$.
Since a filtered colimit of left exact sequences is left exact, this defines the exact form category 
$(\C(\U,\A), \sharp,\can,Q_{\C})$.
When $\U=\A$, we may write $\C(\A)$ in place of $\C(\A,\A)$ and call it the {\em cone category} of $\A$.

The constant diagram functor defines a duality preserving fully faithful exact form functor
\begin{equation}
\label{eqn:ConstDiag}
c:(\U,\sharp,\can,Q) \longrightarrow (\C(\U,\A), \sharp,\can,Q_{\C})
\end{equation}
which is $s$-filtering \cite[Lemma 9.3]{myHermKex}.
In particular, the quotient exact form category $(\C(\U,\A)/\U, \sharp,\can,Q_{S\U})$ is defined; see Lemma \ref{lem:UAformcat}.
We may write $Q_{S\U}$ for the induced quadratic functor $\Q_{\C(\U,\A)/\U}$ on the quotient $\C(\U,\A)/\U$.

\begin{lemma}
\label{lem:SameQuotientCat}
Let $(\U,\sharp,\can,Q)$ be an exact form category, and let $\A\subset \U$ be a fully exact subcategory closed under the duality functor.
Assume that $\A$ and $\U$ are idempotent complete.
Then the inclusion $\A \subset \U$ induces an equivalence of exact form categories
$$\left(\frac{\C(\A,\A)}{\A},\sharp,\can,Q_{\SS\A}\right)  \stackrel{\simeq}{\longrightarrow} \left(\frac{\C(\U,\A)}{\U},\sharp,\can,Q_{\SS\U}\right).$$ 
\end{lemma}

\begin{proof}
By \cite[Lemma 9.4]{myHermKex}, the functor $\C(\A,\A)/\A \to \C(\U,\A)/\U$ is a duality preserving equivalence of exact categories.
Thus we have show that it induces isomorphisms of abelian groups of forms, that is, for $A\in \C(\A,\A)$ we have to show that the map
$$Q_{\SS\A}(A) \to Q_{\SS\U}(A)$$
is an isomorphism.
This is the map
$$\underset{B \stackrel{\sim}{\rightarrowtail} A}{\colim}\ Q_{\C(\A,\A)}(B) \longrightarrow \underset{U \stackrel{\sim}{\rightarrowtail} A}{\colim}\ Q_{\C(\U,A)}(U)$$
where the first colimit is taken over those admissible monomorphisms $B \rightarrowtail A$ in $\C(\A,\A)$ with cokernel isomorphic to an object in $\A$ and the second colimit is taken over those admissible monomorphisms $U \rightarrowtail A$ in $\C(\U,\A)$ with cokernel isomorphic to an object in $\U$.
Since both indexing categories are left filtering posets, it suffices to show that for every object in the right indexing category, there is an object in the left indexing category that maps to it.

So, let $U \rightarrowtail A$ be an object of the right indexing category, that is, an admissible monomorphism in $\C(\U,\A)$ with quotient a constant diagram $A/U=cV$ with in $V \in \U$.
The quotient map $A \twoheadrightarrow cV$ factors as $A \twoheadrightarrow cA^n \to cV$ for some $n \in \N$ where the first map is an admissible epimorphism in $\C(\A,\A)$ \cite[Proof of Lemma 9.3]{myHermKex}.
Since $\U$ is idempotent complete, the map $A^n \to V$ is an admissible epimorphism in $\U$ \cite[Lemma 8.5 (g)]{myHermKex}.
Let $B$ be the kernel of $A \twoheadrightarrow cA^n$ in $\C(\A,\A)$ and let $W$ be the kernel of $A^n \twoheadrightarrow V$ in $\U$.
Then $B$ is an admissible subobject of $U$ in $\C(\U,\A)$ such that $U/B = cW \in \U$,
and the admissible monomorphism $B \rightarrowtail A$ in $\C(\A,\A)$ has quotient $cA^n \in \A$.
In other words, $B \rightarrowtail A$ is an object of the left indexing category and a subobject of $U \rightarrowtail A$.
\end{proof}

\begin{lemma}
\label{lem:CAAisFlasque}
Let $(\A,\sharp,\can,Q)$ be an exact form category with strong duality.
Then there is a non-singular exact form functor $(T,\ffi_q,\ffi): \C\A \to \C\A$ and an
isometry of form functors $$id \perp (T,\ffi_q,\ffi) \cong (T,\ffi_q,\ffi).$$
In particular, the $K$-theory and Grothendieck-Witt spaces of $\C\A$ are contractible.
\end{lemma}

\begin{proof}
In \cite[Lemma 9.5]{myHermKex} we have constructed an exact form functor $(T,\ffi):\C\A \to \C\A$, $\ffi: T\sharp \cong \sharp T$ and an isometry $1 \perp (T,\ffi) \cong (T,\ffi)$ of form functors of underlying exact categories with duality.
Here we need to extend $(T,\ffi)$ to a form functor $(T,\ffi_q,\ffi)$ between form categories
such that the isometry above defines an isometry $1 \perp (T,\ffi_q,\ffi) \cong (T,\ffi_q,\ffi)$ respecting quadratic forms.

Recall \cite[Lemma 9.5]{myHermKex} the duality preserving functor $[-1]:\C_0\A \to \C_0\A$ defined by $(A[-1])_i = A_{i-1}$, $(A[-1])^i=A^{i-1}$, $i\geq 1$, and $(A[-1])_0=(A[-1])^0=0$ with the same structure maps as $A$ (shifted by $1$) and zero elsewhere.
On quadratic forms, we define $[-1]$ by 
$$Q(A) \to Q(A[-1]): (\xi,\ffi) \mapsto (\xi[-1],\ffi[-1])$$
where $(\xi[-1])_i=\xi_{i-1}$, $(\ffi[-1])_i = \ffi_{i-1}$ and zero elsewhere, $i\geq 1$.
Note that for all $\xi \in Q_{\C_0}(A)$ we have 
\begin{equation}
\label{eqn:ConstrOfTfiiqEtc}
e_{[1]}^{\bullet}(\xi) = (e^{[1]})^{\bullet}(\xi[-1]) \in Q_{\C_0} (A[-1]^{[1]})
\end{equation}
where
$$\xymatrix{
A & (A[-1])^{[1]} \ar[l]_{\hspace{-3ex}e_{[1]}}  \ar[r]^{\hspace{2ex}e^{[1]}} & A[-1].}$$

The form functor $(T,\ffi_q,\ffi):\C_0\A \to \C_0\A$ is the orthogonal sum 
$$(T,\ffi_q^0,\ffi) = \perp_{r\geq 0}[-1]^r$$
where $[-1]^r$ is the $r$-th iteration of $[-1]$.
For instance, 
$$(TA)_i=\bigoplus_{0\leq r \leq i}A_{r},\hspace{3ex} (TA)^i=\bigoplus_{0\leq r \leq i}A^{r},\hspace{3ex}  (\ffi_q^0(\xi))_i = \perp_{0\leq r \leq i}\xi_r$$
with the understanding that $A_i$, $A^i$, $\xi_i$ are zero when $i<0$.
By definition, we have an isometry $1 \perp T[-1] \cong T$.

In \cite[Lemma 9.5]{myHermKex} we proved that $(T,\ffi)$ defines a form functor $(C\A,\sharp,\can) \to (\C\A,\sharp,\can)$ using a natural transformation
$\alpha_A:(TA)^{[1]} \to T(A^{[1]})$ making the following diagram commute
$$\xymatrix{TA \hspace{2ex} &&  \ar[ll]_{T(e^{[1]})} \hspace{2ex}
T(A^{[1]}) & & \\
     & (TA)^{[1]} \ar[ul]^{e^{[1]}} \ar[ur]^{\alpha_A} & &
      \hspace{4ex} (T(A^{[1]}))^{[1]}. \ar[ll]^{(T(e^{[1]}))^{[1]}} \ar[ul]_{e^{[1]}}
}$$
We set $\alpha_A(0)=1_{TA}$, $\alpha_A(1) = \alpha_A$ and define recursively the natural transformation 
$$\alpha_A(n) = \alpha_{A^{[n-1]}} \circ (\alpha_A(n-1))^{[1]}:(TA)^{[n]} \to T(A^{[n]}).$$
The diagram then shows that $\alpha(i) \circ e = (Te) \circ \alpha(i+1)$, and we obtain the map
of colimits
$$\xymatrix{
\colim_nQ_{\C_0}(A^{[n]}) \ar[r]^{\ffi^0_q} & \colim Q_{\C_0}T(A^{[n]}) \ar[rr]^{\colim(\alpha(n)^{\bullet})} &&\colim Q_{\C_0}((TA)^{[n]})}$$
which defines $\ffi_q:Q_{\C}(A) \to Q_{\C}(TA)$ and completes the definition of the form functor $(T,\ffi_q,\ffi):\C\A \to \C\A$.
In view of the equality (\ref{eqn:ConstrOfTfiiqEtc}), we have an isometry of form functors $e^{[1]}\circ (e_{[1]})^{-1}: 1 \to [-1]$ which turns the isometry $1 \perp T[-1] \cong T$ into the isometry $1 \perp T \cong T$ of form functors $\C\A \to \C\A$.
\end{proof}

\begin{theorem}
\label{thm:coneSeq}
Let $(\U,\sharp,\can,Q)$ be an exact form category with strong duality and let $\A \subset \U$ be a fully exact subcategory closed under the duality.
Consider $\A$ as an exact form category by restriction of $\sharp$ and $Q$ to $\A$.
Assume that $\A$ and $\U$ are idempotent complete.
Then the commutative square 
$$\xymatrix{{\A} \ar[r] \ar[d] & {\U} \ar[d]\\
{\C(\A,\A)} \ar[r] & {\C(\U,\A)}}$$
of exact form categories and fully faithful exact form functors
induces a homotopy cartesian square of 
Grothendieck-Witt spaces with contractible lower left corner.
\end{theorem}

\begin{proof}
By Theorem \ref{thm:cor:filt_locn}, 
taking vertical quotients yields a map of homotopy fibrations of Grothendieck-Witt spaces which is an equivalence on base spaces (Lemma \ref{lem:SameQuotientCat}).
This proves that the square in the theorem  is homotopy cartesian after taking Grothendieck-Witt spaces.
The lower left corner is contractible, by Lemma \ref{lem:CAAisFlasque}.
\end{proof}

\part{Form categories with weak equivalences}
\label{part:FormCatw}

\section{Form categories with weak equivalences}

Recall \cite{myMV} that an {\it exact category with weak equivalences} is a pair $(\E,w)$ 
with $\E$ an exact category and $w \subset \Mor\E$ a set of
morphisms, called weak equivalences, which contains all identity morphisms, is
closed under isomorphisms, retracts, push-outs along admissible monomorphisms, pull-backs
along admissible epimorphisms, composition and the $2$ out of three property for
composition (if $2$ of the $3$ maps among $a$, $b$, $ab$ are in $w$ then so is
the third).
A weak equivalence is usually depicted as $\stackrel{_{\sim}}{\to}$ in diagrams.
A functor $F:\A \to \B$ between exact categories with weak equivalences 
$(\A,w)$ and $(\B,w)$ is called {\em exact} if it preserves admissible exact sequences and weak equivalences.
For an exact category with weak equivalences $(\E,w)$, we will write $w\E$ for
the subcategory of $\E$ whose arrows are the weak equivalences.
An object $E$ in $\E$ is called {\em acyclic} if the map $0 \to E$ is a weak equivalence.
We denote by $\E^w\subset \E$ the full subcategory of acyclic objects.
By the following lemma, the subcategory is closed under extensions in $\E$ and thus inherits the structure of an exact category from $\E$ in which a sequence is admissible exact in $\E^w$ if and only if it is in $\E$.
For a refinement of the following Lemma; see Lemma \ref{lem:WeakEqVsAcyclicCone}.

\begin{lemma}
\label{lem:WeakEqsAndExts}
Let $(\E,w)$ be an exact category with weak equivalences.
Then the following hold.
An admissible epimorphism in $\E$ is a weak equivalence if and only if its kernel is acyclic.
An admissible monomorphism in $\E$ is a weak equivalence if and only if its cokernel is acyclic.
If in a map of admissible exact sequences in $\E$
\begin{equation}
\label{eqn:MapOfExSeq}
\xymatrix{
X_{-1} \xymono[r] \ar[d]^{f_{-1}} & X_{0} \xyepi[r] \ar[d]^{f_{0}} & X_{1} \ar[d]^{f_{1}}\\
Y_{-1} \xymono[r]  & Y_{0} \xyepi[r]  & Y_{1} 
}
\end{equation}
the maps $f_{-1}$ and $f_1$ are weak equivalences then so is $f_0$.
\end{lemma}

\begin{proof}
Consider composible weak equivalences $f:X \stackrel{\sim}{\to} Y$, $g:Y \stackrel{\sim}{\twoheadrightarrow} Z$ with $g$ an admissible epimorphism.
Let $i: A \rightarrowtail Y$ be the kernel of $g$.
Then $gf$ is a weak equivalence, and we have the following pull-back square
$$\xymatrix{
A\oplus X \xyepi[d]_{(0 \ 1)} \ar[r]_{\sim}^{(i\ f)} & Y \xyepi[d]^g \\
X \ar[r]^{\sim}_{gf} & Z}$$
The canonical inclusion $X \to A\oplus X$ is a weak equivalence since its composition with the weak equivalence $(i\ f)$ is the weak equivalence $f$.
As a retract of $X \to A \oplus X$, the map $0 \to A$ is a weak equivalence.
Applied to $f=id$, this shows that $\ker(g)$ is weakly equivalent to $0$ if $g$ is an admissible epimorphism which is a weak equivalence.
Conversely, if $\ker(g)$ is weakly equivalent to $0$ and if $g$ is an admissible epimorphism then $g$ is the push-out of the weak equivalence $A \to 0$ along the admissible monomorphism $A \to Y$.

The second statement is similar and we omit the details.

For the last statement, we simply remark that $f_0$ factors as $X_0 \to P \to Y_0$ where 
$X_0 \to P$ is the push-out of the upper left corner of the diagram and $P \to Y_0$ is the pull-back of the lower right corner. 
As both maps are weak equivalences, so is their composition $f_0$.
\end{proof}

\begin{definition}
An {\em exact form category with weak equivalences} is a tuple 
$$(\E,w,\sharp,\can,Q)$$
 where $(\E,w,\sharp,\can)$ is an exact category with weak equivalences and duality \cite[Section 2.3]{myMV} and the quadruple $(\E,\sharp,\can,Q)$ is an exact form category \cite[Definition 2.21]{myForm1}. 
In particular, $\sharp: \E^{op} \to \E$ is exact and preserves weak equivalences, $\can_X:X \to X^{\sharp\sharp}$ is a natural weak equivalence, and $Q$ is quadratic left exact, that is, sends admissible exact sequences to left exact sequences as in Definition \ref{dfn:ExFormCat} (\ref{item3:dfn:ExFormCat}).
\end{definition}

\begin{example}
Any exact form category $(\E,\sharp,\can,Q)$ with strong duality defines an exact form category with weak equivalences $(\E,i,\sharp,\can,Q)$ with set of weak equivalences the isomorphisms in $\E$.
\end{example}

Recall from \cite[Definition 2.7]{myForm1} the category $\Quad(\E,Q)$ of quadratic forms in ($\E,w,\sharp,\can,Q)$.
Its objects are pairs $(X,\xi)$ where $X\in \E$ is an object of $\E$ and $\xi\in Q(X)$ is a form on $X$.
Arrows $f:(X,\xi) \to (Y,\psi)$ are those maps $f:X \to Y$ in $\E$ for which $\xi=f^{\bullet}(\psi)$.
Composition of arrows is composition in $\E$.

\begin{definition}[Quadratic spaces]
\label{dfn:QspaceWeEq}
A {\em quadratic space} in an exact form category with weak equivalences $(\E,w,\sharp,\can,Q)$ is a form $(X,\xi)$ in $(\E,Q)$  such that the associated symmetric form $\rho(\xi):X \to X^{\sharp}$ is a weak equivalence.
We denote by 
$$w\Quad(\E,Q)$$ 
the category of quadratic spaces in $(\E,w,\sharp,\can,Q)$ and arrows those maps in $\Quad(\E,Q)$ which are weak equivalences in $(\E,w)$.
\end{definition}

\begin{definition}
A form functor $(F,\ffi_q,\ffi): (\A,w,*,\can,Q_{\A}) \to (\B,w,*,\can,Q_{\B})$ of exact form categories with weak equivalences is called {\em non-singular} if the duality compatibility map $\ffi:F* \to *F$ is a weak equivalence.
We will denote by 
$$(\Fun_{ex}(\A,w;\B,w), \sharp, \can, Q)$$
the fully exact form subcategory of $(\Fun_{ex}(\A,\B), \sharp, \can, Q)$ of Lemma \ref{lem:FunExIsFormCat} consisting of those functors that preserve weak equivalences.
Note that an exact form functor is non-singular if and only if it defines a quadratic space in the exact form category with weak equivalences $(\Fun_{ex}(\A,w;\B,w), w,\sharp \can Q)$ where a map of exact functors $F\to G :(\A,w) \to (\B,w)$ is a weak equivalence if $F(A) \to G(A)$ is a weak equivalence for all objects $A$ of $\A$.
\end{definition}

Let  $(\E,w)$ be an exact category with weak equivalences.
Recall from \cite[Section 1.3]{wald:spaces} the simplicial exact category with weak equivalences $n\mapsto (S_n\E,w)$ and the
$K$-theory space $K(\E,w)=\Omega |wS_{\bullet}\E|$ of $(\E,w)$.
If $(\E,w,\sharp,\can,Q)$ is an exact form category with weak equivalences, the category
$S_n\E$ is equipped with the structure of exact category with weak equivalences and duality $(S_n\E,w,\sharp,\can)$ \cite[Section 2.6]{myMV} and a compatible structure of exact 
form category $(S_n\E,\sharp,\can,Q)$ (\S \ref{sec:FormQ} and \cite[Section 6]{myForm1}) making $(S_n\E,w,\sharp,\can,Q)$ into an exact form category with weak equivalences. 
The edge-wise subdivision $n\mapsto  \RR_n\E = S_{[n]^{op}[n]}\E   = S_{2n+1}\E$
defines the simplicial exact form category with weak equivalences
$$n\mapsto (\RR_n\E,w,\sharp,\can,Q).$$
The associated categories of quadratic spaces define a simplicial category 
$$n \mapsto w\Quad(\RR_{n}\E,Q).$$

\begin{definition}[Grothendieck-Witt space]
\label{dfn:GWSpaceII}
Let $(\E,w,\sharp,\can,Q)$ be an exact form category with weak equivalences.
The composition $w\Quad(\RR_{\bullet}\E,Q) \to w\RR_{\bullet}\E \to wS_{\bullet}\E$
of simplicial categories, in which 
the first arrow is the forgetful functor $(X,\xi) \mapsto X$, 
and the second is the canonical map $X^e_{\bullet} \to X_{\bullet}$ of
simplicial objects induced by $[n] \subset [n]^{op}[n]:i \mapsto i$
yields a map of simplicial categories
\begin{equation}
\label{eqn:RtoS}
w\Quad(\RR_{\bullet}\E,Q) \to wS_{\bullet}\E
\end{equation}
and thus a map on classifying spaces
\begin{equation}
\label{eqn:GWdfn}
|w\Quad(\RR_{\bullet}\E,Q)| \to |wS_{\bullet}\E|
\end{equation}
whose homotopy fibre (with respect to
a zero object  of $\E$ as base point of $wS_{\bullet}\E$) is defined to be the
{\it  Grothendieck-Witt space} 
$$GW(\E,w,Q)$$ of $(\E,w,\sharp,\can,Q)$.
Thus, we have a homotopy fibration
$$GW(\E,w,Q) \to |w\Quad(\RR_{\bullet}\E,Q)| \to |wS_{\bullet}\E|.$$
Inclusion of degree-zero simplicies defines a map of simplicial categories
$$w\Quad(\E,Q) = w\Quad(\RR_0\E,Q) \to w\Quad(\RR_{\bullet}\E,Q)$$
where the left simplicial category is simplicially constant.
The composition
\begin{equation}
\label{eqn:ERES}
w\Quad(\E,Q) \to w\Quad(\RR_{\bullet}\E,Q) \to wS_{\bullet}\E
\end{equation}
factors through $wS_0\E=0$ und we obtain 
a universal map 
\begin{equation}
\label{eqn:EGWunivMap}
|w\Quad(\E,Q)| \to GW(\E,w,Q)
\end{equation}
from the category of quadratic spaces to its Grothendieck-Witt space.
By the main theorem of \cite{myForm1}, this map is a group completion if $w$ is the set of isomorphisms and if every admissible exact sequence in $\E$ splits.

Finally, we define the
{\em higher Grothendieck-Witt groups} of $(\E,w,\sharp,\can,Q)$ as the homotopy
groups $$GW_i(\E,w,Q) = \pi_i GW(\E,w,Q), \phantom{1234}
i\geq 0.$$ 
\end{definition}

We give a presentation of the group $GW_0(\E,w,Q)$ in Theorem \ref{thm:Presentation} below.
For the next lemma, recall from \cite[Definition 2.7]{myForm1} the definition of a natural transformation between form functors, that is, a map between quadratic spaces in the exact form category with weak equivalences $(\Fun_{ex}(\A,w;\B,w), w,\sharp \can Q)$.

\begin{lemma}
\label{lem:Lemma2MV}
Let $f:(F,\ffi_q,\ffi) \stackrel{\sim}{\to} (G,\psi_q,\psi)$ be a natural weak equivalence between non-singular exact form functors $$(F,\ffi_q,\ffi), (G,\psi_q,\psi): (\A,w,\sharp,\can,Q_{\A}) \to (\B,w,\sharp,\can,Q_{\B})$$ of exact form categories with weak equivalences.
Then on Grothendieck-Witt spaces, $(F,\ffi_q,\ffi)$ and $(G,\psi_q,\psi)$ induce homotopic maps 
$$GW(F,\ffi_q,\ffi) \sim GW(G,\psi_q,\psi): GW(\A,w,Q_{\A}) \to GW(\B,w,Q_{\B}).$$
\end{lemma}

\begin{proof}
The natural weak equivalence $f$ defines natural transformations of the two functors after applying the $w\Quad(\RR_n,Q)$- and $wS_{\bullet}$-constructions and thus functors $[1]\times w\Quad(\RR_n\A,\Q_{\A}) \to  w\Quad(\RR_n\B,\Q_{\B})$ and $[1]\times wS_{\bullet}\A \to wS_{\bullet}\B$ compatible with the simplicial structure and forgetful maps.
This  defines homotopies between the $w\Quad(\RR_{\bullet},Q)$. $wS_{\bullet}$-constructions and leads upon taking homotopy fibres
to the homotopy $I \times GW(\A,w,Q_{\A}) \to GW(\B,w,Q_{\B})$ between the maps induced by $(F,\ffi_q,\ffi)$ and $(G,\psi_q,\psi)$.
\end{proof}

\begin{definition}[Hyperbolic form category]
\label{dfn:HypCatFunEtc}
Let $(\E,w)$ be an exact category with weak equivalences.
Recall from \cite[Section 2.9]{myMV} the hyperbolic exact category with weak equivalences and duality $(\H\E,w)$.
As exact category with weak equivalences it is $(\E\times \E^{op},w\times w^{op})$, and it has a strict duality $(X,Y)^{\sharp}=(Y,X)$, $(f,g)^{\sharp}=(g,f)$.
We consider $\H\E$ as an exact form category with functor of quadratic forms the symmetric forms.
Thus, $Q(X,Y) = \E(X,Y)$ and the functorial Mackey-functor of forms is
$$\xymatrix{
\H\E[(X,Y),(Y,X)] \ar[rr]^{\hspace{5ex}(f,g)\mapsto f+g} && Q(X,Y) \ar[rr]^{\hspace{-6ex}f \mapsto (f,f)} && \H\E[(X,Y),(Y,X)].
}$$
In \cite[Proposition 1]{myMV} we constructed a homotopy equivalence 
\begin{equation}
\label{eqn:GWHEisKE}
GW(\H\E,w) \stackrel{\sim}{\longrightarrow} K(\E,w)
\end{equation}
which makes the following diagram commute
$$\xymatrix{
w\Quad(\H\E) \ar[r] \ar[d]^{\wr} & GW(\H\E,w) \ar[d]^{\wr}\\
w\E \ar[r] & K(\E,w)
}$$
where the horizontal maps are the universal maps and the left vertical map is the weak equivalence $(X,Y,f) \mapsto X$ of \cite[Lemma 3]{myMV}.

Let $(E,w,\sharp,\can,Q)$ be an exact form category with weak equivalences.
The {\em hyperbolic} and {\em forgetful} form functors of \S \ref{sec:AddtyExForm} preserve weak equivalences and define the sequence of non-singular exact form functors
$$(\H\E,w) \stackrel{H}{\longrightarrow} (\E,w,\sharp,\can,Q) \stackrel{F}{\longrightarrow} (\H\E,w).$$
In view of the homotopy equivalence (\ref{eqn:GWHEisKE}), we obtain hyperbolic and forgetful maps on Grothendieck-Witt spaces
$$K(\E,w) \simeq GW(\H\E,w) \stackrel{H}{\longrightarrow} GW(\E,w,Q) \stackrel{F}{\longrightarrow} GW(\H\E,w) \simeq K(\E,w).$$
\end{definition}

\section{Additivity Theorems and Presentation of $GW_0$}

Let $(\E,w,\sharp,\can,Q)$ be an exact form category with weak equivalences.
Recall that the poset $[n]$ has a unique strict duality $i \mapsto n-i$ which gives the category $S_n\A$ the structure of an exact form category $(S_n\A,\sharp,\can,Q)$; see \S \ref{sec:FormQ} and \cite[Section 6]{myForm1}.
For any $i\leq j$, the map $[0] \to \Ar[n]: 0 \mapsto (i\leq j)$ induces by restriction
the exact evaluation functor $\ev_{i,j}: (S_n\E,w) \to (\E,w): X \mapsto X_{ij}$ evaluating at $i\leq j$.
The map $[0] \to \Ar[n]: 0 \mapsto (i\leq j)$ preserves dualities if $j=n-i$ in which case restriction along $[0] \to \Ar[n]$ defines an non-singular exact form functor $e_{i,j}:(S_n\E,w,Q)\to (\E,w,Q): X \mapsto X_{ij}$.
Recall from Definition \ref{dfn:HypCatFunEtc} the forgetful map $F:GW(\E,w,Q) \to K(\E,w)$.

\begin{theorem}[Additivity]
\label{thm:AddtyForWeakExFunctors}
Let $(\E,w,\sharp,\can,Q)$ be an exact form category with weak equivalences.
Then the form functor
$$(\ev_{12},\ev_{01}\circ F): (S_3\E,w,Q) \longrightarrow (\E,w,Q) \times (\H\E,w)$$
induces a weak equivalence of spaces
$$GW(S_3\E,w,Q) \longrightarrow GW(\E,w,Q) \times K(\E,w).$$
\end{theorem}

\begin{proof}
In view of Theorem \ref{thm:exAddty},
the proof is the same as \cite[Proof of Theorem 4] {myMV} replacing Lemmas 4 and 5 of 
{\em loc.\ cit.\ }by Lemmas \ref{lem:strictification} and \ref{lem:simplResLem} below.
\end{proof}

Denote by $\FormCatW$ the category of exact form categories with weak equivalences where maps are the non-singular exact form functors. 
Composition is composition of form functors.
Denote by $\FormCatW^{str}$ the subcategory of exact form categories with weak equivalences that have a strict duality where maps are the non-singular exact form functors that commute with dualities and have the identity as duality compatibility maps.
Composition is composition of form functors.
The inclusion defines the functor
$$\lax:\FormCatW^{str} \to \FormCatW.$$

\begin{lemma}[Strictification lemma]
\label{lem:strictification}
There are a (strictification) functor 
 $$\str: \FormCatW \to \FormCatW^{str}: (\A,w,*,\eta,Q) \mapsto (\A^{\str}_w,w,\sharp,id,Q^{\str})$$
and natural transformations
$(\Sigma,\sigma_q,\sigma):id \to \lax\circ\str$ and $(\Lambda,\lambda_q,\Lambda):\lax\circ\str
\to id$ such that the 
compositions $(\Sigma,\sigma_q,\sigma)\circ (\Lambda,\lambda_q,\lambda)$ and
$(\Lambda,\lambda_q,\lambda)\circ (\Sigma,\sigma_q,\sigma)$ are weakly equivalent to the
identity form functor.  
In particular, for any exact form category with weak equivalences 
$(\E,w,*,\can,Q)$, we have a homotopy equivalence of Grothendieck-Witt spaces
$$GW(\Sigma,\sigma_q,\sigma):GW(\E,w,*,\can,Q) \stackrel{\sim}{\longrightarrow}
GW(\E^{\str}_w,w,\sharp,id,Q^{str}).$$ 
\end{lemma}

\begin{proof}
The category with duality and weak equivalences $(\A^{\str}_w,w,\sharp,id)$ was defined in \cite[Proof of Lemma 4]{myMV}. 
We need to equip it with a quadratic functor $Q^{\str}:(\A^{\str}_w)^{op} \to \Ab$.
Recall from {\em loc.\ cit.\ }that the objects of $\A^{\str}_w$ are triples $(A,B,f)$ with $A$, $B$ objects of
$\A$ and  $f:A \stackrel{\sim}{\to} B^*$ a weak equivalence in $\A$.   
A morphism from $(A_0, B_0,f_0)$ to $(A_1, B_1,f_1)$ 
  is a pair $(a,b)$ of morphisms $a:A_0 \to A_1$ and $b:B_1 \to B_0$ in $\A$
  such that $f_1 a = b^* f_0$.
The dual of $(A,B,f)$ is $(B,A,f^*\can)$.
We set 
$$Q^{\str}(A,B,f) = \{(\xi,a) \in Q(A)\times \A(A,B)|\ f^* \can\circ a = a^* f = \rho(\xi)\}$$
and define transfer and restriction as
$$\tau_{\A^{\str}}(a,b) =( \tau_{\A}(b^*f),a+b),\hspace{3ex}\rho_{\A^{\str}}(\xi,a) = a.$$
One checks that this makes $\A^{\str}$ into an exact form category with weak equivalences.
On the level of exact categories with duality and weak equivalences, the natural transformations
$(\Sigma,\sigma)$ and $(\Lambda,\lambda)$ where defined in {\em loc.\ cit.\ }.
On quadratic forms we define them as
$$\renewcommand\arraystretch{1.5}
\begin{array}{cccccccc}
\sigma_q:&Q(A) &\to & Q^{\str}(A,A^*,\can)&: & \xi & \mapsto & (\xi,\rho(\xi))\\
\lambda_q:&Q^{\str}(A,B,f) &\to & Q(A)&: & (\xi,a) &\mapsto& \xi.\end{array}$$
The rest is as in \cite[Proof of Lemma 4]{myMV}.
\end{proof}

Let $(\E,w,\sharp,\can,Q)$ be an exact form category with weak equivalences.
Recall from \S \ref{sec:FormQ}, or \cite[Lemma 6.2]{myForm1}, that for a poset $\P$ with strict duality, the category of functors $\Fun(\P,\E)$ is canonically an exact form category.
A sequence of functors $\P \to \E$ is called exact if it is admissibly exact in $\E$ when evaluated at each object of $\P$.
This makes $\Fun(\P,\E)$ into an exact form category with weak equivalences.

As in \cite{myMV}, we let $\Fun_w(\P,E) \subset \Fun(\P,E)$ be the full subcategory of those functors $F:\P \to \E$ which send every arrow in $\P$ to a weak equivalence in $\E$.
By restriction of structure, $\Fun_w(\P,\E)$ is an exact form category with weak equivalences.
In particular, for $n\in \N$, we obtain the exact form category with weak equivalences 
$\Fun_w(\underline{n},E)$ where $\underline{n}=[n]^{op}[n]=\{n'< \cdots 0'<0<1 < \cdots <n\}$.
For a functor 
$X \in \Fun_w(\underline{n},E)$ the abelian group $Q(X)$ of quadratic forms on $X$ consists of those pairs $(\xi,\ffi)$ where $\xi\in Q(X_{0'})$ is a quadratic form on $X_{0'}\in \E$ and $\ffi:X \to X^*$ is a symmetric form such that $\rho(\xi): X_{0'} \to (X_{0'})^*$ is the diagonal arrow in the commutative diagram
$$\xymatrix{
X_{0'} \ar[rr]^{X(0'\to 0)} \ar[d]_{\ffi_{0'}} && X_0 \ar[d]^{\ffi_0}\\
(X_0)^* \ar[rr]_{X(0' \to 0)^*} && (X_{0'})^*.
}
$$
Varying $n\in \N$, we obtain a simplicial exact form category with weak equivalences $n \mapsto \Fun_w(\underline{n},E)$ where $\underline{n}=[n]^{op}[n]$.

\begin{lemma}[Simplicial Resolution Lemma]
\label{lem:simplResLem}
Let $(\E,w,*,\eta,Q)$ be an exact form category with weak equivalences and strong
duality. 
Then there are homotopy equivalences
$$
\renewcommand\arraystretch{1.5} 
\begin{array}{rcl}
|w\Quad(\RR_{\bullet}\E,Q)| & \simeq & |n \mapsto
i\Quad (\RR_{\bullet}\Fun_w(\underline{n},\E))|,\\
GW(\E,w,Q) & \simeq & | n \mapsto GW(\Fun_w(\underline{n},\E),i,Q)|
\end{array}$$
which are functorial for exact form functors $(F,\ffi_q,\ffi)$ 
for which $\ffi$ is an isomorphism.
\end{lemma}

\begin{proof}
For a category $\C$, denote by $i\C$ the subcategory of isomorphisms in $\C$.
We have an equivalence of categories
\begin{equation}
\label{eqn:ResLemEqn}
i\Fun([n]^{op},\ w\Quad(\E,Q)) \stackrel{\sim}{\longrightarrow} i\Quad(\Fun_w(\underline{n},\E),Q)
\end{equation}
which sends an object 
$$(X_n,\xi_n) \stackrel{f_n}{\to} (X_{n-1},\xi_{n-1})
\stackrel{f_{n-1}}{\to} \cdots \stackrel{f_{1}}{\to}(X_0,\xi_0)$$ of the left
hand category to the object 
$$X_n \stackrel{f_n}{\to} X_{n-1}
\stackrel{f_{n-1}}{\to} \cdots \stackrel{f_{1}}{\to} X_0 \stackrel{\rho(\xi_0)}{\longrightarrow}
X_0^* \stackrel{f_{1}^*}{\to} X_{1}^* \stackrel{f_{2}^*}{\to} \cdots
\stackrel{f_{n}^*}{\to} X_n^*$$
equipped with the form $(\xi_0,\ffi)$ where $\ffi=(\eta_{X_n},...,\eta_{X_0},1,...,1)$.
A map $(g_n,...,g_0)$ (which is an isomorphism) is sent
to the map $(g_n,...,g_0,(g_0^*)^{-1},...,(g_n^*)^{-1}$).
Together with the weak equivalence
$$|\C| \stackrel{\sim}{\longrightarrow} | n \mapsto i\Fun([n]^{op},\C)|$$
\cite[Proof of Lemma 5, Equation (6)]{myMV} for any small category $\C$, this yields the horizontal homotopy equivalences in the commutative diagram
$$\xymatrix{
|w\Quad(\RR_{\bullet}\E,Q)| \ar[r]^{\hspace{-7ex} \sim} \ar[d]_F & |n \mapsto
i\Quad(\RR_{\bullet}\Fun_w(\underline{n},\E),Q)| \ar[d]^F\\
|w\Quad(\RR_{\bullet}\H\E) \ar[r]^{\hspace{-8ex} \sim} & |n \mapsto
i\Quad \RR_{\bullet}\H\Fun_w(\underline{n},\E)|.
}$$
The rest is as in \cite[Lemma 5]{myMV}.
\end{proof}

Recall that a quadratic space in an exact form category with weak equivalences $(\E,w,\sharp,\can,Q)$ is a pair $(X,\xi)$ of an object $X$ of $\E$ equipped with a form $\xi \in Q(X)$ such that $\rho(\xi):X \to X^{\sharp}$ is a weak equivalence.

\begin{theorem}[Presentation of $GW_0$]
\label{thm:Presentation}
Let $(\E,w,\sharp,\can,Q)$ be an exact form category with weak equivalences.
Then the abelian group 
$$GW_0(\E,w,Q)=\pi_0GW(\E,w,Q)$$ has the following presentation.
Generators are the classes $[X,\xi]$ of quadratic spaces $(X,\xi)$ in $(\E,w,Q)$.
The relations are as follows.
\begin{enumerate}
\item
\label{cor:itm0:Kh_0}
$[X \oplus Y,\xi \perp \zeta]  = [X,\xi] +[Y,\zeta] $
\item
\label{cor:itm1:Kh_0}
If $g:X \to Y$ is a weak equivalence, then $[Y,\xi]=[X,g^{\bullet}(\xi)]$.
\item
\label{cor:itm2:Kh_0}
If $(E_{\bullet},\xi,\ffi)$ is a quadratic space in the category
of exact sequences $S_2\E$, that is, $\ffi$ is a map 
$$\xymatrix{
\hspace{2ex}E_{\bullet} \ar[d]^{\ffi}_{\wr}:  & E_{-1}
\xymono[r]^i\ar[d]^{\ffi_{-1}}_{\wr} & E_0 \xyepi[r]^p\ar[d]^{\ffi_{0}}_{\wr}
&  E_1 \ar[d]^{\ffi_{1}}_{\wr} \\
\hspace{2ex}E_{\bullet}^{\sharp}: & E_1^{\sharp} \xymono[r]_{p^{\sharp}} & E_0^{\sharp} \xyepi[r]_{i^{\sharp}} &
E_{-1}^{\sharp}}$$ 
of admissible exact sequences in $\E$ with
$(\ffi_{-1},\ffi_0,\ffi_1)=(\ffi_1^{\sharp},\ffi_0^{\sharp},\ffi_{-1}^{\sharp})\circ \can $ a symmetric
weak equivalence, $\xi \in Q(E_0)$ a quadratic form on $E_0$ satisfying $i^{\bullet}(\xi)=0$ and $\rho(\xi)=\ffi_0$
 then $$[E_0,\xi] = [H(E_{-1})].
$$
\end{enumerate}
\end{theorem}

\begin{proof}
This is {\em mutatis mutandis} the same as \cite[Proposition 3]{myMV} reducing the proof of the theorem to the presentation of $GW_0$ for exact form categories with strong duality (and weak equivalences the isomorphisms) of \cite[Theorem 6.5]{myForm1} by means of the Stricification and Simplicial Resolution Lemmas \ref{lem:strictification} and \ref{lem:simplResLem}.
\end{proof}

We finish the section with a consequence of Additivity for the homotopy type of the $n$-th iteration of  the $\RR_{\bullet}$-construction.
For an exact form category with weak equivalences  $\E=(\E,w,\sharp,\can,Q)$,
we write $\RR^{(n)}_{\bullet\cdots \bullet}\E$ for the $n$-th iterate of the $\RR_{\bullet}$-construction.
This is the $n$-simplicial exact form category with weak equivalences 
$$\RR^{(n)}_{k_1,...,k_n}\E = \RR_{k_1}\RR_{k_2}\dots \RR_{k_n}\E,\vspace{2ex}$$
which is the full exact form subcategory with weak equivalences of
$$\Fun(\Ar(\underline{k_1})\times \Ar(\underline{k_2})\times \cdots \times \Ar(\underline{k_n}), \E)$$
consisting of those functors $A: \Ar(\underline{k_1})\times \Ar(\underline{k_2})\times \cdots \times \Ar(\underline{k_n}) \longrightarrow \E$
such that for all $r,s=1,...,n$ and all $i_r \leq j_r \in \underline{k_r}$, $r\neq s$, the functor
$$A_{(i_1,j_1), \dots (i_{s-1},j_{s-1}), \bullet, (i_{s+1},j_{s+1}),\dots (i_{n}, j_{n}}): \Ar(\underline{k_s})\to \E$$
is an object of $\RR_{k_s}\E$.
Weak equivalences and admissible exact sequences are those maps and sequences of diagrams which are weak equivalences and admissible exact sequences when evaluated at every object of $\Ar(\underline{k_1})\times \cdots \times \Ar(\underline{k_n})$.
We denote by $\RR^{(n)}_{\bullet}\E$ the diagonal of the multi-simplicial exact form category with weak equivalences $\RR^{(n)}_{\bullet\cdots \bullet}\E$.

\begin{corollary}
\label{cor:wAdd}
Let $(\E,w,\sharp,\can,Q)$ be an exact form category with weak equivalences.
Then for all $n\geq 1$ the form functor
$$(\ev_{0'0},\prod_{r=0}^{n-1}\ev_{r,r+1}\circ F): (\RR_n\E,w,Q) \longrightarrow (\E,w,Q) \times (\H\E,w)^n$$
induces a weak equivalence of spaces after application of the functor $|w\Quad\RR^{(n)}_{\bullet}|$.
\end{corollary}

\begin{proof}
This immediately follows from Theorem \ref{thm:AddtyForWeakExFunctors} and is {\em mutatis mutandis} the same as \cite[Proposition A.8]{myJPAA}.
\end{proof}

If we replace $\E$ with $\RR_{\bullet}^{(n)}\E$ in the sequence (\ref{eqn:ERES}) we obtain the sequence of pointed topological spaces
\begin{equation}
\label{eqn:AddtyFib1}
|w\Quad\RR_{\bullet}^{(n)}(\E,Q)| \longrightarrow |w\Quad\RR_{\bullet}\RR_{\bullet}^{(n)}(\E,Q)| \longrightarrow |wS_{\bullet}\RR_{\bullet}^{(n)}\E|
\end{equation}
with trivial composition.

\begin{proposition}
Let $\E=(\E,w,\sharp,\can,Q)$ be an exact form category with weak equivalences.
Then for all $n\geq 1$, the sequence
(\ref{eqn:AddtyFib1})
of pointed topological spaces
is a homotopy fibration.
\end{proposition}

\begin{proof}
The proof is the same as that of \cite[Proposition 4.1]{myJPAA} replacing the reference to \cite[\S 3 Remark 7]{myMV} with Corollary \ref{cor:wAdd}.
\end{proof}

Composing the second map in (\ref{eqn:AddtyFib1}) with the homotopy equivalence 
$$|wS_{\bullet}\RR_{\bullet}^{(n)}\E| \stackrel{\sim}{\to} |wS_{\bullet}S_{\bullet}^{(n)}\E|$$ of \cite[\S 2.4 Lemma 1]{myMV}, we obtain the following.

\begin{corollary}
\label{cor:AddtyExFib1}
Let $\E=(\E,w,\sharp,\can,Q)$ be a an exact form category with weak equivalences.
Then for every integer $n\geq 1$, the sequence
$$|w\Quad(\RR_{\bullet}^{(n)}\E,Q)| \longrightarrow |w\Quad(\RR_{\bullet}^{(1+n)}\E,Q)| \longrightarrow |wS_{\bullet}^{(1+n)}\E|$$
is a homotopy fibration of pointed topological spaces.
\end{corollary}

\section{The Fibration Theorem}

For an exact category with weak equivalences $(\E,w)$, 
the full subcategory $\E^w$ of $w$-acyclic objects is closed under extensions in $\E$, and thus inherits an
exact structure from $\E$ such that the inclusion $\E^w \subset \E$ preserves and detects admissible
exact sequences. 

We recall from \cite[Section 4]{myMV} the definition of a symmetric cone.

\begin{definition}
\label{dfn:symCone}
Let $(\E,w,*,\can)$ be an exact category with weak equivalences and duality.
A {\em symmetric cone} on $(\E,w,*,\can)$ is given by the following data:
\begin{enumerate}
\item
exact functors $P:\E \to \E^w$, and $C:\E \to \E^w$,
\item
a natural admissible epimorphism $p_E: PE \twoheadrightarrow E$ and a natural admissible monomorphism
$i_E:E \rightarrowtail CE$,
\item
a natural transformation $\gamma:P* \to *C$,
 $\gamma_E: P(E^*) \to (CE)^*$, such that
$i^*_E\circ \gamma_E=p_{E^*}$ for all objects $E$ of $\E$.
\end{enumerate}
The symmetric cone is called {\em strong} if $\gamma$ is a natural isomorphism.
\end{definition}
It is convenient to define $\tilde{\gamma}_E:C(E^*) \to (PE)^*$ by 
$\tilde{\gamma}_E = P(\can_E)^*\circ \gamma_{E^*}^*\circ \can_{C(E^*)}$.
Then we have
\begin{equation}
\renewcommand\arraystretch{2}
\label{eqn:gammaTildeGamma}
\begin{array}{ll}
i^*_E\gamma_E=p_{E^*}, & \tilde{\gamma}_E = P(\can_E)^*\circ \gamma_{E^*}^*\circ \can_{C(E^*)},\\
p_E^*=\tilde{\gamma}_E\circ i_{E^*}, & \gamma_E=C(\can_E)^*\circ
\tilde{\gamma}_{E^*}^* \circ \can_{P(E^*)},\\
\gamma_{E^*} \circ P(\can_E) = \tilde{\gamma}_{E}^*\circ \can_{PE}, &
\tilde{\gamma}_{E^*}\circ C(\can_E) = \gamma_E^*\circ \can_{CE}.
\end{array}
\end{equation}

\begin{lemma}
\label{lem:WeakEqVsAcyclicCone}
Let $(\E,w,*,\can)$ be an exact category with weak equivalences and duality which has a symmetric cone $(C,P)$.
If in a map of admissible exact sequences (\ref{eqn:MapOfExSeq}),
two of the maps $f_{-1},f_0, f_1$ are weak equivalences then so is the third.
In particular, an arrow $f:X \to Y$ in $\E$ is a weak equivalence if and only if its cone
$$C(f)=\colim\left(
\xymatrix{
C(X) & X \xymono[l]_{\hspace{3ex}i_X} \ar[r]^f & Y}\right)
$$
is acyclic.
Similarily, 
$f:X \to Y$ in $\E$ is a weak equivalence if and only if its path object
$$P(f)=\lim\left(
\xymatrix{
 X \ar[r]^f & Y &\xyepi[l]_{p_Y} P(Y) }\right)
$$
is acyclic.
\end{lemma}

\begin{proof}
The case where $f_{-1}$ and $f_1$ are weak equivalences is Lemma \ref{lem:WeakEqsAndExts}.
So, assume $f_{-1}$ and $f_0$ are weak equivalences.
Using the same factorisation of $f_0$ as in Lemma \ref{lem:WeakEqsAndExts} we can assume that $f_{-1}$ is the identity and thus, $f_0$ is the pull-back of $f_1$.
Factoring $f_1$ into a weak equivalence followed by an admissible epimorphism $X_1 \to X_1 \oplus P(Y_1) \to Y_1$ and taking the pull-back of that composition along $Y_0 \twoheadrightarrow Y_1$ allows us to assume that $f_1$ is an admissible epimorphism.
Then $f_0$ is an admissible epimorphism which is a weak equivalence.
By Lemma \ref{lem:WeakEqsAndExts} $\ker(f_0)$ is acyclic.
Since $\ker(f_0) =\ker(f_1)$, the map $f_1$ is a weak equivalence again by Lemma \ref{lem:WeakEqsAndExts}.
The case when $f_0$ and $f_1$ are weak equivalences is similar, and we omit the details.

The statement about the cone $C(f)$ of $f$ follows from the above applied to the following map of admissible exact sequences in which the middle vertical arrow is a weak equivalence
$$\xymatrix{
X \xymono[r] \ar[d]^{f} & Y \oplus C(X) \xyepi[r] \ar[d]^{(1\ 0)} & C(f) \ar[d]\\
Y \xymono[r]_1  & Y \xyepi[r]  & 0.
}$$
Similarily for the path object $P(f)$.
\end{proof}

The main result of this section is the following.

\begin{theorem}[Fibration]
\label{thm:ChgOfWkEq}
Let $(\E,w,\sharp,\can,Q)$ be an exact form category with weak equivalences.
Assume that the underlying  exact category with weak equivalences and duality $(\E,w,\sharp,\can)$ has a strong symmetric cone.
Let $v$ be another set of weak equivalences in $\E$ containing $w$ and which
is closed under the duality.
Consider $(\E^v,w)$, $(\E^v,v)$, $(\E,v)$ as exact form categories with weak equivalences by restricting  $(\sharp,\can,Q)$.
Then the commutative square of duality preserving full inclusions of exact form categories
\begin{equation}
\label{eqn:thm:ChgOfWkEq}
\xymatrix{
(\E^v,w,Q) \ar[r] \ar[d] & (\E^v,v,Q) \ar[d]\\
(\E,w,Q)  \ar[r] & (\E,v,Q)}
\end{equation}
induces a 
homotopy cartesian square of associated
Grothendieck-Witt spaces.
Moreover, the upper right
corner has contractible Grothendieck-Witt space. 
\end{theorem}

\begin{remark}
\label{rmk:strongCone}
Theorem \ref{thm:ChgOfWkEq} was proved for symmetric forms in \cite[Theorem 6]{myMV}.
In {\em loc. cit. }we did not assume the symmetric cone to be strong but this is used in the proof as we will see below.
This does not affect any applications of \cite[Theorem 6]{myMV} in \cite{myMV} as all exact categories with weak equivalences and duality
the theorem was applied to had a strong symmetric cone.
\end{remark} 

We will reduce the proof of the Theorem \ref{thm:ChgOfWkEq} to idempotent complete exact form categories.
Recall that $\tilde{\E}$ denotes the idempotent completion of an exact category $\E$; see Example\ref{ex:IdempCompl}. 
If $(\E,w)$ is an exact category with weak equivalences, then $\tilde{\E}$ is equipped with the set of weak equivalences which are precisely the retracts of weak equivalences in $\E$.
We need the following.

\begin{lemma}
\label{lem:EwCofinalI}
Let $(\E,w,\sharp,\can,Q)$ be an exact form category with weak equivalences and strong
duality which has a strong symmetric cone.
Then the
commutative diagram of exact form categories (the weak equivalences being the isomorphisms)
$$
\xymatrix{(\E^w,i,Q) \ar[r] \ar[d] & (\tilde{\E}^w,i,Q) \ar[d] \\
(\E,i,Q) \ar[r] & (\tilde{\E},i,Q)}
$$
induces a homotopy cartesian square of Grothendieck-Witt spaces.
\end{lemma}

\begin{proof}
Using the Cofinality Theorem \ref{thm:CofinalityExCats} and Lemma \ref{lem:beforeExCof}, the proof is now the same as in \cite[Lemma 6]{myMV}.
\end{proof}

Recall that the category $\Mor\E= \Fun([1],\E)$ of morphisms in an exact form category $(\E,*,\can,Q)$, 
is equipped with the quadratic functor $Q$ which at the object $g:X \to Y$ of $\Mor\E$ consists of all pairs $(\xi,a)$ where $\xi\in Q(X)$ and $a:X \to Y^*$ such that $ag^*=a^*\can_Y  g = \rho(\xi)$.
By restriction, this makes $\Mor_w\E=\Fun_w([1],\E)$ into an exact form category.
As in \cite[Section 4]{myMV} the following proposition is the key to proving
Theorem \ref{thm:ChgOfWkEq}.

\begin{proposition}
\label{prop:E=Morw}
Let $(\E,w,\ast,\can,Q)$ be an exact form category with weak equivalences and strong
duality. 
Assume that $(\E,w,\ast,\can)$ has a strong symmetric cone.
Then the commutative diagram 
\begin{equation}
\label{eqn:prop:E=Morw}
\xymatrix{({\E}^w,Q) \ar[r] \ar[d] &  ({\Mor_w}{\E}^w,Q) \ar[d]\\
           ({\E},Q) \ar[r]^{\hspace{-4ex}I}_{\hspace{-2ex}E \mapsto 1_E} & ({\Mor}_w{\E},Q)}
\end{equation}
of full  inclusions of exact form categories with strong duality (all weak
equivalences being isomorphisms)
induces a homotopy cartesian square of associated Grothendieck-Witt spaces.
\end{proposition}

\begin{proof}
By Lemma \ref{lem:EwCofinalI}, we can assume that all categories in diagram (\ref{eqn:prop:E=Morw}) are idempotent complete.
We will write $\C_0(\E,w)$ and $\C(\E,w)$ for $\C_0(\E,\E^w)$ and $\C(\E,\E^w)$.
As in \cite[Proof of Proposition 4]{myMV}, 
we will construct a commutative diagram of
exact form categories with strong duality and non-singular exact form functors
\begin{equation}
\label{eqn:pf:E=Morw}
\xymatrix{{\E}^w \ar[r] \ar[d] &  {\Mor_w}{\E}^w \ar[d] \ar[rr] && {\C}({\E}^w)
  \ar[d]\ar[r] & {\C}({\Mor_w}{\E}^w) \ar[d]\\
   {\E} \ar[r]^{\hspace{-1.5ex}I} & {\Mor}_w{\E} \ar[rr]^{(F,\ffi_q,\ffi)} && {\C}({\E},w)
   \ar[r]^{\hspace{-3ex}\C(I)} & {\C}({\Mor}_w{\E},w)} 
\end{equation}
where all form categories are equipped with isomorphisms as weak equivalences and all quadratic functors are induced by the given quadratic functor $Q:\E^{op} \to \Ab$ as explained before the statement of Proposition \ref{prop:E=Morw} and  in \S \ref{sec:ConeConstr}.
The left square is (\ref{eqn:prop:E=Morw}), and the right square is obtained from (\ref{eqn:prop:E=Morw}) by
functoriality of the cone.
We will show:
\vspace{1ex}

\begin{itemize}
\item[(\dag)] ($GW$ applied to) the compositions
$(F,\ffi_q,\ffi)\circ I$ and $\C(I)\circ(F,\ffi_q,\ffi)$ 
are homotopic to the constant diagram inclusions $\Gamma_0:\E \to \C(\E,w)$ and $\Gamma: \Mor_w\E \to C(\Mor_w\E,w)$, denoted $c$ in diagram (\ref{eqn:ConstDiag}),  
in such a way that the homotopies restricted to $\E^w$
and $\Mor_w\E^w$ have images in the Grothendieck-Witt spaces of $\C(\E^w)$ and
$\C(\Mor_w\E^w)$, respectively.
\end{itemize}
As in \cite[Proof of Proposition 4]{myMV}, this implies Proposition \ref{prop:E=Morw}.
Indeed, by Theorem \ref{thm:coneSeq}, diagram (\ref{eqn:pf:E=Morw}) is the inclusion of fibre into total space of a sequence of fibrations with base the sequence
$$\C(\E,w) \stackrel{\C(I)}{\longrightarrow} \C(\Mor_w\E,w) \stackrel{(F,\ffi_q,\ffi)}{\longrightarrow} \C(\E,w) \stackrel{\C(I)}{\longrightarrow} \C(\Mor_w\E,w).$$
By (\dag), the composition of the first two and the composition of the last two maps are homotopic to the identity (after applying $GW$).
This shows that $GW$ applied to these maps are isomorphisms in the homotopy category of spaces.
In particular, all squares in (\ref{eqn:pf:E=Morw}) are homotopy cartesian after applying $GW$.

We now construct diagram (\ref{eqn:pf:E=Morw}).
The middle lower horizontal form functor is actually the composition of two form functors 
$$\xymatrix{{\Mor}_w{\E} \ar[rr]^{(F,\ffi_q,\ffi)} && {\C_0}({\E},w) \ar[r] & {\C}({\E},w)}$$
where the second is the localization form functor.
The functor $F$ and the duality compatibility map $\ffi$ were defined in {\em loc.\ cit.} 
Recall that $F$ sends an object 
$g:X \to Y$ of $\Mor_w\E$ to the object $F(g)$ of of $\C_0(\E,w)$ given by the diagram $U_{\bullet}(g) \to U^{\bullet}(g)$,
$$\xymatrix{
X\xymono[r] \ar[d]^g & X \oplus PY \ar[d]^{\left(\begin{smallmatrix}g&
      p \\i & 0\end{smallmatrix}\right)} \xymono[r] 
& X \oplus PY \oplus C X 
\ar[d]^{\left(\begin{smallmatrix}g&p & 0 \\i & 0 & 1\\
      0&1&0\end{smallmatrix}\right)} \xymono[r] 
&X \oplus PY \oplus C X \oplus PY 
\ar[d]
      \xymono[r] & \cdots\\
 Y & Y\oplus CX \xyepi[l]
& Y\oplus CX \oplus PY \xyepi[l] & Y\oplus CX \oplus PY \oplus CX \xyepi[l]
& \cdots \xyepi[l]
}$$
which in degree $i = 0,1,2,3,\dots$ is
$$\xymatrix{
U_i\ar@{}[r]|-{=} \ar[d]& X \ar[d]_g \xymono[rd] \ar@{}[r]|-{\oplus}& PY   \ar@{}[r]|-{\oplus} \xyepi[dl]|!{[l];[d]}\hole \ar[dr] & CX\ar@{}[r]|-{\oplus} \ar[dl]|!{[l];[d]}\hole \ar[dr] & PY\ar@{}[r]|-{\oplus} \ar[dl]|!{[l];[d]}\hole \ar[rd] & CX \ar[dl]|!{[l];[d]}\hole  \ar@{}[r]|-{\cdots} & i+1\text{ summands}\\
U^i\ar@{}[r]|-{=}  & Y \ar@{}[r]|-{\oplus}& CX\ar@{}[r]|-{\oplus} & PY\ar@{}[r]|-{\oplus} & CX\ar@{}[r]|-{\oplus} & PY \ar@{}[r]|-{\cdots} &  i+1\text{ summands}
}$$
with matrix representation
$$\begin{pmatrix}
g & p & 0 & \cdots & 0 & 0\\
i  & 0 & 1 & 0 & \vdots &0\\
0 & 1 & 0 & \ddots & 0& \vdots \\
\vdots & 0 & \ddots & \ddots &1 &0\\
0 &\vdots  & 0 & 1 & 0 & 1\\
0 &0 & \cdots  & 0 & 1 & 0
\end{pmatrix},
$$
that is, the only non-zero partial maps are as indicated with solid arrows in the diagram where $X \rightarrowtail CX$, $PY \twoheadrightarrow Y$ are the natural maps $i_X$, $p_Y$, and $CX \to CX$, $PY \to PY$ are identity maps.
The maps $U_i \to U_{i+1}$ and $U^{i+1} \to U^i$ are the canonical inclusions into the first $i+1$ summands and the canonical projections onto the first $i+1$ summands.
They have cokernel and kernel in $\E^w$, respectively.
Recall that $g$ is a weak equivalence.
Therefore, the map $U_i \to U^{i+1}$ is an admissible monomorphism with cokernel in $\E^w$
since its matrix representation (left matrix of (\ref{eqn:UUmono})) can be transformed to the right matrix in (\ref{eqn:UUmono}) by elementary row operation which shows that its cokernel is isomorphic to the cone $C(g)$ which is in $\E^w$, by Lemma \ref{lem:WeakEqVsAcyclicCone}:
\begin{equation}
\label{eqn:UUmono}
\begin{pmatrix}
g & p & 0 & \cdots & 0 \\
i  & 0 & 1 & 0 & \vdots\\
0 & 1 & 0 & \ddots & 0 \\
\vdots & 0 & \ddots & \ddots &1\\
0 &\vdots  & 0 & 1 & 0 \\
0 &0 & \cdots  & 0 & 1 
\end{pmatrix},
\hspace{10ex}
\begin{pmatrix}
g & 0 & 0 & \cdots & 0 \\
i  & 0 & 0 & 0 & \vdots\\
0 & 1 & 0 & \ddots & 0 \\
\vdots & 0 & \ddots & \ddots &0\\
0 &\vdots  & 0 & 1 & 0 \\
0 &0 & \cdots  & 0 & 1 
\end{pmatrix}.
\end{equation}
Similarily, the map $U^{i+1} \to U^i$ is an admissible epimorphism with kernel in $\E^w$ because
its matrix representation (left matrix of (\ref{eqn:UUepi})) can be transformed to the right matrix in (\ref{eqn:UUepi}) by elementary column operation which shows that its kernel is isomorphic to the path object $P(g)$ which is in $\E^w$, by Lemma \ref{lem:WeakEqVsAcyclicCone}:
\begin{equation}
\label{eqn:UUepi}
\begin{pmatrix}
g & p & 0 & \cdots & 0 & 0\\
i  & 0 & 1 & 0 & \vdots &0\\
0 & 1 & 0 & \ddots & 0& \vdots \\
\vdots & 0 & \ddots & \ddots &1 &0\\
0 &\vdots  & 0 & 1 & 0 & 1
\end{pmatrix},
\hspace{10ex}
\begin{pmatrix}
g & p & 0 & \cdots & 0 & 0\\
0  & 0 & 1 & 0 & \vdots &0\\
0 & 0 & 0 & \ddots & 0& \vdots \\
\vdots & 0 & \ddots & \ddots &1 &0\\
0 &\vdots  & 0 & 0 & 0 & 1
\end{pmatrix}.
\end{equation}
Hence $F(g) \in \C_0(\E,w)$.

A map $(a,b): g_1 \to g_2$ in $\Mor_w(\E)$ is a pair of maps $a:X_1 \to X_2$, $b:Y_1 \to Y_2$ in $\E$ such that $g_2a=bg_1$.
Under the functor $F$ it is sent to the map given by the diagonal matrices
$(a,Pb,Ca,Pb,Ca,\dots): U_i(g_1) \to U_i(g_2)$ and 
$(b,Ca,Pb,Ca,Pb,\dots): U^i(g_1) \to U^i(g_2)$.
The duality compatibility map $\ffi: F(g^*) \to F(g)^*$ at $g\in \Mor_w\E$ is given by the diagonal matrices 
$$\ffi_i = (1,\gamma,\tilde{\gamma},\gamma,\tilde{\gamma},\dots):U_i(g^*) \stackrel{\cong}{\longrightarrow} (U^i(g))^*$$
$$\ffi^i = (1,\tilde{\gamma},\gamma,\tilde{\gamma},\gamma, \dots):U^i(g^*) \stackrel{\cong}{\longrightarrow} (U_i(g))^*$$
which are natural isomorphisms\footnote{Here is where we need $\gamma$ and $\tilde{\gamma}$ to be isomorphisms; see Remark \ref{rmk:strongCone}.} since $\gamma$ and $\tilde{\gamma}$ are.
In particular, a symmetric map $(a,\bar{a}): g \to g^*$ in $\Mor_w\E$, given by $a:X \to Y^*$ with $\bar{a}=a^*\can_Y:Y \to X^*$ such that $\bar{a}g = g^*a$,
is sent to the symmetric map $\alpha = \ffi \circ F(a,\bar{a})$ in $\C_0(\E,w)$ given by the diagonal matrices
$$\alpha_i = (a,\gamma_XP(\bar{a}),\tilde{\gamma}_YC(a),\gamma_XP(\bar{a}),\dots):U_i(g) \to (U^i(g))^*,$$
$$\alpha^i = (\bar{a},\tilde{\gamma}_YC(a),\gamma_XP(\bar{a}),\tilde{\gamma}_YC(a),\dots):U^i(g) \to (U_i(g))^*.$$
On quadratic forms, we define the form functor 
$$\ffi_q:Q(g) \longrightarrow Q(F(g)):(\xi,a) \mapsto (\xi_{\bullet},\alpha)$$
as follows.
Recall that an element of $Q(g)$ is a pair $(\xi,a)$ with $\xi\in Q(X)$ and $(a,\bar{a}=a^*\can):g \to g^*$ a symmetric map in $\Mor_w\E$ such that $\rho(\xi)=g^*a=\bar{a}g$.
It is sent to the element $(\xi_{\bullet},\alpha) \in Q(F(g))$ where $\alpha:F(g) \to F(g)^*$ was defined above and
$\xi_{\bullet}=(\xi_i)_{i\geq 0}$ is the compatible family of quadratic forms $\xi_i \in Q(U_i)$ defined by
$$
\xi_i= \left(
\renewcommand\arraystretch{1.1}
\begin{array}{l} X \oplus PY \oplus CX \oplus PY \cdots \\ \downarrow 1  \\ X \end{array}\right)^{\bullet}\xi +
\tau\left(
\renewcommand\arraystretch{1.5}
\begin{array}{l} X\hspace{1ex}  \oplus   \hspace{2ex}PY \hspace{1.5ex} \oplus  \hspace{1ex}CX  \hspace{2ex}\oplus \hspace{1.5ex} PY \cdots \\  \hspace{4ex}\swarrow \hspace{-.5ex} \text{\tiny$\bar{a}p$}  \hspace{4ex}\swarrow  \hspace{-.5ex} \text{\tiny$\tilde{\gamma}Ca$} \hspace{4ex} \swarrow  \hspace{-.5ex}  \text{\tiny$\gamma P\bar{a}$} \\ X^*  \oplus  (PY)^*  \oplus  (CX)^*  \oplus  (PY)^* \cdots  \end{array}
\right).
$$
For the entry of $\tau$ we obtain, using (\ref{eqn:gammaTildeGamma}), the lower triangular matrix with first off diagonal entries
 $$(\bar{a}p_Y,\tilde{\gamma}_YC(a),\gamma_X P(\bar{a}),\dots)^*\can = (\tilde{\gamma}_YC(a)i_X, \gamma_XP(\bar{a}), \tilde{\gamma}_YC(a),\dots).$$
 It follows that $\rho(\xi_i) = \alpha^i\circ (U_i\to U^i)$, and $(\xi_{\bullet},\alpha)$ is indeed an element of $Q(F(g))$.
This finishes the definition of the non-singular exact form functor $(F,\ffi_q,\ffi)$.

It remains to show ($\dag$), that is, we have to show that the natural transformations defined in \cite[Proof of Proposition 4]{myMV} respect quadratic forms.
Consider the natural transformations of functors $\Mor_w\E \to \C_0(\E,w)$ sending the object $g:X \to Y$ of $\Mor_w\E$ to 
$$X \stackrel{j}{\longrightarrow} F(g) \stackrel{q}{\longrightarrow} Y$$
where $X$ and $Y$ are considered as constant diagrams and the map $j$ is the map from the initial object of the diagram $F(g)$ to $F(g)$ and the map $q$ is the map from $F(g)$ to its terminal object.

When restricted to $\E \subset \Mor_w\E$, the natural transformation $j$ respects quadratic forms and thus is a natural transformation
$j: \Gamma_0 \to (F,\ffi_q,\ffi)\circ I$ of non-singular exact form functors, that is, a morphism of quadratic spaces in $\Fun_{ex}(\E,\C_0(\E,w))$.
By  \cite[Lemma 2.29]{myForm1}, this defines an orthogonal sum decomposition of non-singular exact form functors
$$(F,\ffi_q,\ffi)\circ I = \Gamma_0 \perp \Gamma_0^{\perp}.$$
As noted in {\em loc.\ cit.,\ }$\Gamma_0^{\perp}$ has image in $\C_0(\E^w)$.
By Lemma \ref{lem:CAAisFlasque}, the Grothendieck-Witt space of $\C(\E^w)$ is contractible.
It follows that the non-singular exact form functors $(F,\ffi_q,\ffi)\circ I$ and $\Gamma_0: \E \to \C(\E,w)$ induce homotopic maps on Grothendieck-Witt spaces.
When restricted to $\E^w$, both form functors and the natural transformation have images in $\C(\E^w)$.
This proves ($\dag$) for the first composition.

Write $\A_0$ and $\A$ for the exact form category of exact functors from $\Mor_w\E$ to $\C_0(\Mor_w\E,w)$ and $\C(\Mor_w\E,w)$, respectively.
Localisation defines the canonical exact form functor $\A_0 \to \A$.
As in {\em loc.\ cit.,\ }we will construct a quadratic space 
\begin{equation}
\label{eqn:SomeSqInFibThm}
\xymatrix{
j \xymono[r] \xyepi[d]&
C(I)\circ F \xyepi[d]\\
\Gamma \xymono[r] &
q}
\end{equation}
in $(S_3\A,Q)$
whose evaluation at the upper right corner $(03)$ is $C(I)\circ (F,\ffi_q,\ffi)$, evaluation at the lower left corner $(12)$ is the constant diagram form functor $\Gamma:\Mor_w\E \to \C(\Mor_w\E,w)$ and its evaluation $\ker(j \twoheadrightarrow \Gamma)$ at $(01)$ has image in $\C(\Mor_w\E^w)$.
By the Additivity Theorem \ref{thm:AddtyForExFunctors} (\ref{prop:AddtyForFunctors:2}) and the contractibility of the Grothendieck-Witt space of $\C(\Mor_w\E^w)$, the following non-singular exact form functors induce homotopic maps after application of $GW$
$$C(I)\circ (F,\ffi_q,\ffi) \sim \Gamma \perp H(\ker(j \twoheadrightarrow \Gamma)) \sim \Gamma.$$
This will prove ($\dag$) for the second composition.

As a symmetric space in $S_3\A$, diagram  (\ref{eqn:SomeSqInFibThm}) was constructed in {\em loc.\ cit.}
As a quadratic space in $(S_3\A,Q)$, diagram (\ref{eqn:SomeSqInFibThm}) is the image under the canonical map 
$\A_0\to \A$ 
of a quadratic space in $\Fun([1]\times [1],\A_0)$ which under the identification $\C_0(\Mor_w\E,w) = \Mor_w\C_0(\E,w)$ sends $(g:X \to Y) \in \Mor_w\E$ to 
$$
 \xymatrix{ 
( X \stackrel{j}{\to} F(g) ) \ar[r]^{(j,1)} \ar[d]_{(1,q)} &
(F(g) \stackrel{1}{\to} F(g) ) \ar[d]^{(1,q)}\\
( X \stackrel{g}{\to}  Y ) \ar[r]^{(j,1)} &
(F(g) \stackrel{q}{\to}  Y) 
}$$
equipped with the symmetric form of {\em loc.\ cit.}
We have to check that the quadratic forms on $(F(g) \stackrel{1}{\to} F(g) )$ and 
$( X \stackrel{g}{\to}  Y )$ considered as functors in $g$ agree when restricted to $( X \stackrel{j}{\to} F(g) )$.
The first restriction yields the form functor which on quadratic forms is the composition
$$\renewcommand\arraystretch{2}
\begin{array}{ccccccc}
Q\left(X \stackrel{g}{\to} Y\right) & \stackrel{\ffi_q}{\longrightarrow} & Q(F(g)) & \stackrel{I}{\longrightarrow} & Q\left(F(g) \stackrel{1}{\to} F(g)\right) & \stackrel{(j,1)^{\bullet}}{\longrightarrow}& Q\left(X \stackrel{j}{\to} F(g)\right)\\
\xi,a &  \mapsto & (\xi_{\bullet},\alpha) & \mapsto  & (\xi_{\bullet},\alpha), (\alpha,\alpha)  & \mapsto & \xi, (\alpha j, j^*\alpha)
\end{array}
$$
where the second entry in each column of the lower row always refers to the symmetric form.
The second restriction yields the form functor which on quadratic forms is 
$$\renewcommand\arraystretch{2}
\begin{array}{ccc}
Q\left(X \stackrel{g}{\to} Y\right) & \stackrel{(1,q)^{\bullet}}{\longrightarrow} & Q\left(X \stackrel{g}{\to} F(g)\right)\\
\xi,a  & \mapsto & \xi, (q^*a,\bar{a}q).
\end{array}
$$
We have to check that these two maps are equal.
Since $(\alpha j, j^*\alpha)$ and $(q^*a,\bar{a}q)$ are symmetric maps, it suffices to check the equality $\alpha j = q^* a$ of maps $X \to F(g)^*$.
Since $X \in \C_0(\E,w)$ is a constant diagram, any map $X \to F(g)^*$ in $\C_0(\E,w)$ factors through the initial object $Y^*$ of $F(g)^*$, and the equality $\alpha j = q^* a$ is equivalent to the equality $(\alpha j)_0 = (q^* a)_0$.
The last equality holds because $(\alpha j)_0 = \alpha_0\circ  j_0 = a \circ 1$
and $ (q^* a)_0 =  (q^*)_0 \circ a_0=1 \circ a$.
\end{proof}

Next we prove a variant of Theorem \ref{thm:ChgOfWkEq}.

\begin{proposition}
\label{prop:ChgOfWkEqIsoVariant}
Let $(\E,w,\sharp,\can,Q)$ be an exact form category with weak equivalences and strong duality.
Assume that the underlying  exact category with weak equivalences and duality $(\E,w,\sharp,\can)$ has a strong symmetric cone. 
Then the commutative square of duality preserving inclusions of exact form categories with weak equivalences 
$$\xymatrix{
(\E^w,i,Q) \ar[r] \ar[d] & (\E^w,w,Q) \ar[d]\\
(\E,i,Q)  \ar[r] & (\E,w,Q)}
$$
induces a 
homotopy cartesian square of associated
Grothendieck-Witt spaces where the weak equivalences of the left two form categories are the isomorphisms.
Moreover, the upper right
corner has contractible Grothendieck-Witt space. 
\end{proposition}

\begin{proof}
Replacing \cite[Proposition 4]{myMV} with Proposition \ref{prop:E=Morw}, the proof is {\em mutatis mutandis} the same as that of \cite[Proposition 5]{myMV}.
Let $\Fun^1_w(\underline{n},E) \subset \Fun_w(\underline{n},E)$ be the full form subcategory of those objects $A: \underline{n} \to \E$ for which $A_p \to A_q$ is an admissible monomorphism and $A_{q'} \to A_{p'}$ is an admissible epimorphism for $0 \leq p \leq q \leq n$.
The main point is to show that this form functor inclusion induces an equivalence of Grothendieck-Witt spaces.
The inverse is provided by the difference $GW(F)-GW(G)$ for two form functors $F,G:\Fun_w(\underline{n},E) \to \Fun^1_w(\underline{n},E)$ defined in \cite[Proof of Proposition 5]{myMV}.
We need to say what $F$ and $G$ do on quadratic forms.
Note that a quadratic form on $A: \underline{n} \to \E$ is a pair $(q,\beta)$ where $q\in Q(E_{0'})$ and $\beta:A \to A^{\sharp}$ such that $\rho(q) = \beta_{0}\circ A_{0'\leq 0}$ and $\beta = \beta^{\sharp}\can$.
As in \cite{myMV} we illustrate the case $n=1$.
Recall that $F$ sends an object $E_{1'} \stackrel{\sim}{\rightarrow} E_{0'} \stackrel{\sim}{\rightarrow} E_{0} \stackrel{\sim} {\rightarrow} E_{1}$ to $E_{1'} \oplus P(E_{0'})\stackrel{\sim}{\twoheadrightarrow} E_{0'} \stackrel{\sim}{\rightarrow} E_{0} \stackrel{\sim} {\rightarrowtail} E_{1} \oplus CE_0$, and $G$ sends that object to $PE_{0'} \stackrel{\sim}{\twoheadrightarrow} 0 \stackrel{\sim}{\rightarrow} 0 \stackrel{\sim} {\rightarrowtail} E_{1} CE_0$.
On quadratic forms, the form functors $F$ and $G$ are defined by 
$$\renewcommand\arraystretch{1.5}
\begin{array}{ccc}
Q\left(E_{1'} \stackrel{a}{\to} E_{0'} \stackrel{b}{\to} E_{0} \stackrel{c}{\to} E_{1}\right) &\longrightarrow&
Q\left(E_{1'} \oplus PE_{0'}\stackrel{(a,p)}{\twoheadrightarrow} E_{0'} \stackrel{b}{\to} E_{0} \stackrel{\left(\begin{smallmatrix}c\\ i\end{smallmatrix}\right)}{\rightarrowtail} E_{1}\oplus CE_0\right)\\
(q,\beta) & \mapsto & (q,\beta^{F})
\end{array}
$$
and
$$\renewcommand\arraystretch{1.5}
\begin{array}{ccc}
Q\left(E_{1'} \stackrel{a}{\to} E_{0'} \stackrel{b}{\to} E_{0} \stackrel{c}{\to} E_{1}\right) &\longrightarrow&
Q\left(PE_{0'}\twoheadrightarrow 0 \to  0 \rightarrowtail  CE_0\right)\\
(q,\beta) & \mapsto & (0,\beta^{G})
\end{array}
$$
where $\beta^F$ and $\beta^G$ are as in \cite{myMV}, that is, 
$$\beta^F_{1'} =\left( \begin{smallmatrix}\beta_{1'} & 0 \\ 0 &  \gamma_{E_0}\circ P\beta_{0'}\end{smallmatrix}\right),\hspace{2ex}
\beta^F_{0'} = \beta_{0'},\hspace{2ex} \beta^F_0=\beta_0,\hspace{2ex} 
\beta^F_{1} =\left( \begin{smallmatrix}\beta_{1} & 0 \\ 0 &  \tilde{\gamma}_{E_{0'}}\circ C\beta_{0}\end{smallmatrix}\right)$$
and
$$
\beta^G_{1'}=\gamma_{E_0}\circ P\beta_{0'},\hspace{2ex} \beta^G_1 = \tilde{\gamma}_{E_{0'}}\circ C\beta_0.
$$
\end{proof}

\begin{proof}[Proof of Theorem \ref{thm:ChgOfWkEq}]
The proof now is the same is that of \cite[Theorem 6]{myMV} replacing Lemma 2, Lemma 4 and Proposition 5 of \cite{myMV} with Lemma \ref{lem:Lemma2MV}, Lemma \ref{lem:strictification} and Proposition \ref{prop:ChgOfWkEqIsoVariant}.
\end{proof}

\begin{theorem}[Cofinality]
\label{thm:CofinilityForWeakEq}
Let $(\E,w,\sharp,\can,Q)$ be an exact form category with weak equivalences which has a strong symmetric cone.
Let $A \subset K_0(\E,w)$ be a subgroup closed under the duality action on
$K_0(\E,w)$, and let $\E_A \subset \E$ be the full 
subcategory of those objects whose class in $K_0(\E,w)$ belongs to $A$.
Then the category $\E_A$ inherits the structure of an exact form category with weak
equivalences from $(\E,w,\sharp,\can,Q)$, and
the induced map on Grothendieck Witt spaces
$$GW(\E_A,w,Q) \longrightarrow GW(\E,w,Q)$$ is an isomorphism on
$\pi_i$, $i\geq 1$, and a monomorphism 
on $\pi_0$.
\end{theorem}

\begin{proof}
The proof is the same as that of \cite[Theorem 7]{myMV} replacing Lemma 2, Lemma 4, Proposition 5 and the reference to \cite[Corollary 5.2]{myHermKex} in \cite{myMV} with Lemma \ref{lem:Lemma2MV}, Lemma \ref{lem:strictification}, Proposition \ref{prop:ChgOfWkEqIsoVariant} and Theorem \ref{thm:CofinalityExCats}.
\end{proof}

To finish the section we record the following corollary of Theorem \ref{thm:ChgOfWkEq}.

\begin{corollary}
\label{cor:RndotFib}
Assume the hypothesis of Theorem \ref{thm:ChgOfWkEq}.
Then for all $n\geq 1$, the commutative square (\ref{eqn:thm:ChgOfWkEq}) of exact form categories with weak equivalences induces a homotopy cartesian square with contractible upper right corner after application of the functor
$|w\Quad\RR^{(n)}(\phantom{E})|$.
\end{corollary}

\begin{proof}
This follows from Theorem \ref{thm:ChgOfWkEq} in view of the definition of the Grothendieck-Witt space (in case $n=1$), Corollary \ref{cor:AddtyExFib1} (for $n\geq 2$) and its $K$-theoretic analogue \cite[Theorem 1.6.4]{wald:spaces}, \cite[Proposition 1.5.3]{wald:spaces} and the fact that barycentric subdivision does not change homotopy type \cite[Lemma 1]{myMV}.
\end{proof}

\section{From exact categories to chain complexes}
\label{sec:FromExToChain}

Let $(\E,{\sharp},\can,Q)$ be an exact form category with strong duality.
We will make the category $\Ch_b(\E)$ of bounded chain complexes 
in $\E$ into an exact form category with weak equivalences $(\Ch_b\E,\quis,\sharp,\can,Q)$ such that $(\E,Q)$ and $(\Ch_b\E,\quis,Q)$ have homotopy equivalent Grothendieck-Witt spaces (Theorem \ref{prop:gilletWald}).
The structure of exact category with weak equivalences and duality on $\Ch_b\E$ is standard \cite[\S 6]{myMV}, and our task is to equip $\Ch_b(\E)$ with the correct quadratic functor; see Definition \ref{dfn:quadChain}.

In this paper, we use homological indexing.
So, an object of $\Ch_b(\E)$ is a chain complex in $\E$, that is, a diagram in $\E$
$$(E,d):\hspace{4ex} \cdots \to E_{n+1}
\stackrel{d_{n+1}}{\longrightarrow} E_n \stackrel{d_{n}}{\longrightarrow} E_{n-1} \to \cdots,$$
$n\in \Z$, such that $E_n=0$ for $|n| >> 0$ and $d_nd_{n+1}=0$ for all $n\in \Z$.
A morphism of chain complexes, also called chain map,  $f:(A,d) \to (B,d)$ is a family of maps $f_n:A_n \to B_n$ in $\E$ such that $f_{n-1}d_n=d_nf_n$ for all $n\in \Z$.
I will denote by $\E((A,d),(B,d))$, or simply $\E(A,B)$ if the differentials $d$ are understood, the abelian group of chain maps $(A,d) \to (B,d)$ in $\Ch_b\E$.
A sequence 
$(A,d) \to (B,d) \to
(C,d)$ of chain complexes in $\E$ is called {\em exact} if the sequence $A_n \rightarrowtail B_n \twoheadrightarrow
C_n$ is admissible exact in $\E$ for all $n\in \Z$.
Call a chain complex $(E,d)$ in $\E$ {\em strictly
  acyclic} if every 
differential $d_n$ is the composition $E_n \twoheadrightarrow \im d_n
\rightarrowtail E_{n-1}$  of an admissible epimorphism followed by an admissible monomorphism, and the
sequences $\im d_{n+1} \rightarrowtail E_n \twoheadrightarrow \im d_n$  are
admissible exact in $\E$.
A chain complex is called {\em acyclic} if it is homotopy equivalent to a
strictly acyclic chain complex.
A chain map $f: (A,d) \to (B,d)$  is a {\em quasi-isomorphism}  if its cone $C(f)$ is acyclic where
$C(f)_n = B_n \oplus A_{n-1}$ with differential $C(f)_n \to C(f)_{n-1}$ given by
$$\left(\begin{smallmatrix}d_n & f \\ 0 & -d_{n-1}\end{smallmatrix}\right): B_n \oplus A_{n-1} \to B_{n-1} \oplus A_{n-2}.$$
Write $\quis$ for the set of quasi-isomorphisms in $\Ch_b\E$.
Then the pair
$$(\Ch_b(\E),\quis)$$
is an exact category with weak equivalences.

The duality functor $({\sharp},\can)$ on $\E$ extends to a duality functor $({\sharp},\can)$ on $\Ch_b\E$ where
$(E,d)^{\sharp}=(E^{\sharp},d^{\sharp})$, $f^{\sharp}:(Y,d)^{\sharp} \to (X,d)^{\sharp}$ and $\can_E:(E,d) \to (E,d)^{{\sharp}{\sharp}}$ are given by
\renewcommand\arraystretch{1.5} 
$$\begin{array}{lcllcl}
(E^{\sharp})_i &=&(E_{-i})^{\sharp},&
 (f^{\sharp})_i&=&(f_{-i})^{\sharp},\\
 (d^{\sharp})_i&=& (-1)^{i+1}(d_{-i+1})^{\sharp},&
 (\can_{E})_i&=& (-1)^{i}\can_{E_i}.
\end{array}
$$

\begin{remark}
We explain the choice of signs in the formula for the duality above.
If $(\E,\otimes, [\phantom{A},\phantom{B}])$ is a closed symmetric monoidal category with evaluation $e:[A,B]\otimes A \to B$, coevaluation $\triangledown: A \to [B, A\otimes B]$ and symmetry $c:A\otimes B \to B\otimes A$.  
Then the category of bounded chain complexes $\Ch_b\E$ in $\E$ has a standard structure of a closed symmetric monoidal category where evaluation and coevaluation are induced by $e$ and $\triangledown$ without introducing signs, and the symmetry follows the Koszul sign rule.
A standard example is $\E=R\proj$, the category of finitely generated projective modules over a commutative ring $R$; see for instance \cite[\S 1.1]{myJPAA}.

An object $A$ of a closed symmetric monoidal category defines a duality $X \mapsto [X,A]$ with double dual identification 
$$X \stackrel{\triangledown}{\longrightarrow} [[X,A],X \otimes [X,A]] \stackrel{[1,c]}{\longrightarrow}  [[X,A],[X,A]\otimes X] \stackrel{[1,e]}{\longrightarrow}  [[X,A],A] .$$

With our sign choices, if the duality  $X \mapsto [X,A]_{\E}$ on $\E$ comes from a closed symmetric monoidal structure as explained above, then its extension to $\Ch_b\E$ is the duality $X \to [X,A]_{\Ch_b\E}$  coming from the  closed symmetric monoidal structure on $\Ch_b\E$ where  $A$ is considered as a chain complex concentrated in degree $0$. 
\end{remark}

Now we come to the crucial definition of quadratic forms on chain complexes.

\begin{definition}
\label{dfn:quadChain}
Let $(\E,{\sharp},\can,Q)$ be an exact form category with strong duality.
A {\em quadratic form on a chain complex} $(E,d) \in \Ch_b\E$ is a pair $(\xi,\ffi)$, where 
$\ffi:(E,d) \to (E,d)^{\sharp}$ is a map of complexes and
$\xi\in Q(E_0)$ is a quadratic form on $E_0\in \E$  satisfying 
$$\ffi^{\sharp}\can_E = \ffi,\hspace{4ex} d_1^{\bullet}(\xi) = 0,\hspace{4ex}\rho(\xi) = \ffi_0.$$
As usual, we denote by $Q(E)$ the abelian group of quadratic forms on $E=(E,d)$.
Transfer and restriction
$$\E(E,E^{\sharp}) \stackrel{\tau}{\longrightarrow} Q(E) \stackrel{\rho}{\longrightarrow} \E(E,E^{\sharp})$$
are defined by
$\tau(f) = (\tau(f_0),f+f^{\sharp}\can_E)$ and $\rho(\xi,\ffi) = \ffi$.
Note that $\tau(f)$ 
is indeed a quadratic form on the chain complex $E$ since
$$d_1^{\bullet}(\tau(f_0)) = \tau((d_1)^{\sharp}f_0d_1) = \tau(-f_{-1}d_0d_1)=0$$ and
$\rho(\tau(f_0)) = f_0+(f_0)^{\sharp}\can_{E_0}$. 
\end{definition}

\begin{lemma}
Let $(\E,{\sharp},\can,Q)$ be an exact form category with strong duality.
Then the tuple 
$$(\Ch_b\E,\quis,{\sharp},\can,Q)$$
 is an exact form category with weak equivalences and strong duality which has a strong symmetric cone.
\end{lemma}

\begin{proof}
It is standard that $(\Ch_b\E,\quis,{\sharp},\can)$ is an exact category with weak equivalences and strong duality which has a strong symmetric cone \cite[\S 6]{myMV}.
Thus, it remains to show that the functor $(Q,\tau,\rho)$ is quadratic left exact.
For the convenience of the reader we will also spell out the symmetric cone below.

It is clear that $\tau$ and $\rho$ are $C_2$-equivariant, functorial in $E\in \Ch_b\E$ and that $\rho\tau=1+\sigma$.
For chain maps $f,g:(A,d) \to (B,d)$ we have
$$\renewcommand\arraystretch{1.5}
\begin{array}{rl}
& Q(f+ g)^{\bullet}(\xi,\ffi) - Q(f)^{\bullet}(\xi,\ffi) - Q(g)^{\bullet}(\xi,\ffi)\\
 = & (\tau(f_0^{\sharp}\circ \rho(\xi)\circ g_0), f^{\sharp}\ffi g + g^{\sharp}\ffi f)\\
 = & \tau(f^{\sharp}\circ \rho(\xi,\ffi)\circ g).
\end{array}$$
Thus, $(\Ch_b\E,\sharp,\can,Q)$ is an additive form category.

We need to check that $Q:(\Ch_b\E)^{op} \to \Ab$ is quadratic left exact (\cite[Definition 2.22]{myForm1}).
Let $X \stackrel{i}{\rightarrowtail} Y \stackrel{p}{\twoheadrightarrow} Z$ be an admissible exact sequence of bounded chain complexes in $\E$.
Then $p^{\bullet}:Q(Z) \to Q(Y):(\xi,\ffi) \mapsto (p_0^{\bullet}\xi,p^{\sharp}\ffi p)$ is injective because $p_0^{\bullet}$ and $\ffi \mapsto p^{\sharp}\ffi p$ are.
Let $(\xi,\ffi)\in Q(Y)$ such that $i^{\sharp} \circ \ffi=0$ and $i_0^{\bullet}(\xi)=0$.
Then there are unique $\zeta\in Q(Z_0)$ and $\psi:Z \to Z^{\sharp}$ such that $\xi = p_0^{\bullet}(\zeta)$ and $\ffi = p^{\sharp}\psi p$ since $Q:\E^{op} \to\Ab$ and $(\Ch_b\E)^{op} \to\Ab: E \mapsto \E(E,E^{\sharp})^{C_2}$ are quadratic left exact functors.
Note that $\psi^{\sharp}\can_Z=\psi$, $d_1^{\bullet}(\zeta)=0$ and $\rho(\zeta)=\psi_0$ since
$\psi\mapsto p^{\sharp}\psi p$ is $C_2$-equivariant and injective, $p_1^{\bullet}d_1^{\bullet}(\zeta) = d_1^{\bullet}p_0^{\bullet}(\zeta) = d_1^{\bullet} (\xi) = 0$ and $p_1^{\bullet}$ is injective, and $p_0^{\sharp} \rho(\zeta)p_0 = \rho(p_0^{\bullet}\zeta) = \rho(\xi) = \ffi_0 = p_0^{\sharp} \psi_0 p_0$.
Thus, $(\zeta,\psi) \in Q(Z)$ satisfies $p^{\bullet}(\zeta,\psi) = (\xi,\ffi)$ and $Q:(\Ch_b\E)^{op} 
\to \Ab$ is quadratic left exact.

A conceptual proof of the existence of a strong symmetric cone was given in \cite[\S 7.5]{myMV}. 
Here, we spell out the formulas in detail (switching the order of tensor product compared to {\em loc.\ cit.}).
For a chain complex $E=(E,d)$, cone $CE$ and path object $PE$ are defined as usual by $(CE)_n=E_n\oplus E_{n-1}$ and $(PE)_n = E_n \oplus E_{n+1}$ with differentials
$$\left(\begin{smallmatrix}d_n & 1 \\ 0 & -d_{n-1}\end{smallmatrix}\right): (CE)_n \to (CE)_{n-1} \hspace{3ex}\text{and}\hspace{3ex}
\left(\begin{smallmatrix}d_n & 0 \\ 1 & -d_{n+1}\end{smallmatrix}\right): (PE)_n \to (PE)_{n-1}.
$$
They are equipped with the natural maps $i:E \rightarrowtail CE$ and $p:PE \twoheadrightarrow E$ which in degree $n$ are
$$\left(\begin{smallmatrix} 1 \\ 0 \end{smallmatrix}\right): E_n \to E_n \oplus E_{n-1}\hspace{4ex}\text{and}\hspace{4ex}
(1\ \ 0) : E_n \oplus E_{n+1} \to E_n.$$
Finally, the natural transformation $\gamma_E: P\sharp E \to \sharp CE$ in degree $n$ is the isomorphism
$$\left(\begin{smallmatrix}1 & 0 \\ 0 & (-1)^n \end{smallmatrix}\right): (E_{-n})^{\sharp} \oplus (E_{-n-1})^{\sharp} \stackrel{\cong}{\longrightarrow} (E_{-n})^{\sharp} \oplus (E_{-n-1})^{\sharp}.$$
\end{proof}

For an exact form category with strong duality $(\E,{\sharp},\can,Q)$, 
the inclusion $G:\E \to \Ch_b\E$ as complexes concentrated in degree $0$, defines a fully faithful 
exact form functor
\begin{equation}
\label{eqn:GilletWald}
(G,1,1): (\E,i,\sharp,\can,Q) \to (\Ch_b\E,\quis,\sharp,\can,Q).
\end{equation}

\begin{lemma}
\label{lem:GilWladSurj}
Assume that $\E$ is semi-idempotent complete.
Then the form functor (\ref{eqn:GilletWald}) induces a surjective map of Grothendieck-Witt groups
$$GW_0(\E,i,Q) \twoheadrightarrow GW_0(\Ch_b\E,\quis,Q).$$
\end{lemma}

\begin{proof}
Recall from Theorem \ref{thm:Presentation} that $GW_0(\Ch_b\E,\quis,Q)$ is generated by symbols $[E,\xi,\ffi]$ where
 $E$ is a bounded chain complex in $\E$  and $(\xi,\ffi)\in Q(E)$ is a quadratic form on $E$ with $\ffi$ a quasi-isomorphism.

Let $n\geq 2$ be an integer and let $[E,\xi,\ffi]$ be one of those generators.
Assume that $E_i=0$ for $i>n$ and $i<-n$.
Since the cone of $\ffi$ is acyclic (hence strictly acyclic, by
semi-idempotent completeness of $\E$), the map 
$\left(\begin{smallmatrix}d_{n}\\
    \phi_{n}\end{smallmatrix}\right):E_{n} \to E_{n-1}\oplus E_{-n}^{\sharp}$ 
is an admissible monomorphism.
Define a complex $\tilde{E}$ supported in $[-n,n]$ by
$$E_{n} \stackrel{\left(\begin{smallmatrix}d_{n}\\
    \ffi_{n}\end{smallmatrix}\right)}{\rightarrowtail} 
E_{n-1} \oplus E_{-n}^{\sharp} \stackrel{\left(\begin{smallmatrix}d &
    0 \end{smallmatrix}\right)}{\longrightarrow} E_{n-2}
\stackrel{d}{\to} \cdots
\stackrel{d}{\to}  E_{-n+2} 
\stackrel{\left(\begin{smallmatrix}d\\
    0\end{smallmatrix}\right)}{\longrightarrow} E_{-n+1} \oplus E_{-n} 
\stackrel{\left(\begin{smallmatrix}d_{-n+1}&
    1\end{smallmatrix}\right)}{\twoheadrightarrow} E_{-n}.
$$
The complex $\tilde{E}$ is equipped with a non-singular quadratic form
$(\xi,\tilde{\ffi})$ where  $\tilde{\ffi}=\ffi_i$ for all $i$ except in degrees $i=-n+1,
n-1$ where $\tilde{\ffi}_{n-1} = \ffi_{n-1} \oplus 1$ and $\tilde{\ffi}_{-n+1} = \ffi_{-n+1} \oplus \can_{E_{-n}}$.
Note that $(\tilde{E},\xi,\tilde{\ffi})$ is indeed non-singular since the cone of $\tilde{\ffi}$ is homotopy equivalent to the cone of $\ffi$ which is acyclic.
The bicartesian square
$$\xymatrix{
 \text{\small $(E_n \to E_{n-1} \oplus E_{-n}^{\sharp} \to \cdots E_{-n+1} \to E_{-n}) $}
  \xymono[r] \xyepi[d] & \tilde{E} \xyepi[d]\\
 E \xymono[r] & 
  \text{\small $(E_n \to E_{n-1} \to \cdots E_{-n+1} \oplus E_{-n} \to E_{-n}) $}
 }$$
 carries an obvious symmetric form induced by $\ffi$.
 On the upper left corner, this form is $\ffi$ in each degree, except in degree $n-1$ where it is $\ffi_{n-1} \oplus 1$.
 The cone of that chain map is homotopy equivalent to the cone of $\ffi$ and thus acyclic.
 In particular the symmetric form on the square is non-singular. 
 In the square, the quadratic forms $(\xi,\tilde{\ffi})$ and $(\xi,\ffi)$ on $\tilde{E}$ and $E$ restrict to the same form on the upper left corner.
 Therefore,
 $$[E,\xi,\ffi] = [\tilde{E},\xi,\tilde{\ffi}] -[H(E_{-n})[-n]] \in GW_0(\Ch_b\E,\quis,Q).$$
 
Similarly, define the complex $E'$ supported in $[-n+1, n-1]$ by
$$0 \to \coker {\left(\begin{smallmatrix}d_{n}\\
    \ffi_{n}\end{smallmatrix}\right)}
\stackrel{\left(\begin{smallmatrix}d &
    0 \end{smallmatrix}\right)}{\longrightarrow} E_{n-2}
\stackrel{d}{\to} \cdots
\stackrel{d}{\to}  E_{-n+2} 
\stackrel{\left(\begin{smallmatrix}d\\
    0\end{smallmatrix}\right)}{\longrightarrow}
\ker{\left(\begin{smallmatrix}d_{-n+1}&
    1\end{smallmatrix}\right)} \to 0.
$$
We have a diagram of complexes 
$$E' \stackrel{\sim}{\longtwoheadleftarrow} \left(
\text{\footnotesize $ E_n  \stackrel{\left(\begin{smallmatrix}d_{n}\\
    \ffi_{n}\end{smallmatrix}\right)}{\longrightarrow}  E_{n-1} \oplus E_{-n}^{\sharp} \to \cdots  \to E_{-n+2} 
\stackrel{\left(\begin{smallmatrix}d\\
    0\end{smallmatrix}\right)}{\longrightarrow}
\ker{\left(\begin{smallmatrix}d_{-n+1}&
    1\end{smallmatrix}\right)} \to 0$ }\right)
    \stackrel{\sim}{\longrightarrowtail}
    \tilde{E}$$
where the left arrow is the canonical quotient map and the right arrow is the canonical inclusion.
Both arrows are quasi-isomorphisms since kernel and cokernel are acyclic.
Moreover, the quadratic form 
$(\xi,\tilde{\ffi})$ on $\tilde{E}$ restricted to the middle complex  uniquely descends to a quadratic form $(\xi',\ffi')$ on $E'$, by quadratic left exactness of $Q$.
It follows that 
$$[E,\xi,\ffi] =  [\tilde{E},\xi,\tilde{\ffi}] -[H(E_{-n})] =  [E',\xi',\ffi'] -[H(E_{-n})] \in GW_0(\Ch_b\E,\quis,Q)$$
The same argument works for $n=1$ provided we define $\tilde{E}$ concentrated in $[-1,1]$ and $E'$ concentrated in degree $0$ as 
$$
 E_1 \stackrel{\left(\begin{smallmatrix}d_1 \\ 0 \\ \ffi_1 \end{smallmatrix}\right)} {\longrightarrowtail}  E_0 \oplus E_{-1} \oplus E_{-1}^{\sharp} \stackrel{(d_0 \ 1\ 0)}{\longtwoheadrightarrow} E_{-1} \hspace{3ex} \text{and}\hspace{3ex}
  H_0(\tilde{E}). $$
The chain complex $\tilde{E}$ is equipped with the quadratic form $(\xi \perp h_{E_{-1}}, \tilde{\ffi})$ where $h_{E_{-1}}$ is the hyperbolic quadratic form on $E_{-1}$,  $\tilde{\ffi}_1=\ffi_1$, 
 $\tilde{\ffi}_0=\ffi_0 \perp \rho(h_{E_{-1}})$ and $\tilde{\ffi}_{-1}=\ffi_{-1}$.
The quadratic form on $E'$ is the unique form which agrees with the form on $\tilde{E}$ when restricted to the middle complex in the diagram of quasi-isomorphisms
$$E' \stackrel{\sim}{\longtwoheadleftarrow} \left( E_1 \rightarrowtail \ker(d_0 \ 1\ 0) \to 0 \right) \stackrel{\sim}{\longrightarrowtail} \tilde{E}.$$

To finish the proof, by induction, any generator $[{E},\xi,{\ffi}]$ is a sum of hyperbolic spaces and a space $[E',\xi',\ffi']$ with $E'$ concentrated in degree $0$. 
Then the quasi-isomorphism $\ffi'$ is an isomorphism, and $[E,\xi,\ffi]$ is in the image of $GW_0(\E,i,Q) \to GW_0(\Ch_b\E,\quis,Q)$.
\end{proof} 

The following is the main result of this paper.
In the case of symmetric forms it was proved in 
\cite[Proposition 6]{myMV}.

\begin{theorem}
\label{prop:gilletWald}
For an exact form category with strong duality $(\E,{\sharp},\can,Q)$, 
the form functor (\ref{eqn:GilletWald})
induces a homotopy equivalence of Grothendieck-Witt spaces
$$GW(\E,i,Q) \stackrel{\sim}{\longrightarrow}
GW(\Ch_b\E,\quis,Q).$$ 
\end{theorem}

\begin{proof}
By Lemma \ref{lem:strongCof} below, we can assume that $\E$ is semi-idempotent complete and that therefore acyclic complexes in $\Ch_b\E$ are strictly acyclic.
Let $\Ac_b\E\subset \Ch_b\E$ be the full exact form subcategory of acyclic complexes, and
for $n\geq 0$, let $\Ac_{[-n,n]}\E \subset \Ac_b\E$ and $ \Ch_{[-n,n]}\E \subset \Ch_b\E$ be the full exact form subcategories of complexes $E$ satisfying $E_i=0$ for $i>n$ and $i<-n$.
For $n\geq 1$, consider the commutative diagram of exact form categories with weak equivalences the isomorphisms
\begin{equation}
\label{eqn:Agreement1}
\xymatrix{
\Ac_{[-n+1,n-1]}\E \ar[r] \ar[d] & \Ac_{[-n,n]}\E \ar[r] \ar[d] & \H \E \ar@{=}[d] \\
 \Ch_{[-n+1,n-1]}\E \ar[r] & \Ch_{[-n,n]}\E \ar[r]  & \H \E 
 }
 \end{equation}
 where the form functors in the left square are the obvious inclusions and the right horizontal form functors send a complex $E$ to $(E_n,E_{-n}^{\sharp})$ with duality compatibility $(1,\can_{E_n}): (E_{-n}^{\sharp}, E_n^{\sharp\sharp}) \to (E_{-n}^{\sharp},E_n)$.
 The rows of diagram (\ref{eqn:Agreement1}) induce (split) homotopy fibrations of $GW$-spaces in view of the Additivity Theorem \ref{thm:AddtyForWeakExFunctors}
 applied to the bicartesian squares of form functors
 sending $A \in \Ac_{[-n,n]}\E$ and $B\in \Ch_{[-n,n]}\E$ to
 $$
\xymatrix{
(A_n \to A_{n-1}\to  \cdots \to \ker d_{-n+1} \to 0) \xymono[r] \xyepi[d] &  
(A_n \to A_{n-1} \to \cdots \to A_{-n+1} \to A_{-n})  \xyepi[d] \\
(0 \to  \im d_{n-1}\to \cdots \to \ker d_{-n+1} \to 0) \xymono[r] & (0 \to  \im d_{n-1} \to \cdots A_{n-1}  \to A_n)
}$$
and
$$
\xymatrix{
(0  \to B_{n-1}\to  \cdots \to B_{-n+1} \to B_{-n}) \xymono[r] \xyepi[d] &  
(B_n \to B_{n-1} \to \cdots \to B_{-n+1} \to B_{-n})  \xyepi[d] \\
(0 \to  B_{n-1}\to \cdots \to B_{-n+1} \to 0) \xymono[r] & (B_n \to  B_{n-1} \to \cdots B_{n-1}  \to 0)
}$$
equipped with the unique forms which are the identity on the upper right corners.
Note that the condition $d_1^{\bullet}(\xi)=0$ is used to make sense of the diagram of form functors for $A \in \Ac_{[-1,1]} = S_2\E$.
It follows that the left square in (\ref{eqn:Agreement1}) is homotopy cartesian.
Then all squares in
$$
\xymatrix{
\Ac_{[0,0]}\E \ar[r] \ar[d] & \cdots \ar[r] &  \Ac_{[-n+1,n-1]}\E \ar[r] \ar[d] & \Ac_{[-n,n]}\E \ar[r] \ar[d] & \colim_n \Ac_{[-n,n]}\E \ar[d] \\
 \Ch_{[0,0]}\E \ar[r] & \cdots \ar[r] &  \Ch_{[-n+1,n-1]}\E \ar[r]  & \Ch_{[-n,n]}\E \ar[r]  & \colim_n\Ch_{[-n,n]}\E
 }
$$
including the outer square are homotopy cartesian.
In other words, the left square in
$$\xymatrix{
    0 \ar[r]\ar[d] &  (\Ac^b\E,i) \ar[d] \ar[r]& (\Ac^b\E,\quis) \ar[d]\\
   (\E,i) \ar[r]  & (\Ch^b\E,i) \ar[r] & (\Ch^b\E,\quis).}$$
induces a homotopy cartesian square of $GW$-spaces.
By Proposition \ref{prop:ChgOfWkEqIsoVariant}, the right square also induces a homotopy cartesian square of $GW$-spaces.
Hence, the total square does, too.
Since the GW-space of the upper right corner is contranctible, the homotopy fibre of the lower horizontal map is trivial, that is, the map
$GW_i(\E,i,Q) \to GW_i(\Ch_b\E,\quis,Q)$ is an isomorphism for $i>0$ and a monomorphism for $i=0$.
By Lemma \ref{lem:GilWladSurj}, the map is also surjectve for $i=0$.
This finishes the proof of the theorem.
\end{proof}

\begin{lemma}
\label{lem:strongCof}
Let $(\E,w,\sharp,\can,Q)$ be an exact form category with weak equivalences and strong
duality, then the inclusion $\E \subset \tilde{\E}^{0}$ into its semi-idempotent completion induces a
homotopy equivalence of Grothendieck-Witt spaces
$$GW(\E,w,Q) \to GW(\tilde{\E}^{0},w,Q).$$
\end{lemma}

\begin{proof}
Using Theorem \ref{thm:CofinalityExCats} and Example \ref{ex:IdempCompl}, the proof now is the same as \cite[Lemma 12]{myMV}.
\end{proof}

\section{The Grothendieck-Witt spectrum and the Bott sequence}

In this section we will generalise various results from \cite{myJPAA}.
We will construct  semi-nonconnective delooping of the $GW$-space giving rise to a spectrum $GW$ generalising \cite[\S 5]{myJPAA}.
This allows us to generalize the Bott sequence \cite[Theorem 6.1]{myJPAA}.

Let $(\E,w,\sharp,\can,Q)$ be an exact form category with weak equivalences.
Recall that the category  $\Mor\E=\Fun([1],\E)$ of morphisms in $\E$ is canonically equipped with the structure of an exact form category; see the paragraph before Proposition \ref{prop:E=Morw}.
The dual of an object $f:X \to Y$ of $\Mor\E$ is $f^{\sharp}:Y^{\sharp} \to X^{\sharp}$, the double dual identification of $f$ is $(\can_X,\can_Y)$,
the set of quadratic forms $Q(f:X\to Y)$ on $f$ is the set of pairs $(\xi,a)$ where $\xi\in Q(X)$ is a quadratic form on $X$ and $a: X \to Y^{\sharp}$ is an arrow in $\E$ such that $f^{\sharp}a= \rho(\xi)$.
A sequence in $\Mor\E$ is admissibly exact if the sequences of source objects and target objects are admissibly exact in $\E$. 
A morphism in $\Mor\E$ is a weak equivalence it the maps of source and target objects are weak equivalences in $\E$.
If $\E$ has a strong symmetric cone, then so has $\Mor\E$, by functoriality of the symmetric cone of $\E$.

Consider the exact form functor
\begin{equation}
\label{eqn:FunisH}
F: \Mor\E \to \H\E: (f:X\to Y) \mapsto (X,Y^{\sharp})
\end{equation}
with duality compatibility map $(1,\can):F(f^{\sharp})\to F(f)^{\sharp}$ and map on quadratic forms
$$Q_{\Mor\E}(f) \to Q_{\H\E}(X,Y^{\sharp})=\E(X,Y^{\sharp}): (\xi,a) \mapsto a.$$

\begin{lemma}
\label{lem:FunisH}
Let $\E=(\E,w,\sharp,\can,Q)$ be an exact form category with weak equivalences.
Then for any integer $n\geq 1$, the functor (\ref{eqn:FunisH}) induces a homotopy equivalence of pointed topological spaces
$$|w\Quad\RR_{\bullet}^{(n)}(\Mor\E,Q) | \stackrel{\sim}{\longrightarrow} |w\Quad\RR_{\bullet}^{(n)}(\H\E)|.$$
\end{lemma}

\begin{proof}
Since the lemma holds for $K$-theory (same proof as below), we are reduced to showing that (\ref{eqn:FunisH}) induces a homotopy equivalence on Grothendieck-Witt spaces, in view of Corollary \ref{cor:AddtyExFib1}.
The inverse of (\ref{eqn:FunisH}) on Grothendieck-Witt spaces is given by the form functor
\begin{equation}
\label{eqn:FunisHinverse}
(G,\psi): \H\E \to \Mor\E: (X,Y)\mapsto (0:X \to Y^{\sharp})
\end{equation}
with duality compatibility map $(\can_{Y},1_{X^{\sharp}}): G(Y,X) \to G(X,Y)^{\sharp}$
and map on quadratic forms
$$Q_{\H\E}(X,Y) = \E(X,Y) \longrightarrow Q_{\Mor\E}(X \stackrel{0}{\to} Y^{\sharp}): a \mapsto (0,a)$$
This is because $(1,\can): id_{\H\E} \to F\circ G$ is a natural weak equivalence of form functors and thus induces a homotopy of associated maps on Grothendieck-Witt spaces by Lemma \ref{lem:Lemma2MV}..
Moreover, $G\circ F$ is naturally weakly equivalent to the duality preserving  form functor $(f:X\to Y) \mapsto (0:X\to Y)$ which on quadratic forms is $Q(f:X \to Y) \to Q(0:X \to Y): (\xi,a) \mapsto (0,a)$.  Up to homotopy, this form functor induces the same map on Grothendieck-Witt spaces as the identity functor, by the Additivity Theorem, 
 in view of the natural short exact sequence in $\Mor\E$
$$(0 \to Y) \to (f:X\to Y) \to (X \to 0).$$
\end{proof}

Consider the sequence of functors
\begin{equation}
\label{eqn:FunisK}
\renewcommand\arraystretch{1.5}
\begin{array}{ccccc}
w\Quad\RR^{(n)}_{\bullet}(\Mor\E,Q) & \stackrel{(E,\xi)\mapsto E}{\longrightarrow} & w\RR^{(n)}_{\bullet}\Mor\E& &\\
& \longrightarrow & w\RR^{(n)}_{\bullet}\E & 
\stackrel{E\mapsto E\iota}{\longrightarrow} & wS^{(n)}_{\bullet}\E
\end{array}
\end{equation}
in which the non-labelled map is the functor $\Mor\E = \Fun([1],\E) \to \Fun([0],\E) = \E$ induced by the inclusion $[0] \to [1]: 0 \mapsto 0$ and the last functor is induced by the inclusion $\iota: [p] \subset [p]^{op}[p]$.

\begin{corollary}
\label{cor:RFunIsSdot}
For every integer $n\geq 1$ and every exact form category with weak equivalences $\E=(\E,w,\sharp,\can,Q)$, the composition of the maps in (\ref{eqn:FunisK}) induces a homotopy equivalence of pointed spaces
$$ |w\Quad \RR^{(n)}_{\bullet}(\Mor\E,Q)| \stackrel{\sim}{\longrightarrow} 
|wS^{(n)}_{\bullet}\E|.$$
\end{corollary}

\begin{proof}
The composition in (\ref{eqn:FunisK}) factors as
{\small
$$w\Quad\RR^{(n)}_{\bullet}(\Mor\E,Q)
\stackrel{F}{\longrightarrow}
w\Quad \RR^{(n)}_{\bullet}\H\E = (\H w\RR^{(n)}_{\bullet}\E)_h
\longrightarrow w\RR^{(n)}_{\bullet}\E \longrightarrow 
wS^{(n)}_{\bullet}\E$$}
where $\C_h$ denotes the category of symmetric forms in a category with duality $\C$ as defined in \cite[\S 2.1]{myMV}.
The first map induces a homotopy equivalence by Lemma \ref{lem:FunisH}, the second map, by \cite[\S 2 Lemma 3]{myMV}, and the third map, by \cite[\S 2 Lemma 1]{myMV}.
\end{proof}

Let $\E=(\E,w,\sharp,\can,Q)$ be an exact form category with weak equivalences which has a strong symmetric cone.
We denote by 
$$\E^{[1]} = (\Mor\E,w_{cone},\sharp,\can,Q)$$
 the exact form category with weak equivalences which has underlying exact form category $(\Mor\E,\sharp,\can,Q)$ and where the weak equivalences are the maps $f \to g$ in $\Mor\E$ which induce a weak equivalence $C(f) \to C(g)$ in $\E$ on cones.
Denote by $\cone:\Mor\E \to \E^{[1]}$ the change of weak-equivalence form functor
\begin{equation}
\label{eqn:ConeChofWEq}
\cone:(\Mor\E,w,\sharp,\can,Q) \to  (\Mor\E,w_{cone},\sharp,\can,Q)
\end{equation}
which is the identity on underlying categories and quadratic functors.
The functor $I:\E \to \Mor\E$ sending an object $E$ to the identity arrow $1_E$ on $E$ makes
$(\E,w,\sharp,\can,Q)$ into a fully faithful exact form subcategory of $(\Mor\E,w,\sharp,\can,Q)$ with image in the form subcategory $(\Mor\E,w,\sharp,\can,Q)^{\cone}$ of those objects $f$ of $\Mor\E$ for which $0\to C(f)$ is a weak equivalence in $\E$.

\begin{lemma}
\label{lem:EandMorCone}
Let $(\E,w,\sharp,\can,Q)$ be an exact form category with weak equivalences which has a strong symmetric cone. 
Then the fully faithful inclusion of form categories
$$I: (\E,w,\sharp,\can,Q) \longrightarrow (\Mor\E,w,\sharp,\can,Q)^{\cone}$$
induces homotopy equivalences of $GW$-spaces and  for $n\geq 1$ after applying the functor $|w\Quad\RR^{(n)}_{\bullet}(\phantom{E})|$.
\end{lemma}

\begin{proof}
On Grothendieck-Witt spaces, the form functor
$$F: (\Mor\E)^{\cone} \to \E: (f:X\to Y) \mapsto X$$
with duality compatibility
$f^{\sharp}:F(\sharp f)=Y^{\sharp} \to \sharp F(f) = X^{\sharp}$ which on abelian groups of quadratic forms is the map
$$Q_{\Mor\E}(f) \to Q_{\E}(F(f)): (\xi,a) \mapsto \xi$$
defines an inverse up to homotopy of $I$ in view of the 
identity $F\circ I= id$ and the
natural weak equivalences of form functors $I\circ F \to id$ 
 given by $(1_X,f)$ for $(f:Y\to X) \in \Mor\E$; see Lemma \ref{lem:WeakEqVsAcyclicCone}.
Similarly, $I$ induces an equivalence of $K$-theory spaces, and hence of $|w\Quad\RR_{\bullet}^{(n)}(\phantom{\A})|$ spaces in view of Corollary \ref{cor:AddtyExFib1}.
\end{proof}

By functoriality, if $\E = (\E,w,\sharp,\can,Q)$ has a strong symmetric cone then $\E^{[1]}$ also has a strong symmetric cone.
Setting $\E^{[0]}=\E$,  we can define recursively for $r\geq 0$ the exact form categories with weak equivalences
$$\E^{[r+1]} = (\E^{[r]})^{[1]}.$$

Let $\A=(\A,w,\sharp,\can,Q)$ be an exact form category with weak equivalences and strong symmetric cone.
Then we have the commutative diagram
\begin{equation}
\label{eqn:GWnLoopGWA}
\xymatrix{
w\Quad\left(\A,Q\right)
\ar[rr] \ar[d] & & 0 \ar[d]\\
w\Quad \RR_{\bullet}\left(\A,Q \right)\hspace{1ex}
\ar@{^(->}[r]^{\hspace{-3ex}I} \ar[d] &
w\Quad\RR_{\bullet}  \left(\Mor \A,Q\right)
\ar[d]^{\Cone} \ar[r]^{\hspace{6ex}\sim}_{\hspace{8ex}\text{(\ref{eqn:FunisK})}} &
wS_{\bullet}\A\\
w\Quad \RR_{\bullet} \left((\A^{[1]})^w,Q\right)\hspace{1ex}
\ar@{^(->}[r] &
w\Quad \RR_{\bullet} \left(\A^{[1]},Q\right).
&
}
\end{equation}

Let $\E=(\E,w,\sharp,\can,Q)$ be an exact form category with weak equivalences and strong symmetric cone.
For $n\geq 0$, replacing $\A$ above with 
$$\RR^{(n)}_{\bullet}\E^{[n]}=\RR^{(n)}_{\bullet}\left(\E^{[n]},w,\sharp,\can,Q\right)$$
 we obtain the commutative diagram of (multi-) simplicial categories
{\small
\begin{equation}
\label{eqn:GWnLoopGWn+1}
\xymatrix{
w\Quad\RR^{(n)}_{\bullet}\left(\E^{[n]},Q\right)
\ar[rr] \ar[d] & & 0 \ar[d]\\
w\Quad \RR^{(1+n)}_{\bullet}\left(\E^{[n]},Q \right)\hspace{1ex}
\ar@{^(->}[r]^{\hspace{-3ex}I} \ar[d] &
w\Quad \RR^{(1+n)}_{\bullet}\left(\Mor \E^{[n]},Q\right)
\ar[d]^{\Cone} \ar[r]^{\hspace{6ex}\sim}_{\hspace{8ex}\text{(\ref{eqn:FunisK})}} &
wS^{(1+n)}_{\bullet}\E^{[n]}\\
w\Quad \RR^{(1+n)}_{\bullet}\left((\E^{[1+n]})^w,Q\right)\hspace{1ex}
\ar@{^(->}[r] &
w\Quad \RR^{(1+n)}_{\bullet}\left(\E^{[1+n]},Q\right).
&
}
\end{equation}
}

\begin{proposition}
\label{prop:GWtoGWtildeDiagram}
Let $\E=(\E,w,\sharp,\can,Q)$ be an exact form category with weak equivalences and strong symmetric cone.
Then, after topological realization, in diagram (\ref{eqn:GWnLoopGWn+1}), 
the lower square is homotopy cartesian
with contractible lower left corner for $n\geq 0$, and
the upper square 
is homotopy cartesian
for $n\geq 1$.
\end{proposition}

\begin{proof}
For $n\geq 1$, the upper square is homotopy cartesian, in view of Corollary \ref{cor:AddtyExFib1}.
The lower square is homotopy cartesian, by Corollary \ref{cor:RndotFib} applied to the change of weak equivalence form functor (\ref{eqn:ConeChofWEq}) where we use Lemma \ref{lem:EandMorCone} to identify the upper left corner of the square.
\end{proof}

Diagram (\ref{eqn:GWnLoopGWn+1}) defines canonical maps for $n\geq 0$
\begin{equation}
\label{eqn:GWbondingMaps}
\left|w\Quad\RR^{(n)}_{\bullet}\left(\E^{[n]},Q\right)\right| \longrightarrow 
\Omega \left|w\Quad\RR^{(1+n)}_{\bullet}\left(\E^{[1+n]},Q\right)\right|.
\end{equation}
Actually, the map is a functorial zigzag of maps where the arrows in the wrong direction are homotopy equivalences.

\begin{definition}[Grothendieck-Witt spectrum]
\label{dfn:GWSpectrum}
Let $\E=(\E,w,\sharp,\can,Q)$ be an exact form category with weak equivalences and strong symmetric cone.
The {\em Grothendieck-Witt spetrum of $\E$} is the sequence of spaces
$$GW(\E,w,Q) = \left\{GW(\E,w,Q)_0, \ GW(\E,w,Q)_1,\ GW(\E,w,Q)_2,\ \dots \right\}
$$
where
$$\ GW(\E,w,Q)_n =  |w\Quad\RR^{(n)}_{\bullet}\left(\E^{[n]},Q\right)|$$
together with the bonding maps defined in  (\ref{eqn:GWbondingMaps}).
For $n\geq 0$ we set
$$GW^{[n]}(\E,w,Q)=GW(\E^{[n]},w,Q).$$
As usual, we define the $n$-th shifted Grothendieck-Witt groups as the homotopy groups of the Grothendiek-Witt spectrum
$$GW^{[n]}_i(\E,w,Q) = \pi_iGW^{[n]}(\E,w,Q)$$
for $i\in \Z$.
By Theorem \ref{thm:OmegaInfOfGW} below, these groups coincide with the homotopy groups of the Grothendieck-Witt space for $i\geq 0$.
\end{definition}

\begin{definition}
Let $\E=(\E,\sharp,\can,Q)$ be an exact form category with strong duality (equipped with isomorphisms as weak equivalences), we set
$$GW^{[n]}(\E,Q) = GW^{[n]}(\Ch_b\E,\quis,\sharp,\can,Q)$$
where the form category with weak equivalences $(\Ch_b\E,\quis,\sharp,\can,Q)$ was defined in \S \ref{sec:FromExToChain}.
By Theorem \ref{thm:OmegaInfOfGW} below and Theorem \ref{prop:gilletWald}, its infinite loop space is the Grothendieck-Witt space of \cite[Definition 6.3]{myForm1} which is the group completion of $|w\Quad \left(\E,Q\right)|$ in case all admissible exact sequences in $\E$ split 
\cite[Theorem 6.5]{myForm1}.
\end{definition}

\begin{theorem}
\label{thm:OmegaInfOfGW}
Let $\E=(\E,w,\sharp,\can,Q)$ be an exact form category with weak equivalences and strong symmetric cone.
Then the maps in (\ref{eqn:GWbondingMaps}) are homotopy equivalences for $n\geq 1$.
For $n=0$, the map in (\ref{eqn:GWbondingMaps}) factors as 
$$|w\Quad\left(\E,Q\right)| \longrightarrow GW(\E,w,Q) \stackrel{\simeq}{\longrightarrow} \Omega |w\Quad\RR_{\bullet}\left(\E^{[1]},Q\right)|$$
where the first map is the universal map (\ref{eqn:EGWunivMap}) and the second is a homotopy equivalence.
In particular, the spectrum $GW(\E,w,Q)$ of Definition \ref{dfn:GWSpectrum} is a positive $\Omega$-spectrum with associated infinite loop space the Grothendieck-Witt space $GW(\E,w,Q)$ of Definition \ref{dfn:GWSpaceII}.
\end{theorem}

\begin{proof}
For $n\geq 1$ this is Proposition \ref{prop:GWtoGWtildeDiagram}, and for $n=0$, this is Proposition \ref{prop:GWtoGWtildeDiagram} together with the definition of the Grothendieck-Witt space.
\end{proof}

\begin{lemma}
Let $\E = (\E,w,\sharp,\can,Q)$ be an exact form category with weak equivalences equipped with a strong symmetric cone $(C,P)$.
Then the exact functor functor
\begin{equation}
\label{eqn:EandE[1]inKthCone}
(\Mor\E,w_{cone} )\to  (\E,w): f \mapsto C(f)
\end{equation}
induces an equivalence of $K$-theory spaces with inverse given by 
\begin{equation}
\label{eqn:EandE[1]inKth}
(\E,w) \to (\Mor\E,w_{cone}):X \mapsto (0 \to X).
\end{equation}
Moreover, the map (\ref{eqn:EandE[1]inKth}) induces the negative of the map
\begin{equation}
\label{eqn:EandE[1]inKthNegative}
(\E,w) \to (\Mor\E,w_{cone}):X \mapsto (X \to 0)
\end{equation}
in $K$-theory.
\end{lemma}

\begin{proof}
The composition $\E \to \Mor\E \to \E$ is the identity, and the composition $\Mor\E \to \E \to \Mor\E$ is naturally weakly equivalent to the identity by the following zigzag of $w_{cone}$-weak equivalences
$$\xymatrix{
X \xymono[r]^{i_X} \ar[d]_f & C(X) \ar[d] & 0 \ar[d] \ar[l]\\
Y \xymono[r] & C(f) & C(f). \ar[l]_1}$$
The last claim follows from Additivity in $K$-theory applied to the admissible exact sequence of functors $(\E,w) \to (\Mor\E,w_{cone})$
$$\xymatrix{
0 \xymono[r] \ar[d] & X \xyepi[r]^1 \ar[d]^1 & X \ar[d] \\
X \xymono[r]_1 & X \ar[r] & 0}
$$
since the middle functor lands in the subcategory of $w_{cone}$-acyclic objects.
\end{proof}

It follows that $\E^{[1]}$ and $\E$ have the same $K$-theory space, and by iteration, so do $\E^{[n]}$ and $\E$.
In particular, for every exact form category with weak equivalences and strong symmetric cone we have the forgetful map
$$F:GW^{[n]}(\E,w,Q) \stackrel{\text{Dfn \ref{dfn:HypCatFunEtc}}}{\longrightarrow} K(\E^{[n]},w) \stackrel{\text{(\ref{eqn:EandE[1]inKthCone})}}{\longrightarrow} K(\E,w)$$
and the hyperbolic map
$$H:K(\E,w) \stackrel{\text{(\ref{eqn:EandE[1]inKth})}}{\longrightarrow} K(\E^{[n+1]},w) \stackrel{\text{Dfn \ref{dfn:HypCatFunEtc}}}{\longrightarrow} GW^{[n+1]}(\E,w,Q).$$

Consider the commutative square of exact form categories with weak equivalences.
\begin{equation}
\label{eqn:BottSquare}
\xymatrix{
(\E,w,Q) \ar[r]^I \ar[d] & (\Mor\E,w,Q) \ar[d]\\
(\Mor\E,w_{cone},Q)^{w_{cone}} \ar[r] & (\Mor\E,w_{cone},Q)
}
\end{equation}
The following generalizes \cite[Theorem 6.1]{myJPAA}.

\begin{theorem}[Algebraic Bott sequence]
\label{thm:BottSeq}
Let $\E = (\E,w,\sharp,\can,Q)$ be an exact form category with weak equivalences equipped with a strong symmetric cone.
Then diagram (\ref{eqn:BottSquare}) induces a homotopy fibration of $GW$-spectra for all $n\geq 0$
$$GW^{[n]}(\E,w,Q) \stackrel{F}{\longrightarrow} K(\E,w) \stackrel{H}{\longrightarrow} GW^{[n+1]}(\E,w,Q).$$
\end{theorem}

\begin{proof}
By Corollary  \ref{cor:RndotFib} for $\E^{[n]}$ in place of $\E$ and $\RR^{(n)}$,
diagram (\ref{eqn:BottSquare}) induces a homotopy cartesian square of $GW$-spectra with contractible lower left corner.
In particular, we have a homotopy fibration of spectra
\begin{equation}
\label{eqn:PreGWspBottFib}
GW(\E,w,Q) \to GW(\Mor\E,w,Q) \to GW^{[1]}(\E,w,Q).
\end{equation}
By Lemma \ref{lem:FunisH}, the form functor $(\Mor\E,w,Q) \to \H\E$ of (\ref{eqn:FunisH}) induces an equivalence of $GW$-spectra with inverse the form functor (\ref{eqn:FunisHinverse}).
One observes that the composition
$$(\E,w,Q) \stackrel{I}{\longrightarrow} (\Mor\E,w,Q) \stackrel{\text{ (\ref{eqn:FunisH})}}{\longrightarrow} (\H\E,w)$$
is the forgetful functor of Definition \ref{dfn:HypCatFunEtc}, and the composition
$$(\H\E,w) \stackrel{\text{ (\ref{eqn:FunisHinverse})}}{\longrightarrow}  (\Mor\E,w,Q) \longrightarrow (\Mor\E,w_{cone},Q)$$
is the composition
$$(\H\E,w) \stackrel{\text(\ref{eqn:EandE[1]inKthNegative})}{\longrightarrow} \H(\Mor\E,w_{cone}) \stackrel{H}{\longrightarrow} (\Mor\E,w_{cone},Q).$$
Thus, the homotopy fibration of spectra (\ref{eqn:PreGWspBottFib}) yields the homotopy fibration
$$GW(\E,w,Q) \stackrel{F}{\longrightarrow} K(\E,w) \stackrel{-H}{\longrightarrow} GW^{[1]}(\E,w,Q).$$
Replacing $H$ with $-H$ yields the theorem for $n=0$.
Replacing $\E$ with $\E^{[n]}$ yields the general case.
\end{proof}

\bibliographystyle{plain}

\newcommand{\etalchar}[1]{$^{#1}$}

\end{document}